\documentstyle[12pt,A4]{article}
\catcode `\@=11
\@addtoreset{equation}{section}

\def\harr#1#2{\smash{\mathop{\hbox to .5in{\rightarrowfill}}
\limits^{\scriptstyle#1}_{\scriptstyle#2}}}

\def\harrl#1#2{\smash{\mathop{\hbox to .5in{\leftarrowfill}}
\limits^{\scriptstyle#1}_{\scriptstyle#2}}}

\def\qed{\vrule height 6pt width 6pt depth 6pt}

\newcommand{\be}{\begin{equation}}
\newcommand{\ee}{\end{equation}}
\newcommand{\R}{\rm I \mkern -3mu R}
\newcommand{\N}{\rm I \mkern -3mu N}
\newtheorem{thm}{Theorem}[section]
\newtheorem{rem}{Remark}[thm]
\newtheorem{lemma}[thm]{Lemma}
\newtheorem{cor}[thm]{Corollary}
\newtheorem{prop}[thm]{Proposition}
\begin{document}
\begin{titlepage}
\begin{center}
{\bf \Large{The Variational Sequence on Finite Jet Bundle 
~\\Extensions and the Lagrangian Formalism}}
\end{center}
\vskip 1.0truecm
\centerline{D. R. Grigore
\footnote{e-mail: grigore@theor1.ifa.ro, grigore@roifa.ifa.ro}}
\vskip5mm
\centerline{Dept. of Theor. Phys., Inst. Atomic Phys.}
\centerline{Bucharest-M\u agurele, P. O. Box MG 6, ROM\^ANIA}
\vskip 2cm
\bigskip \nopagebreak
\begin{abstract}
\noindent
The geometric Lagrangian theory (of arbitrary order) is based on the analysis
of some basic mathematical objects such as: the contact ideal, the (exact)
variational sequence, the existence of Euler-Lagrange and Helmholtz-Sonin
forms, etc. In this paper we give new and much simpler proofs for the whole
theory using Fock space methods. Using these results we give the most general
expression for a variationally trivial Lagrangian and the generic expression
for a locally variational system of partial differential equation.
\end{abstract}
\end{titlepage}

\section{Introduction}

The modern Lagrangian theory is formulated in the language of jet bundle 
extensions. There are two different, but closely related approaches in this
framework. One is based on infinite jet bundle extensions \cite{AD}, \cite{An},
\cite{Tu} (and uses the so-called variational bicomplex) and the other is based
on finite jet bundle extensions \cite{Kr1} - \cite{Kr7} (and uses the so-called
variational sequence). The second approach is more in the spirit of the 
original Lagrangian theory as studied in physical literature corresponding 
to a Lagrangian theory of arbitrary, but finite order. On the other hand, the
study of the finite order seems to be much more complicated from the
combinatorial point of view as it is apparent from the proofs appearing in the
references quoted above. In some recent papers we have tried to prove that the
simplest way to disentangle the combinatorial proofs is to use the observation
that the various tensors appearing in the line of argument have some symmetry
or antisymmetry properties which make them as elements in some Fock space. Then
most of the relations to be solved can be written with the help of the creation
and annihilation operators. Using elementary properties of these Fock space
operators one can significantly simplify the proofs. This is one of the main
motivation of writing this paper, i.e to clarify the technical aspects. It
appears that it is convenient to write the paper in a self contained style so
it will also serve as a pedagogical introduction to this field of interest.
In Section 2 we provide the basic construction - the (finite) jet bundle
extension idea and give some relevant formul\ae. In Sections 3 and 4 we study 
the contact and respectively the strong contact forms on a jet bundle 
extension and obtain with our method the structure of the contact ideal. In 
section 5 we establish the existence of the Euler-Lagrange and of 
Helmholtz-Sonin forms adapting to the finite jet bundle extension approach 
some ideas of Anderson and collaborators. The proof of the existence of the 
Helmholtz-Sonin form is new. In Section 6 we sketch the proof 
for the exactness of the variational sequence and we clarify some points about
the characterization of the elements of the variational sequence by forms. 
In Section 7 we give the most general expression for a variationally trivial
Lagrangian and in Section 8 we provide the generic expression for a locally
variational partial differential system of equations using the methods
developed in the preceding sections. These problems have been under intensive
study for sometime in the literature but many of the results obtained here are 
new. In particular, we obtain in Section 8 an {\it if and only if} type of
result. Section 9 is devoted to some final comments.

\newpage
\section{Finite Order Jet Bundle Extensions}

The content of this section is standard and is included mainly to fix the
notations. For the sake of completeness, we give however some proofs which are
not trivial.

The kinematical structure of a classical field theory is based on fibred
bundle structures. Let
$
\pi: Y \mapsto X
$
be fibre bundle, where
$X$
and
$Y$
are differentiable manifolds of dimensions 
$
dim(X) = n,\quad dim(Y) = m + n
$
and
$\pi$
is the canonical projection of the fibration. Usually $X$
is interpreted as the ``space-time" manifold and the fibres of $Y$
as the field variables. An {\it adapted chart} to the fibre bundle structure
is a couple
$(V,\psi)$
where $V$ is an open subset of $Y$ and 
$\psi: V \rightarrow \R^{n} \times \R^{m}$
is the so-called {\it chart map}, usually written as
$
\psi = (x^{i},y^{\sigma}) \quad (i = 1,...,n; \sigma = 1,...,m)
$
such that
$(\pi(V), \phi)$
where
$
\phi = (x^{i}) \quad (i = 1,...,n)
$
is a chart on $X$ and the canonical projection has the following expression:
$
\pi(x^{i},y^{\sigma}) = (x^{i}).
$
If 
$p \in Y$
then the real numbers
$x^{i}(p), \quad y^{\sigma}(p)$
are called the {\it (fibred) coordinates} of $p$. For simplicity we will give 
up the attribute {\it adapted} sometimes in the following. Also we will refer
frequently to the first entry $V$ of $(V,\psi)$ as a chart.

Next, one considers the $r$-{\it jet bundle extensions}
$
J^{r}_{n}Y \mapsto X \quad (r \in \N).
$
The construction is the following (see for instance \cite{Kr3}). 

\begin{thm}
Let 
$x \in X$,
and
$y \in \pi^{-1}(x)$.
We denote by
$\Gamma_{(x,y)}$
the set of sections
$
\gamma: U \rightarrow Y
$
such that: (i) $U$ is a neighbourhood of $x$; (ii) 
$\gamma(x) = y$.
We define on 
$\Gamma_{(x,y)}$
the relationship
``$\gamma \sim \delta$"
{\it iff} there exists a chart
$(V,\psi)$
on $Y$ such that $\gamma$ and $\delta$ have the same partial derivatives up to
order $r$ in the given chart i.e.
\be
{\partial^{k} \over \partial x^{i_{1}}...\partial x^{i_{k}}} 
\psi \circ \gamma \circ \phi^{-1} (\phi(x)) =
{\partial^{k} \over \partial x^{i_{1}}...\partial x^{i_{k}}} 
\psi \circ \delta \circ \phi^{-1} (\phi(x)), 
\quad k \leq r.
\ee

Then this relationship is chart independent and it is an equivalence relation.
\label{equivalence}
\end{thm}

{\bf Proof:}
(i) We first prove the chart independence. Let
$
(V,\psi), (\bar{V},\bar{\psi})
$
be a two (adapted) chart on this bundle such that
$
V \cap \bar{V} \not= \emptyset
$
and let
$
(\pi(V),\phi), (\pi(\bar{V}),\bar{\phi})
$
the corresponding charts on 
$X$.
If we define
$
f: \R^{n} \rightarrow \R^{n}
$
by:
$$
f \equiv \bar{\phi} \circ \phi^{-1}
$$
then the function
$
\bar{\psi} \circ \psi^{-1}: \R^{m+n} \rightarrow \R^{m+n}
$
has the following expression:
\be
\bar{\psi} \circ \psi^{-1} = (f^{i},F^{\sigma})
\ee
where 
$F$ 
is a smooth function from 
$\R^{n}$ 
into 
$\R^{m+n}$.
Let
$\gamma \in \Gamma_{(x,y)}$; 
then in adapted coordinates we have the following expression:
\be
\psi \circ \gamma \circ \phi^{-1}(x^{i}) = (x^{i},g^{\sigma}(x)).
\ee

Similarly we have in the other chart:
$$
\bar{\psi} \circ \gamma \circ \bar{\phi}^{-1}(x^{i}) = 
(x^{i},\bar{g}^{\sigma}(x))
$$
where:
\be
\bar{g}^{\sigma}(x) = F^{\sigma}(f^{-1}(x),g^{\sigma}(f^{-1}(x))).
\ee

We must show that if
$g^{\sigma}$
and
$g'^{\sigma}$
have the same partial derivatives up to order $r$ then
$\bar{g}^{\sigma}$
and
$\bar{g'}^{\sigma}$
have the same property. 
To prove this one firstly gets from the expression of
$\bar{g}^{\sigma}$
above that:
\be
{\partial \bar{g}^{\sigma} \over \partial x^{i}} = 
(\Delta_{i} F^{\sigma})(f^{-1}(x),g^{\sigma}(f^{-1}(x)))
\label{Delta}
\ee
where the differential operator
$\Delta_{i}$
acts in the space of smooth functions 
$f: \R^{m+n} \rightarrow \R$
according to:
\be
\Delta_{i} \equiv P^{j}_{i} \left({\partial \over \partial x^{j}} +
{\partial g^{\sigma} \over \partial x^{j}} 
{\partial \over \partial y^{\sigma}} \right)
\ee
and
$P^{j}_{i}$
is the inverse of the matrix
${\partial f^{l} \over \partial x^{p}}$:
\be
P^{j}_{i} {\partial f^{i} \over \partial x^{p}} = \delta^{j}_{p}.
\label{x-bar-x}
\ee

From (\ref{Delta}) we obtain by recurrence:
\be
{\partial^{k} \bar{g}^{\sigma} \over \partial x^{i_{1}}...
\partial x^{i_{k}}} = 
(\Delta_{i_{1}}\cdots \Delta_{i_{k}} F^{\sigma})
(f^{-1}(x),g^{\sigma}(f^{-1}(x))).
\label{bar-g}
\ee

But it is clear that the expression
$\Delta_{i_{1}}\cdots \Delta_{i_{k}} F^{\sigma}$
contains the derivatives of the function
$g^{\sigma}$
only up to order $r$. So, the chart independence of the relation
$\sim$
follows.

(ii) Using the chart independence proved above, it easily follows that
$\sim$
is an equivalence relation. 
$\qed$

A $r$-{\it order jet with source} $x$ {\it and target} $y$
is, by definition, the equivalence class of some section
$\gamma$
with respect to the equivalence relationship defined above and it is 
denoted by
$
j^{r}_{x}\gamma.
$

Let us define
$
J^{r}_{(x.y)}\pi \equiv \Gamma_{(x,y)} / \sim
$
Then the
$r$-{\it order jet bundle extension} is, set theoretically
$
J^{r}Y \equiv \bigcup_{x} J^{r}_{(x,y)}\pi.
$

Now we define the following projections:
$
\pi^{r,s}: J^{r}Y \rightarrow J^{s}Y \quad (0 < s \leq r)
$
by
\be
\pi^{r,s}(j^{r}_{x}\gamma) = j^{s}_{x}\gamma,
\ee
$
\pi^{r,0}: J^{r}Y \rightarrow Y
$
given by
\be
\pi^{r,0}(j^{r}_{x}\gamma) = \gamma(x).
\ee
(this is consistent with the identification
$
J^{0}Y \equiv Y
$)
and finally
$
\pi^{r}: J^{r}Y \rightarrow X
$
by
\be
\pi^{r}(j^{r}_{x}\gamma) = x.
\ee
These projections are obviously surjective.

Let
$(V,\psi), \quad \psi = (x^{i},y^{\sigma})$
be a chart on $Y$. Then we define the couple
$(V^{r},\psi^{r})$, 
where:
$V^{r} = (\pi^{r,0})^{-1}(V)$
and
$$
\psi = (x^{i},y^{\sigma},y^{\sigma}_{j},...,y^{\sigma}_{j_{1},...,j_{k}},...,
y^{\sigma}_{j_{1},...,j_{r}}), \quad j_{1} \leq j_{2} \leq \cdots \leq j_{k}, 
\quad k = 1,...,r
$$
where
\begin{eqnarray}
y^{\sigma}_{j_{1},...,j_{k}}(j^{r}_{x}\gamma) = 
\left. {\partial^{k} \over \partial x^{j_{1}} \cdots \partial x^{j_{k}}} 
y^{\sigma} \circ \gamma \circ \phi^{-1} \right|_{\phi(x)}, \quad k = 1,...,r
\nonumber \\
x^{i}(j^{r}_{x}\gamma) = x^{i}(x), \quad 
y^{\sigma}(j^{r}_{x}\gamma) = y^{\sigma}(\gamma(x)).
\end{eqnarray}

Then
$(V^{r},\psi^{r})$
is a chart on
$J^{r}Y$
called {\it the associated chart} of
$(V,\psi)$.

\begin{rem}
The expressions 
$y^{\sigma}_{j_{1},...,j_{k}}(j^{r}_{x}\gamma)$
are defined for {\it all} indices
$j_{1},...,j_{k} = 1,...,n$,
and the restrictions
$j_{1} \leq j_{2} \leq \cdots \leq j_{k}$
in the definition of the charts are in order to avoid overcounting and are a
result of the obvious symmetry property:
\be
y^{\sigma}_{j_{P(1)},...,j_{P(k)}}(j^{r}_{x}\gamma) =
y^{\sigma}_{j_{1},...,j_{k}}(j^{r}_{x}\gamma), 
\label{sym}
\ee
for any permutation
$P \in {\cal P}_{k}, \quad k = 2,...,r.$
\end{rem}

Now we have the following result.

\begin{thm}

If a collection of (adapted) charts
$(V,\psi)$
are the elements of a differentiable atlas on
$Y$
then 
$(V^{r},\psi^{r})$
are the elements of a differentiable atlas on
$
J^{r}_{n}(Y)
$
which admits a fibre bundle structure over $Y$.
\label{smooth}
\end{thm}

{\bf Proof:} 

As in theorem \ref{equivalence} we consider two non-overlapping charts
$
(V,\psi), (\bar{V},\bar{\psi})
$
and let
$
(\pi(V),\phi)
$, 
$
(\pi(\bar{V}),\bar{\phi})
$
the corresponding charts on 
$X$.
If we define
$
f: \R^{n} \rightarrow \R^{n}
$
and
$
F: \R^{m+n} \rightarrow \R^{m}
$
as in the proof of this theorem, then we have
\be
\bar{\psi}^{r} \circ (\psi^{r})^{-1} = 
(f^{i},F^{\sigma},F^{\sigma}_{j},...,F^{\sigma}_{j_{1},...,j_{r}})
\label{transf}
\ee
where
$
F^{\sigma}_{j_{1},...,j_{k}}, \quad j_{1} \leq j_{2} \leq \cdots \leq j_{k},
\quad k = 1,...,r
$
are functions depending on the variables
$(x^{i},y^{\sigma},y^{\sigma}_{j},...,y^{\sigma}_{j_{1},...,j_{k}},...,
y^{\sigma}_{j_{1},...,j_{r}}).
$

We must prove that these functions are smooth. To this purpose we obtain an
rather explicit formula for them. Inspecting the notations from theorem
\ref{equivalence} it is clear that 
\be
y^{\sigma}_{j_{1},...,j_{k}} = 
\left. {\partial g^{\sigma} \over \partial x^{j_{1}}...\partial
x^{j_{k}}}\right|_{\phi(x)} 
\label{g-sigma}
\ee
and
$$
\bar{y}^{\sigma}_{j_{1},...,j_{k}} = 
\left. {\partial \bar{g}^{\sigma} \over \partial x^{j_{1}}...\partial
x^{j_{k}}}\right|_{\bar{\phi}(x)}. 
$$

But we have an explicit formula for the last expression, namely (\ref{bar-g}). 
By induction one proves that the expression
$\Delta_{i_{1}}\cdots \Delta_{i_{k}} F^{\sigma}$
is a polynomial in the expressions
$$
g^{\sigma}_{j_{1},...,j_{l}} = 
\left. {\partial g^{\sigma} \over \partial x^{j_{1}}...\partial
x^{j_{k}}}\right|_{\phi(x)}, \quad l = 1,...,k
$$
with coefficients depending only on
$(x^{i},y^{\sigma}).$
But we immediately have
$$
\Delta_{j_{1}} g^{\sigma}_{j_{2},...,j_{l}} = 
P^{i}_{j_{1}} g^{\sigma}_{i,j_{2},...,j_{l}}, \quad l = 1,...,k.
$$
So one can replace 
$\Delta_{i} \rightarrow \tilde{\Delta}_{i}$
where
$\tilde{\Delta}_{i}$
is an operator acting on polynomials in
$
g^{\sigma}_{j_{1},...,j_{l}}, \quad l = 1,...,k
$
with coefficients depending only on 
$(x^{i},y^{\sigma})$
according to
$$
\tilde{\Delta}_{j_{1}} g^{\sigma}_{j_{2},...,j_{l}} = 
P^{i}_{j_{1}} g^{\sigma}_{i,j_{2},...,j_{l}}, \quad l = 1,...,k.
$$
and
$$
\tilde{\Delta}_{j} f = \Delta_{j} f
$$
on the coefficients which are functions depending only on
$(x^{i},y^{\sigma})$.

Taking into account (\ref{g-sigma}) it follows that
$
F^{\sigma}_{j_{1},...,j_{k}}, \quad k = 1,...,r
$
is a polynomial in
$
y^{\sigma}_{j_{1},...,j_{l}} 
$
for
$l = 1,...,k$
with smooth coefficients of 
$(x^{i},y^{\sigma})$
so these are smooth functions. The existence of the differentiable structure is
proved. The fibre bundle structure over $X$ is obvious. 
$\qed$

To be able to use the summation convention over the dummy indices we consider
$y^{\sigma}_{j_{1},...,j_{k}}$
for {\it all} values of the indices
$j_{1},...,j_{k} \in \{ 1,...,n\}$
as smooth functions on the chart
$V^{r}$
defined in terms of the independent variables
$y^{\sigma}_{j_{1},...,j_{k}}, j_{1} \leq j_{2} \leq ...\leq j_{k}
\quad k = 1,2,...,r$ 
according to the formula (\ref{sym}) and we make a similar convention for the
partial derivatives
${\partial \over \partial y^{\sigma}_{j_{1},...,j_{k}}}$.

Then we define on the chart
$V^{r}$
the following vector fields:
\be
\partial^{j_{1},...,j_{k}}_{\sigma} \equiv {r_{1}!...r_{n}! \over k!} 
{\partial \over \partial y^{\sigma}_{j_{1},...,j_{k}}}, \quad k = 1,...,r
\label{partial}
\ee
for all values of the indices
$j_{1},...,j_{k} \in \{ 1,...,n\}$.
Here
$r_{l}, \quad l = 1,...,n$
is the number of times the index $l$ enters into the set
$\{ j_{1},...,j_{k}\}$.

One can easily verify the following formul\ae:
\be
\partial^{i_{1},...,i_{k}}_{\sigma} y^{\nu}_{j_{1},...,j_{l}} = 0, \quad
(k \not= l)
\ee
\be
\partial^{i_{1},...,i_{k}}_{\sigma} y^{\nu}_{j_{1},...,j_{k}} =
{\cal S}^{+}_{j_{1},...,j_{k}} 
\delta^{i_{1}}_{j_{1}} \cdots \delta^{i_{k}}_{j_{k}}
\ee
where
${\cal S}^{+}_{j_{1},...,j_{k}}$
is the symmetrization projector operator in the indices
$j_{1},...,j_{k}$
defined by the formula (\ref{sa-indices}) from the Appendix.

Also we have for any smooth function $f$ on the chart $V^{r}$: 
\be
df = {\partial f \over \partial x^{i}} dx^{i} + \sum_{k=0}^{r} 
(\partial^{j_{1},...,j_{k}}_{\sigma}f) dy^{\sigma}_{j_{1},...,j_{k}} = 
{\partial f \over \partial x^{i}} dx^{i} + \sum_{|J| \leq r} 
(\partial^{J}_{\sigma}f) dy^{\sigma}_{J}. 
\label{df}
\ee
In the last formula we have introduced the multi-index notations in an obvious
way. This formula also shows that the coefficients appearing in the definition 
(\ref{partial}) are exactly what is needed to use the summation convention
over the dummy indices without overcounting.

We now define the expressions
\be
d_{i}^{r} \equiv {\partial \over \partial x^{i}} + \sum_{k=0}^{r-1} 
y^{\sigma}_{i,j_{1},...,j_{k}} \partial^{j_{1},...,j_{k}}_{\sigma}
\label{formal}
\ee
called {\it formal derivatives}. When it is no danger of confusion we denote
simply 
$d_{i} = d_{i}^{r}.$

\begin{rem}
The formal derivatives are not vector fields on 
$J^{r}Y$.
\end{rem}

Next one immediately sees that
\be
d_{i} y^{\sigma}_{j_{1},...,j_{k}} = y^{\sigma}_{i,j_{1},...,j_{k}}, 
\quad k = 0,...,r-1
\ee
and one can obtain from theorem \ref{smooth} a recurrence relation for the 
chart transformation formul\ae: 

\begin{lemma}
The following recurrence relation are valid for the chart transformation
expressions (\ref{transf}):
\be
F^{\sigma}_{j_{1},...,j_{k}} = Q^{l}_{j_{1}} d_{l}
F^{\sigma}_{j_{2},...,j_{k}}, \quad k = 1,...,r-1.
\ee

In particular, the functions
$F^{\sigma}_{j_{1},...,j_{k}}$
are polynomial expressions in the variables
$y^{\sigma}_{j},...,y^{\sigma}_{j_{1},...,j_{k}}$
(with smooth coefficients depending on
$(x^{i},y^{\sigma})$) for all $k = 1,...,r$.
\label{tr}
\end{lemma}

From the definition of the formal derivatives it easily follows by direct
computation that:
\be
\left[ \partial^{j_{1},...,j_{k}}_{\sigma} , d_{i} \right] = 
{1\over k} \sum_{l=1}^{k} \delta^{j_{l}}_{i} 
\partial^{j_{1},...,\hat{j}_{l},...,j_{k}}_{\sigma}, \quad k = 0,...,r
\label{commutator}
\ee
where we use Bourbaki conventions
$\sum_{\emptyset} \equiv 0, \quad \prod_{\emptyset} \equiv 1$.

Based on this relation one gets:
\be
[d_{i},d_{j}] = 0
\ee
so one can consistently define for any multi-index
\be
d_{J} \equiv \prod_{i \in J} d_{i}.
\ee

If $\gamma$ is a section of the fibre bundle $Y$ then the map
$X \ni x \mapsto j^{r}_{x}\gamma \in J^{r}Y$
is section of the fibre bundle 
$J^{r}Y$
called the $r$-{\it extension} of $\gamma$.
We denote it by:
$j^{r}\gamma: V \rightarrow J^{r}Y$ .

Then one can define a map
$h: T(J^{r+1}Y) \rightarrow T(J^{r}Y)$
called {\it horizontalization} by
\be
h = h_{j^{r+1}_{x}\gamma} \equiv \left(j^{r}\gamma\right)_{*,x} \circ
\left(\pi^{r+1}\right)_{*,j^{r+1}_{x}\gamma}.
\ee

If
$\xi \in  T_{j^{r+1}_{x}\gamma}(J^{r+1}Y)$
then
$h\xi$
is called its {\it horizontal} component. We will also define 
\be
p\xi \equiv (\pi^{r+1,r})_{*}\xi - h\xi.
\label{hp}
\ee

Let
$\pi_{i}: Y_{i} \rightarrow X_{i}, \quad i = 1,2$
be two fibre bundles. Then a map
$\phi: Y_{1} \rightarrow Y_{2}$
is called a {\it fibre bundle morphism} if there is a map
$\phi_{0}: X_{1} \rightarrow X_{2}$
such that
$$
\pi_{2} \circ \phi = \phi_{0} \circ \pi_{1};
$$
one also says that $\phi$ {\it covers} $\phi_{0}$.

In this case, one can define the $r$-{\it order jet extension} of 
$\phi$ as the the map
$j^{r}\phi: J^{r}Y_{1} \rightarrow J^{r}Y_{2}$
given by
\be
j^{r}\phi (j^{r}_{x}\gamma) \equiv j^{r}_{\phi_{0}(x)} \phi \circ \gamma 
\circ \phi_{0}^{-1}, \quad \forall x \in X_{1}.
\ee

If $\xi$ is a projectable vector field on the fibre bundle Y i.e. there is a
vector field $\xi_{0}$ on $X$ such that
$\xi_{0} = \pi_{*} \xi$,
then the flow 
$\phi^{\xi}_{t}$
associated to $\xi$ covers the flow 
$\phi^{\xi_{0}}_{t}$
associated to $\xi_{0}$
so we can define the $r$-{\it order extension} of $\xi$ by
\be
j^{r}\xi \equiv \left. {d \over dt}\right|_{t=0} j^{r}\phi^{\xi}_{t}.
\ee

The vector field 
$j^{r}\xi$ 
is projectable.

A $\pi$-vertical vector field on the fibre bundle $Y$ is called a {\it
evolution}. In the chart
$(V,\psi)$
an evolution has the expression:
\be
\xi = \xi^{\sigma} \partial_{\sigma}
\label{ev}
\ee
with 
$\xi^{\sigma}$
smooth functions on the chart $V$. One denotes the set of evolutions by
${\cal E}(J^{r}Y)$;
this set is a fibre bundle over
$J^{r}Y$.

We denote the forms of degree $q$ on 
$J^{r}Y$
by
$\Omega^{r}_{q}$.
A form 
$\rho \in \Omega^{r}_{q}$
is called $\pi^{r}$-{\it horizontal} (or {\it basic}) if 
$i_{\xi} \rho = 0$
for any $\pi^{r}$-vertical vector field on
$J^{r}Y$.
In local coordinates such a form has the following expression:
\be
\rho = B_{i_{1},...,i_{q}} dx^{i_{1}} \wedge \cdots \wedge dx^{i_{q}}
\label{hor1}
\ee
with
$B_{i_{1},...,i_{q}}$
smooth symmetric functions on $V^{r}$.
We denote the set of basic forms of degree $q$ by
$\Omega^{r}_{q,X}$.
The elements of
$\lambda \in \Omega^{r}_{n,X}$
are called {\it Lagrange forms}. They have the local expression
\be
\lambda = L \theta_{0}
\label{L}
\ee
where
$L$ is a smooth function on $V^{r}$ and
\be
\theta_{0} \equiv dx^{1} \wedge \cdots \wedge dx^{n}.
\ee

We also define some generalisation of the horizontal forms for $q > n$. We say 
that 
$\rho \in \Omega^{r}_{q} \quad q > n$ 
is {\it horizontal} {\it iff} 
\be
i_{\xi_{1}} \cdots i_{\xi_{q-n+1}}\rho = 0
\ee
for any vector fields $\xi_{1},...,\xi_{q-n+1}$ which are $\pi^{r}$-vertical.
The local expression of such a form in the chart $(V^{r},\psi^{r})$ is:
\be
\rho = \sum_{|J_{1}|,...,|J_{q-n}| \leq r} 
A^{J_{1},...,J_{q-n}}_{\nu_{1},...,\nu_{q-n}} dy^{\nu_{1}}_{I_{1}} \wedge 
\cdots dy^{\nu_{q-n}}_{I_{q-n}} \wedge \theta_{0}
\label{hor2}
\ee
with
$A^{J_{1},...,J_{q-n}}_{\nu_{1},...,\nu_{q-n}} $
smooth functions on $V^{r}$ with appropriate symmetry properties.
We conserve the notation $\Omega^{r}_{q,X}$ for these forms.
\newpage

\section{Contact Forms}

In this section we will give a new proof for the structure formula of the
contact forms. This problem was solved by different methods in \cite{Kr6} 
(a sketch of the proof also appears in \cite{Kr4}). Our proof is based on 
some Fock space methods, as said in the Introduction. These kind of methods 
have been employed recently for related problems \cite{Gr1}, \cite{Gr3}. We 
start with the basic definitions and some elementary properties.

\subsection{Basic Definitions and Properties}

By a {\it contact form} we mean any form
$
\rho \in \Omega^{r}_{q}
$
verifying
\be
\left( j^{r}\gamma\right)^{*} \rho = 0
\ee
for any section $\gamma$. We denote by
$\Omega^{r}_{q(c)}$
the set of contact forms of degree $q \leq n$. If one considers only the 
contact forms on an open set
$V \subset Y$
then we emphasize this by writing
$\Omega^{r}_{q(c)}(V)$.
One immediately notes that
$\Omega^{r}_{0(c)} = 0$
and that for  $q > n$ any $q$-form is contact.
It is also elementary to see that the set of all contact forms is an ideal,
denoted by
${\cal C}(\Omega^{r})$,
with respect to the operation $\wedge$. Because the operations of pull-back and
of differentiation are commuting this ideal is left invariant by exterior 
differentiation:
\be
d{\cal C}(\Omega^{r}) \subset {\cal C}(\Omega^{r}).
\label{dif-inv}
\ee

By elementary computations one finds out that for any chart 
$(V,\psi)$
on $Y$, every element of the set
$\Omega^{r}_{1(c)}(V)$
is a linear combination of the following expressions:
\be
\omega^{\sigma}_{j_{1},...,j_{k}} \equiv d y^{\sigma}_{j_{1},...,j_{k}} -
y^{\sigma}_{i,j_{1},...,j_{k}} d x^{i}, \quad k = 0,...,r-1
\label{o}
\ee
or, in multi-index notations
\be
\omega^{\sigma}_{J} \equiv d y^{\sigma}_{J} - y^{\sigma}_{iJ} d x^{i}, \quad
|J| \leq r-1.
\label{o'}
\ee

From the definition above it is clear that the linear subspace of the 1-forms
is generated by
$d x^{i}, \quad \omega^{\sigma}_{J}, \quad (|J| \leq r-1)$
and
$dy^{\sigma}_{I}, |I| = r$.

Formula (\ref{df}) can be now written as follows: for any smooth function on 
$V^{r}$ we have
\be
df = (d_{i}f) dx^{i} + \sum_{|J| \leq r-1} (\partial^{J}_{\sigma}f) 
\omega^{\sigma}_{J} + \sum_{|I| = r} (\partial^{I}_{\nu}f) d y^{\nu}_{I}.
\label{df'}
\ee

We also have the formula
\be
d \omega^{\sigma}_{J} = - \omega^{\sigma}_{Ji} \wedge  d x^{i}, \quad  
|J| \leq r-2.
\ee

Let us now consider an arbitrary form
$\rho \in \Omega^{r}_{q}, \quad q > 1$.
For any $k = 0,...,q$
we define its {\it contact component of order k} to be the form
$p_{k}\rho \in \Omega^{r+1}_{q}$
given by:
\be
p_{k}\rho(j^{r+1}_{x}\gamma)(\xi_{1},...,\xi_{q}) \equiv {1 \over k!(q-k)!} 
\varepsilon^{j_{1},...,j_{q}} \rho(j^{r}_{x}\gamma)
(p\xi_{j_{1}},...,p\xi_{j_{k}},h\xi_{j_{k+1}},...,h\xi_{q}).
\label{k-contact}
\ee

One usually calls
$h\rho \equiv p_{0}\rho$ 
and
$p\rho \equiv \sum_{k=1}^{q} p_{k}\rho$
the {\it horizontal component} and respectively the {\it contact component} of
$\rho$. 
It is useful to particularize the definition above for the case $k = 0$:
\be
h\rho(\xi_{1},...,\xi_{q}) = \rho(h\xi_{1},...,h\xi_{q}).
\ee

A form
$\rho \in \Omega^{r}_{q}$
is called {\it k-contact} if
$p_{j}\rho = 0, \forall j \not= k$
and one says that it has the {\it contact order greater that k} if
$p_{j}\rho = 0, \forall j \leq k-1$.

Now we have a decomposition property, namely

\begin{prop}
For any
$\rho \in \Omega^{r}_{q}$
the following formula is valid:
\be
(\pi^{r+1,r})^{*}\rho = \sum_{k=0}^{q} p_{k}\rho.
\ee
\label{pi-ro-pi}
\end{prop}

{\bf Proof:}
One starts from the definition of the pull-back and use (\ref{hp}) to get
$$
(\pi^{r+1,r})^{*}\rho(j^{r+1}{x}\gamma)(\xi_{1},...,\xi_{q}) =
\rho(j^{r}{x}\gamma)(h\xi_{1}+p\xi_{1},...,h\xi_{q}+p\xi_{q}).
$$

Then one derives the formula in the statement if one uses the definition 
(\ref{k-contact}) and the following combinatorial lemma:

\begin{lemma}

Let $L$ and $M$ be linear finite dimensional spaces and
$\omega: \underbrace{L \times \cdots \times L}_{q-times} \rightarrow M$
a linear and antisymmetric map. Then the following formula is true:
\be
\omega(a_{1}+b_{1},...,a_{q}+b_{q}) = \varepsilon^{i_{1},...,i_{q}}
\sum_{k=0}^{q} {1 \over k!(q-k)!} 
\omega(a_{i_{1}},...,a_{i_{k}},b_{i_{k+1}},...,b_{i_{q}})
\ee
for all
$a_{1},...,a_{q},b_{1},...,b_{q} \in L$.
\label{decomp}
\end{lemma}

{\bf Proof:} 
Is straightforward by induction on $q$. 
$\qed$

One can use the proposition above to deduce that $\rho$ is contact {\it iff}  
it verifies  
$h\rho = 0$
and that the expressions
$p_{1}\rho,...,p_{q}\rho$
are contact forms in the sense of the definition given at the beginning of this
Section. Finally we remind the fact that the horizontalization operation
$h: \Omega^{r} \rightarrow \Omega^{r+1}$
verifies the properties:
\be
h(\mu + \nu) = h\mu + h\nu, \quad h(\mu \wedge \nu) = h\mu \wedge h\nu
\label{morphism}
\ee
for any forms
$\mu, \nu \in \Omega^{r},$
\be
h d x^{i} = d x^{i}, \quad h d y^{\sigma}_{J} = y^{\sigma}_{Ji} d x^{i} \quad
(|J| \leq r)
\label{h-on-xy}
\ee
and also
\be
h f = (\pi^{r+1,r})^{*} f = f \circ \pi^{r+1,r}
\ee
for all smooth functions on $V^{r}$. Moreover, these three properties given
above determine {\it uniquely} the map $h$.
\newpage

\subsection{The Structure Theorem for Contact Forms}

In this subsection we give a new proof to the following fact.

\begin{thm}

Let
$(V,\psi)$
an adapted chart on the fibre bundle $Y$ and let
$\rho \in \Omega^{r}_{q}(Y), \quad q = 2,...,n$.
Then $\rho$ is contact {\it iff} it has the following expression in the 
associated chart:
\be
\rho = \sum_{|J|\leq r-1} \omega^{\sigma}_{J} \wedge \Phi^{J}_{\sigma} +
\sum_{|I|= r-1} d\omega^{\sigma}_{I} \wedge \Psi^{I}_{\sigma}
\label{contact-q}
\ee
where 
$\Phi^{J}_{\sigma} \in \Omega^{r}_{q-1}$
and
$\Psi^{I}_{\sigma} \in \Omega^{r}_{q-2}$
can be arbitrary forms on $V^{r}$. (We adopt the convention that 
$\Omega^{r}_{q} \equiv 0, \forall q < 0$).
\label{contact}
\end{thm}

{\bf Proof:}
(i) If $\rho$ has the expression (\ref{contact-q}) one uses (\ref{morphism}) to
obtain that it is a contact form. We now concentrate on the converse statement.
Firstly we need a canonical expression for any 
$\rho \in \Omega^{r}_{q}$.
We start from the fact that the forms
$d x^{i}, \quad \omega^{\sigma}_{J} \quad (|J| \leq r-1)$
and
$d y^{\sigma}_{I} \quad (|I| = r)$
are a basis in the linear space of 1-forms. The form $\rho$ is a polynomial of
degree $q$ in these forms (with respect to the product $\wedge$). We separate
all the terms containing at least a factor 
$\omega^{\sigma}_{J}$
and get a decomposition
\be
\rho = \rho_{0} + \rho'
\label{ro}
\ee
where
$\rho_{0}$
has the structure
\be
\rho_{0} = \sum_{|J|\leq r-1} \omega^{\sigma}_{J} \wedge \Phi^{J}_{\sigma} 
\label{ro0}
\ee
and $\rho'$ is a polynomial of degree $q$ only in 
$d x^{i}$
and
$d y^{\sigma}_{I} \quad (|I| = r)$.
It is clear the one can write it as follows:
\be
\rho' = \sum_{s=0}^{q} {1 \over s!(q-s)!} \sum_{|I_{1}|,...,|I_{s}|=r}
A^{I_{1},...,I_{s}}_{\sigma_{1},...,\sigma_{s},i_{s+1},...,i_{q}}
dy^{\sigma_{1}}_{I_{1}} \wedge \cdots \wedge dy^{\sigma_{s}}_{I_{s}} \wedge
dx^{i_{s+1}} \wedge \cdots \wedge dx^{i_{q}},
\label{canonical}
\ee
where 
$A^{I_{1},...,I_{s}}_{\sigma_{1},...,\sigma_{s},i_{s+1},...,i_{q}}$
are smooth functions on $V^{r}$ and can be assumed to verify the following
symmetry property:
\be
A^{I_{P(1)},...,I_{P(s)}}_{\sigma_{P(1)},...,\sigma_{P(s)},
i_{Q(s+1)},...,i_{Q(q)}} = (-1)^{|P|+|Q|}
A^{I_{1},...,I_{s}}_{\sigma_{1},...,\sigma_{s},i_{s+1},...,i_{q}}
\label{symmetry}
\ee
for any permutations
$P \in {\cal P}_{s}, \quad Q \in {\cal P}_{q-s}$.

One must impose the condition
$h\rho = h\rho' = 0$.
Using the relations (\ref{morphism}) and (\ref{h-on-xy}) one gets
$$
h\rho' = \sum_{s=0}^{q} {1 \over s!(q-s)!} \sum_{|I_{1}|,...,|I_{s}|=r}
A^{I_{1},...,I_{s}}_{\sigma_{1},...,\sigma_{s},i_{s+1},...,i_{q}}
y^{\sigma_{1}}_{I_{1}i_{1}} \cdots y^{\sigma_{s}}_{I_{s}i_{s}}
dx^{i_{1}} \wedge \cdots \wedge dx^{i_{q}},
$$
so the relation above is equivalent to
\be
{\cal S}^{-}_{i_{1},...,i_{q}} \sum_{|I_{1}|,...,|I_{s}|=r}
A^{I_{1},...,I_{s}}_{\sigma_{1},...,\sigma_{s},i_{s+1},...,i_{q}}
y^{\sigma_{1}}_{I_{1}i_{1}} \cdots y^{\sigma_{s}}_{I_{s}i_{s}} = 0, \quad
s = 0,...,q.
\ee

Here
${\cal S}^{-}_{i_{1},...,i_{q}}$
is the antisymmetrization operator in the indices
$i_{1},...,i_{q}$
defined by the formula (\ref{sa-indices}) from the Appendix.

We apply to the relation above derivative operators of the type
$\partial^{I}_{\nu}, \quad (|I| = r+1)$
and obtain, taking into account the symmetry property (\ref{symmetry}) the 
following relations (see \cite{Kr4} , \cite{Kr6}):
\be
A_{i_{1},...,i_{q}} = 0
\ee
and
\be
{\cal S}^{-}_{i_{1},...,i_{q}} {\cal S}^{+}_{I_{1}p_{1}} \cdots 
{\cal S}^{+}_{I_{s}p_{s}} 
A^{I_{1},...,I_{s}}_{\sigma_{1},...,\sigma_{s},i_{s+1},...,i_{q}}
\delta^{p_{1}}_{i_{1}} \cdots \delta^{p_{s}}_{i_{s}} = 0, \quad
s = 1,...,q.
\label{A}
\ee

(ii) We analyse in detail the relation above. This is the point where we depart
from the idea of the proof from \cite{Kr6}. Let us introduce the following
tensor spaces:
$$
{\cal H}_{s} \equiv {\cal F}^{(-)}(\R^{n}) \otimes 
\underbrace{{\cal F}^{(+)}(\R^{n}) \otimes \cdots \otimes
{\cal F}^{(+)}(\R^{n})}_{s-times}
$$
where
${\cal F}^{(\pm)}(\R^{n})$
are the symmetric (coresp. $+$) and the antisymmetric (coresp. $-$) Fock
spaces (see the Appendi). We have the well known decomposition in subspaces 
with fixed number of ``bosons" and ``fermions":
$$
{\cal H}_{s} = \oplus_{k=0}^{r} \oplus_{l_{1},...,l_{s} \geq 0}
{\cal H}_{k,l_{1},...,l_{s}}.
$$

We make the convention that
${\cal H}_{k,l_{1},...,l_{s}} \equiv 0$
if any one of the indices
$k,l_{1},...,l_{s}$
is negative or if
$k > n$.
Then we can consider
$A^{I_{1},...,I_{s}}_{\sigma_{1},...,\sigma_{s},i_{s+1},...,i_{q}}$
as the components of a tensor
$$
A_{\sigma_{1},...,\sigma_{s}} \in 
{\cal H}_{q-s,\underbrace{r,...,r}_{s-times}}.
$$  

We can write in an extremely compact way the relation (\ref{A}) if we use the
creation and the annihilation fermionic operators
$a^{*i}, a_{i}, \quad (i = 1,...,n)$ 
and the corresponding creation and annihilation bosonic operators 
$b^{*}_{(\alpha)i}, b^{i}_{(\alpha)} \quad (\alpha = 1,...,s; i = 1,...,n).$
We will using conventions somewhat different from that used in quantum
mechanics (see the Appendix). One introduces the operators
\be
B_{\alpha} \equiv b^{*}_{(\alpha)i} a^{*i}, \quad (\alpha = 1,...,s)
\label{B}
\ee
and proves by elementary computations that the relations (\ref{A}) are
equivalent to:
\be
B_{1}\cdots B_{s} A_{\sigma_{1},...,\sigma_{s}} = 0, \quad (s = 1,...,q).
\label{A-compact}
\ee

Indeed, it is elementary to check that we have for instance:
\be
( B_{1} A)^{\{j_{0},...,j_{k}\},I_{2},...,I_{s}}_{\sigma_{1},...,\sigma_{s},
i_{0},...,i_{l}} = 
{\cal S}^{+}_{j_{0},...,j_{k}} {\cal S}^{-}_{i_{0},...,i_{l}}
A^{\{j_{1},...,j_{k}\},I_{2},...,I_{s}}_{\sigma_{1},...,\sigma_{s},
i_{1},...,i_{l}} \delta^{j_{0}}_{i_{0}}.
\label{B-alpha}
\ee
and similarly for $\alpha = 2,\dots,s$.

(iii) Equations of the type appearing in the (\ref{A-compact}) can be
completely analysed with the help of two lemmas, which are the backbone of our
paper. Firstly we analyse the case $s =  q$. In this case
$A_{\sigma_{1},...,\sigma_{q}} \in 
{\cal H}_{0,\underbrace{r,...,r}_{q-times}}.$

\begin{lemma}
Let
$X \in {\cal H}_{0,r_{1},...,r_{q}}, \quad (r_{1},...,r_{q} \in \N;
\quad q = 2,...,n)$
verifying the relation
\be
B_{1}\cdots B_{q} X = 0. 
\ee
Then $X = 0$. 
\label{s=q}
\end{lemma}

{\bf Proof:}
We use complete induction over $n$. For $n = 2$ we have only the case $q = 2$.
The generic form of $X$ is
$$
X = \sum_{i=0}^{r_{1}} \sum_{j=0}^{r_{2}} c_{i,j} 
(b^{*}_{(1)1})^{i} (b^{*}_{(1)2})^{r_{1}-i} (b^{*}_{(2)1})^{j} 
(b^{*}_{(2)2})^{r_{2}-j} \Omega;
$$
here $\Omega$ is the vacuum vector in ${\cal H}$. We extend the coefficients  
$c_{i,j}$ to all integer values of the indices by taking
$c_{i,j} = 0$ if $i$ (or $j$) is outside the set $\{1,...,r_{1}\}$ (or
$\{1,...,r_{2}\}$).  One computes that the relation 
$B_{1}B_{2}X = 0$ 
is equivalent to
$c_{i-1,j} = c_{i,j-1}$;
this relation gives by recurrence that $c_{i,j} = 0, \quad \forall i,j \in 
{\bf Z}$ i.e. $X = 0$.

We assume the assertion from the statement to be true for a given $n$ and we
prove it for $n+1$. In this case the indices will take the values
$0,1,...,n$
and we must make the substitution
$$
B_{\alpha} \rightarrow \tilde{B}_{\alpha} \equiv B_{\alpha} + 
b^{*}_{(\alpha)0} a^{*0}, \quad (\alpha = 1,...,q).
$$

By hypothesis we have
$$
\prod_{\alpha =1}^{q} \tilde{B}_{\alpha} X = 0.
$$

One immediately shows this to be equivalent to the following two equations:
\be
\prod_{\alpha =1}^{q} B_{\alpha} X = 0
\label{n+1,1}
\ee
and
\be
\sum_{\alpha = 1}^{q} (-1)^{\alpha -1} b_{(\alpha)}^{*0} B_{1} \cdots
\hat{B}_{\alpha} \cdots B_{q} X = 0.
\label{n+1,2}
\ee

The generic expression for $X$ is
\be
X = \sum_{t_{1}=0}^{r_{1}}\cdots \sum_{t_{q}=0}^{r_{q}}
\left( b_{(1)}^{*0}\right)^{t_{1}} \cdots \left( b_{(q)}^{*0}\right)^{t_{q}}
X_{t_{1},...,t_{q}}
\label{X,n+1}
\ee
where
$X_{t_{1},...,t_{q}} \in {\cal H}_{0,r_{1}-t_{1},...,r_{q}-t_{q}}$
are tensors obtained from the vacuum by applying only polynomials in the
creation operators
$b_{(\alpha)}^{*i}, \quad (\alpha = 1,...,n; 1 = 1,...,n)$.

Then the equation (\ref{n+1,1}) is equivalent to
\be
\prod_{\alpha =1}^{q} B_{\alpha} X_{t_{1},...,t_{q}} = 0, \forall 
t_{1},...,t_{q} \in {\bf Z}. 
\ee

If $q \leq n$ we can apply the induction hypothesis and obtain
$X_{t_{1},...,t_{q}} = 0, \forall t_{1},...,t_{q} \in {\bf Z}$
i.e. $X = 0$.
So it remains to study only the case $q = n + 1$. In this case one notices that
the equation (\ref{n+1,1}) becomes an identity because
$\prod_{\alpha =1}^{n+1} B_{\alpha} = 0$;
indeed, in the left hand side, at least one of the operators
$a^{*i}, \quad (i = 1,...,n)$
appears twice. So, in this case we are left with the second equation
(\ref{n+1,2}):
\be
\sum_{\alpha = 1}^{n+1} (-1)^{\alpha -1} b_{(\alpha)}^{*0} B_{1} \cdots
\hat{B}_{\alpha} \cdots B_{n+1} X = 0.
\ee

One must substitute here the generic expression for $X$ (\ref{X,n+1}) and the
following relation is produced:
\be
\sum_{\alpha = 1}^{n+1} (-1)^{\alpha -1} b_{(\alpha)}^{*0} B_{1} \cdots
\hat{B}_{\alpha} \cdots B_{n+1} 
X_{t_{1},...t_{\alpha -1},t_{\alpha}-1,t_{\alpha +1},...,t_{n+1}} = 0,
\forall t_{1},...,t_{n+1} \in {\bf Z}.
\label{t}
\ee

This relation can used to prove that $X = 0$. We outline the argument. First we
take 
$t_{1} = t + 1 \quad (t \geq 0), t_{2} = \cdots = t_{n+1} = 0$
in the relation above and get
$$
B_{2} \cdots B_{n+1} X_{t,\underbrace{0,...,0}_{n-times}} = 0.
$$
The induction hypothesis can be applied and we get
$X_{t,\underbrace{0,...,0}_{n-times}} = 0, \quad \forall t \in {\bf Z}$.
Similarly, one can obtain
$X_{0,...,0,t,0,...,0} = 0, \quad \forall t \in {\bf Z}$
where the index $t$ can be positioned anywhere.

Now one can prove, by induction on $p$ that
$X_{t_{1},,,,.t_{n+1}} = 0$
if at least $n - p + 1$ indices are equal to zero. For $p = 1$ this statement
has just have been proved above. We accept it for a given $p$ and prove it for
$p + 1$. We take in (\ref{t}) 
$t_{p+2} = \cdots t_{n+1} = 0$
and obtain
\be
\sum_{\alpha = 1}^{p+1} (-1)^{\alpha -1} b_{(\alpha)}^{*0} B_{1} \cdots
\hat{B}_{\alpha} \cdots B_{n+1} 
X_{t_{1},...t_{\alpha -1},t_{\alpha}-1,t_{\alpha +1},...,t_{n+1}} = 0.
\label{t-p}
\ee

If
$t_{1} = t, t_{2} = t_{3} = \cdots = t_{p+1} = 1$
one can use the induction hypothesis to get:
$$
X_{t,\underbrace{1,...,1}_{p-times},\underbrace{0,...,0}_{(n-p)-times}} = 0,
\quad \forall t \in {\bf Z}.
$$

Now one uses the relation above and (\ref{t-p}) to prove that
\be
X_{t_{1},...,t_{p+1},\underbrace{0,...,0}_{(n-p)-times}} = 0;
\label{t-p+1}
\ee
this is done by recurrence on
$t_{2} + t_{3} + \cdots + t_{p+1}$.

Evidently the argument leading to the relation (\ref{t-p+1}) works in the same
way for any positioning of the $p + 1$ indices. The induction is finished and
we get
$X_{t_{1},...,t_{n+1}} = 0, \quad \forall t_{1},...,t_{n+1} \in {\bf Z}$
i.e. $X = 0$.
$\nabla$

If we apply this lemma to the equation (\ref{A-compact}) for $s = q$ we get
\be
A^{I_{1},...,I_{q}}_{\sigma_{1},...,\sigma_{q}} = 0.
\ee

(iv) We still have to analyse the case $0 < s < q$ of the equation
(\ref{A-compact}). We can analyse immediately the case $s = 1$ i.e. the
equation
\be
B_{1} X = 0
\label{B1}
\ee
if we use (\ref{BB-star}); we get that
$$
(n+r-q+1) X = BB^{*} X
$$
i.e. $X$ is of the form
\be
X = BX_{1}
\ee
for some $X_{1} \in {\cal H}_{q-2,r-1}$.

(v) We generalize this result to all $s = 1,...,q-1$.

\begin{lemma}
Let
$
X \in {\cal H}_{k,r_{1},...,r_{s}},\quad 0 \leq k < n, \quad 0 < s \leq n,
\quad s + k \leq n
$
Then $X$ verifies the equation
\be
B_{1} \cdots B_{s} X = 0
\label{B-X}
\ee
{\it iff} it is of the form
\be
X = \sum_{\alpha =1}^{s} B_{\alpha} X_{\alpha}
\label{X-B}
\ee
for some
$X_{\alpha} \in 
{\cal H}_{k-1,r_{1},...r_{\alpha -1},r_{\alpha}-1,r_{\alpha +1},...,r_{s}}.
$
\label{s<q}
\end{lemma}

{\bf Proof:}

From (\ref{B2}) the implication (\ref{X-B}) $\Rightarrow$ (\ref{B-X}) is
obvious. We prove now the converse statement, as before, by induction on $n$.
For 
$n = 2$ 
we can have 
$k = 0,1$. 
In the first case the statement is true according to lemma \ref{s=q} and the
second case was analysed before at (iv).  We suppose that the statement is true
for a given $n$ and we prove it for 
$n + 1$. 
With the same notations as in lemma \ref{s=q} we have an equation of the type:
\be
\prod_{\alpha =1}^{s} \tilde{B}_{\alpha} X = 0.
\label{s}
\ee
where
$$
\tilde{B}_{\alpha} \equiv B_{\alpha} + 
b^{*}_{(\alpha)0} a^{*0}, \quad (\alpha = 1,...,q).
$$

The generic form for $X$ is
\be
X = X_{0} + a^{*0} Z
\ee
with
$X_{0} \in {\cal H}_{k,r_{1},...,r_{s}}$
and
$Z \in {\cal H}_{k-1,r_{1},...,r_{s}}$
tensors obtained from the vacuum by applying only polynomials in
$a^{*i} \quad (i = 1,...,n)$
and
$b^{*}_{(\alpha)\mu} \quad (\mu = 0,...,n; \alpha = 1,...,s)$.

The equation (\ref{s}) becomes equivalent to the following two equations:
\be
\prod_{\alpha =1}^{s} B_{\alpha} X_{0} = 0
\label{s1}
\ee
and
\be
\prod_{\alpha =1}^{s} B_{\alpha} Z = \sum_{\alpha =1}^{s} (-1)^{s-\alpha}
b^{*0}_{(\alpha)} B_{1} \cdots \hat{B}_{\alpha} \cdots B_{s} X_{0}.
\label{s2}
\ee

As in lemma \ref{s=q}, the generic form of $X_{0}$ is
\be
X_{0} = \sum_{t_{1}=0}^{r_{1}}\cdots \sum_{t_{q}=0}^{r_{q}}
\left( b_{(1)}^{*0}\right)^{t_{1}} \cdots \left( b_{(q)}^{*0}\right)^{t_{q}}
X_{t_{1},...,t_{q}}.
\label{X0}
\ee

Then the equation (\ref{s1}) is equivalent to
\be
\prod_{\alpha =1}^{s} B_{\alpha} X_{t_{1},...,t_{q}} = 0, \quad \forall
t_{1},...,t_{q} \in {\bf Z}.
\label{s1a}
\ee

We have two distinct cases:

(a) $k < n$ and $s \leq n$

In this case we can apply the induction hypothesis to the relation above and
obtain in the end that $X_{0}$ is of the following form:
$$
X_{0} = \sum_{\alpha =1}^{s} B_{\alpha} X_{\alpha}.
$$

If we introduce this expression into the equation (\ref{s2}) we easily get:
$$
\prod_{\alpha =1}^{s} B_{\alpha} \left( Z - 
\sum_{\alpha =1}^{s} b^{*0}_{(\alpha)} X_{\alpha}\right) = 0 
$$
so again we can apply the induction hypothesis to obtain that $Z$ has the
following structure:
$$
Z = \sum_{\alpha =1}^{s} b^{*0}_{(\alpha)} X_{\alpha} +
\sum_{\alpha =1}^{s} B_{\alpha} Z_{\alpha}.
$$

Now we define
$$
\tilde{X}_{\alpha} \equiv X_{\alpha} - a^{*0} Z_{\alpha}
$$
and obtain from the previous relations that
$$
X = \sum_{\alpha =1}^{s} \tilde{B}_{\alpha} \tilde{X}_{\alpha}
$$
which finishes the proof.

(b) If $k = n$ then from the restrictions on $s$ and $k$ we necessarily have 
$s = 1$ and we can use (iv). If $s = n + 1$ then the same restrictions fix 
$k = 0$ and we can apply lemma \ref{s=q}.
$\nabla$

We can apply the lemma above to the equation (\ref{A-compact}) for the cases
$s = 2,...,q-1$
and
obtain that the tensors
$A_{\sigma_{1},...,\sigma_{s}}$
have the following structure
$$
A_{\sigma_{1},...,\sigma_{s}} = \sum_{\alpha =1}^{s} B_{\alpha}
X^{\alpha}_{\sigma_{1},...,\sigma_{s}}
$$
for some tensors
$X^{\alpha}_{\sigma_{1},...,\sigma_{s}}$.

Now it is the moment to use again full index notation. The relation above means
that the expression
$A^{I_{1},...,I_{s}}_{\sigma_{1},...,\sigma_{s},i_{s+1},...,i_{q}}$
is a sum of terms such that every term has at least a factor of the type
$\delta^{j}_{i}$
where the index $j$ belongs to some $I_{p}$ and the index $i$ is one of the
indices 
$i_{s+1},...,i_{q}$.
That's it if, say, 
$I_{1} = \{j_{1},...,j_{r}\}$
then
$$
A^{\{j_{1},...,j_{r}\},...}_{\sigma_{1},...,\sigma_{s},i_{s+1},...,i_{q}} =
\sum_{u=1}^{r} \sum_{v=s+1}^{q} \delta^{j_{u}}_{i_{v}} 
\tilde{A}^{\{j_{1},...,\hat{j}_{u},...,j_{r}\},...}_{\sigma_{1},...,\sigma_{s},
i_{s+1},...,\hat{i}_{v},...,i_{q}} + \cdots.
$$

If we substitute the preceding relation into the expression of $\rho'$
(see (\ref{canonical})) we obtain a sum of contributions of the type
$$
\tilde{A}^{\{j_{1},...,\hat{j}_{u},...,j_{r}\},...}_{\sigma_{1},...,\sigma_{s},
i_{s+1},...,\hat{i}_{v},...,i_{q}} 
d y^{\sigma_{1}}_{\{j_{1},...,j_{r}\}} \wedge d y^{\sigma_{2}}_{I_{2}} \cdots
\wedge d y^{\sigma_{s}}_{I_{s}} \wedge d x^{i_{s+1}} \wedge \cdots d x^{j_{u}}
\cdots \wedge d x^{i_{q}}
$$
i.e. a sum of terms containing the expression
$
d y^{\sigma_{1}}_{\{j_{1},...,j_{r}\}} \wedge d x^{j_{u}} = 
- d \omega^{\sigma_{1}}_{\{j_{1},...,\hat{j}_{u},...,j_{r}\}}.
$

So the contribution $\rho'$ to the contact form $\rho$ form gives the second
terms from the statement of the theorem (v. the formul\ae\/ (\ref{ro}) and 
(\ref{contact-q}).
$\qed$

Let us note for further use that one can combine lemmas \ref{s=q} and \ref{s<q}
in a single result:

\begin{lemma}
Let
$X \in {\cal H}_{k,r_{1},...,r_{s}}$
Then $X$ verifies the equation (\ref{B-X}) {\it iff} it is of the form
\be
X = \sum_{\alpha =1}^{s} B_{\alpha} X_{\alpha}
\ee
for some
$X_{\alpha} \in 
{\cal H}_{k-1,r_{1},...r_{\alpha -1},r_{\alpha}-1,r_{\alpha +1},...,r_{s}}.
$
\label{BXB}
\end{lemma}

\begin{rem}
One can show in fact that the decomposition of an arbitrary form given by
the formul\ae\/ (\ref{ro}), (\ref{ro0}) and (\ref{canonical}) can be refined
with the help of the so called {\it trace decomposition identity} \cite{Kr7}.
Although we do not need this more refined decomposition we will provide an
alternative proof of this fact, based on the same tricks, in the Appendix. 
This will emphasize once more the power of our method.
\end{rem}
\newpage

\subsection{Some Properties of the Contact Forms}   

We start with the transformation formula for the contact forms. We have:

\begin{prop}
Let
$(V,\psi)$
and
$(\bar{V},\bar{\psi})$
be two overlapping charts on $Y$. Then on
$V^{r} \cap \bar{V}^{r}$
the following formula is true:
\be
\bar{\omega}^{\sigma}_{I} = \sum_{|J| \leq |I|} (\partial^{J}_{\nu}
\bar{y}^{\sigma}_{I}) \omega^{\nu}_{J} \quad (|I| \leq r-1).
\label{o-bar-o}
\ee
\end{prop}

{\bf Proof:}
The proof is based on simple manipulations of the formula (\ref{df}), the
definition (\ref{o}) of the 1-contact forms and use is also made of lemma
\ref{tr}. 
$\qed$

An element 
$T \in \Omega^{s}_{n+1,X}$
is called a {\it differential equation} if 
$i_{\xi} T = 0$
for any $\pi^{s,1}$-vertical vector field. In the chart $V^{s}$ the 
differential equation $T$ has the following expression:
\be
T = T_{\sigma} \omega^{\sigma} \wedge \theta_{0}.
\label{T}
\ee
(see (\ref{hor2})). Using (\ref{o-bar-o}) one can indeed see that $T$ has this
form in any chart; explicitly, the transformation formula is:
\be
T_{\sigma} = {\cal J} (\partial_{\sigma} \bar{y}^{\nu}) \bar{T}_{\nu}
\label{t-bar-t}
\ee
where ${\cal J}$ is the Jacobian of the chart transformation on $X$:
\be
{\cal J} \equiv det\left({\partial \bar{x}^{i} \over \partial x^{j}}\right).
\label{J}
\ee

If $\gamma$ is a section of the fibre bundle $\pi: Y \rightarrow X$ then on 
says that {\it it verifies the differential equation} $T$ {\it iff} we have
\be
(j^{s}\gamma)^{*} i_{Z} T = 0
\ee
for any vector field $Z$ on $J^{s}Y$. In local coordinates we have on $V^{s}$:
\be
T_{\sigma} \circ j^{s}\gamma = 0 \quad (\sigma = 1,...,m).
\label{eq-diff}
\ee

Another important property of the contact ideal is that it behaves naturally
with respect to prolongations. More precisely, let
$\pi_{i}: Y_{i} \rightarrow X_{i}, \quad i = 1,2$
be two fibre bundles and 
$\phi: Y_{1} \rightarrow Y_{2}$
a fibre bundle morphism. Then the prolongation
$j^{r}\phi$
(defined in the end of section 2) verifies:
\be
(j^{r}\phi)^{*}{\cal C}(\Omega^{r}(Y_{1})) \subset {\cal C}(\Omega^{r}(Y_{2})).
\label{pr-cont}
\ee

The proof follows directly from the definition of a contact form. As a
consequence, if $\xi$ is a projectable vector field on the fibre bundle
$Y$, then the Lie derivative of 
$j^{r}\xi$ 
leaves the contact ideal invariant:
\be
L_{j^{r}\xi} {\cal C}(\Omega^{r}(Y)) \subset {\cal C}(\Omega^{r}(Y)).
\label{Lie}
\ee

This formula can be used to find out the explicit expression of 
$j^{r}\xi$ 
\cite{An}. Indeed, if in the chart 
$(V,\psi)$
we have
\be
\xi = a^{i}(x) {\partial \over \partial x^{i}} + b^{\sigma}(x,y)
\partial_{\sigma} 
\ee
with $a^{i}$ and $b^{\sigma}$ smooth function, then $j^{r}\xi$ must have the
following expression in the associated chart 
$(V^{r},\psi^{r})$: 
\be
j^{r}\xi = a^{i}(x) {\partial \over \partial x^{i}} + \sum_{|J| \leq r} 
b^{\sigma}_{J} \partial_{\sigma}^{J}. 
\ee

One imposes an equivalent form of (\ref{Lie}), namely
\be
L_{j^{r}\xi} \omega^{\sigma}_{J} \in {\cal C}(\Omega^{r}(Y)), \quad
|J| \leq r-1.
\label{Lie'}
\ee

The left hand side of this relation can be computed explicitly:
$$
L_{j^{r}\xi} \omega^{\sigma}_{J} = (d_{i} b^{\sigma}_{J} - b^{\sigma}_{Ji} -
y^{\sigma}_{Jl} d_{i}a^{l}) dx^{i} + \sum_{|I|=r} (\partial^{I}_{\nu}
b^{\sigma}_{J}) dy^{\nu}_{I} + {\rm contact~ terms}.
$$

The following recurrence formula for the coefficients 
$b^{\sigma}_{J}$
follows:
\be
b^{\sigma}_{Ji} = d_{i} b^{\sigma}_{J} - y^{\sigma}_{Jl} d_{i}a^{l}, \quad
|J| \leq r-1;
\ee
we also have:
\be
\partial^{I}_{\nu} b^{\sigma}_{J} = 0, \quad |I| = r.
\ee

In particular, if $\xi$ is an evolution i.e. it has the local expression 
(\ref{ev}), then we have
\be
j^{r}\xi = \sum_{|J| \leq r} (d_{J}\xi^{\sigma}) \partial_{\sigma}^{J}.
\label{pr-ev}
\ee

One may wonder what is the expression of the prolongation $j^{r}\phi$ (where 
$\phi$ is a bundle morphism of the fibre bundle $Y$. One can proceed in 
complete analogy with the computations above. If $\phi$ has the following 
expression in the chart
$(V,\psi)$
\be
\phi(x^{i},y^{\sigma}) = (f^{i},F^{\sigma})
\ee
then we must have in the associated chart: 
$(V^{r},\psi^{r})$:
\be
j^{r}\phi(x^{i},y^{\sigma},y^{\sigma}_{j},...,y^{\sigma}_{j_{1},...,j_{r}}) =
(f^{i},F^{\sigma},F^{\sigma}_{j},...,F^{\sigma}_{j_{1},...,j_{r}})
\label{pr-transf}
\ee
where
$
F^{\sigma}_{j_{1},...,j_{k}}, \quad j_{1} \leq j_{2} \leq \cdots \leq j_{k},
\quad k = 1,...,r
$
are smooth functions on the chart $V^{r}$. One starts from an equivalent form
of (\ref{pr-cont}), namely:
\be
(j^{r}\phi)^{*} \omega^{\sigma}_{J} \subset {\cal C}(\Omega^{r}(Y)
\quad |J| \leq r-1
\label{pr-cont'}
\ee
and computes the left hand side:
$$
(j^{r}\phi)^{*} \omega^{\sigma}_{J} = \left( d_{l}F^{\sigma}_{J} - 
F^{\sigma}_{Ji} {\partial f^{l} \over \partial x^{i}} \right) dx^{i} +
\sum_{|J|=r} (\partial^{I}_{\nu} F^{\sigma}_{J}) dy^{\sigma}_{I} + 
{\rm contact~ terms}, \quad |J| \leq r-1.
$$

The condition (\ref{pr-cont'}) above gives us a recurrence formula for the 
functions $F^{\sigma}_{J}$:
\be
F^{\sigma}_{Ji} = Q^{l}_{i} d_{l}F^{\sigma}_{J}\quad |J| \leq r-1;
\ee
we also have
\be
\partial^{I}_{\nu} F^{\sigma}_{J} = 0 \quad |I|=r.
\ee

Let us note that the recurrence formula above formally coincide with the
recurrence formula from lemma \ref{tr}.

We close this subsection reminding another important construction appearing
when one considers the so-called variationally trivial Lagrangians (they will 
be defined in section 7). We introduce the following subset of the space of 
basic forms:
\be
{\cal J}^{r}_{q} \equiv \{\rho \in \Omega^{r}_{q,X} | \exists \nu \in
\Omega^{r-1}_{q}, \quad s.t.\quad \rho = h\nu\}.
\label{i}
\ee

One notes that this subspace is closed with respect to the wedge product
$\wedge$ and also that the operator
$D: {\cal J}^{r}_{q} \rightarrow {\cal J}^{r}_{q+1}$
given by
\be
Dh\nu \equiv hd\nu
\label{D}
\ee
is well defined \cite{AD}. The operator $D$ is called {\it total exterior
derivative}. We list some elementary properties of this operator directly
deductible from the definition.
\be
D \circ D = 0
\ee
and
\be
Dh = hD.
\ee

Moreover, $D$ is a derivation completely determined by the following relations:
\be
Df = (d_{i}f) dx^{i},\quad \forall f \in {\cal J}^{r}_{0},
\ee
and
\be
D(dx^{i}) = 0, \quad Dhdy^{\sigma}_{J} = 0, \quad |J| \leq r.
\ee

We remind the reader that we have remarked before that the operators 
$d_{i}$
are not vector fields. However, we have the following relation, which is the
next best thing except a vector field. If 
$f \in \Omega^{r-1}_{0}$
and we have two overlapping charts 
$(V,\psi)$
and
$(\bar{V},\bar{\psi})$
on $Y$ then we have on the intersection 
$V^{r} \cap \bar{V}^{r}$
the following relation:
\be
\bar{d}_{i}f =  Q^{j}_{i} d_{j}f
\label{d-bar-d}
\ee
(where 
$\bar{d}_{j}$ 
are the formal derivatives in the chart 
$\bar{V}^{r}$ and the matrix $Q$ has been defined previously: it is the inverse
of the Jacobian matrix of the chart transformation - see (\ref{x-bar-x})).
The proof is elementary and consists in expressing the (globally defined)
operator $D$ in both charts.
\newpage

\section{Strongly Contact Forms}

The concept of strongly contact form has been introduced by Krupka \cite{Kr4}.
The idea is to observe that the definition of the contact forms is trivially
satisfied if the degree of the form is $q \geq n+1$. So, it is natural to try a
generalization of the contact forms in this case. It seems plausible to use
instead of the horizontalization operator $h$ some other projection $p_{k}$
from those introduced in the beginning of the preceding section. The proper
definition is the following. Let 
$q = n+1,...,N \equiv dim(J^{r}Y) = m{n+r \choose n}$ 
and let 
$\rho \in \Omega^{r}_{q}$.
One says that $\rho$ is a {\it strongly contact form} {\it iff} its contact
component of order $q - n$ vanishes i.e.
\be
p_{q-n}\rho = 0.
\label{str=cont}
\ee

For a certain uniformity of notations, we denote these forms by 
$\Omega^{r}_{q(c)}$.
We need a structure formula for strongly contact forms, i.e. an analogue of
theorem \ref{contact}. First we need some properties of the projections
$p_{k}$ \cite{Kr4} - \cite{Kr6}. 

\begin{lemma}
If 
$\rho \in \Omega^{r}_{q}$  
and 
$\rho' \in \Omega^{r}_{t}$ 
then the following formula is true:
\be
p_{k}(\rho \wedge \rho') = \sum_{l+s=k} p_{l}\rho \wedge p_{s}\rho' \quad
\forall k \geq 0
\ee
\label{p-wedge}
where we make the convention that
\be
p_{k}\rho \equiv 0, \quad {\rm if} \quad k > {\rm deg}(\rho),\quad {\rm or} 
\quad k < 0.
\ee
\label{p-k}
\end{lemma}

{\bf Proof:} 
Is based on induction on $q$. For 
$q = 1$ 
one starts from the definition (\ref{k-contact}) of the operator 
$p_{k}$
and from the definition for the wedge product, which in our case is:
$$
(\rho \wedge \rho')(\xi_{0},...,\xi_{t}) = \sum_{i=0}^{t} (-1)^{i}
\rho(\xi_{i}) \rho'(\xi_{0},...,\hat{\xi}_{i},...,\xi_{t}).
$$

Then, one supposes the formula true for 
$1, 2,\dots,q$ 
and proves it for 
$q + 1$.
One does not loose generality if one supposes that $\rho$ is of the form
$\rho = \rho_{1} \wedge \rho_{2}$
with 
${\rm deg}(\rho_{1}) = 1$
and
${\rm deg}(\rho_{2}) = q$
and the result for
$q + 1$
is obvious.
$\qed$

As a corollary, we have

\begin{cor}
If 
$\rho_{i} \in \Omega^{r}_{q_{i}}, \quad i = 1,...,l$
then the following formula is true:
\be
p_{k}(\rho_{1} \wedge \dots \wedge \rho_{l}) = \sum_{s_{1}+\dots +s_{l}=k}
p_{s_{1}}\rho_{1} \wedge \dots \wedge p_{s_{l}} \rho_{l}.
\ee
In particular, if the order of contactness of the forms
$\rho_{i} \in \Omega^{r}_{q_{i}}, \quad i = 1,...,l$
is equal to 1 and if
$\rho' \in \Omega^{r}_{q}$ 
is arbitrary, then the following formula is true:
\be
p_{k}(\rho_{1} \wedge \cdots \wedge \rho_{l} \wedge \rho') =
(\pi^{r+1,r})^{*}\rho_{1} \wedge \cdots \wedge (\pi^{r+1,r})^{*}\rho_{l}
\wedge p_{k-l}\rho', \quad \forall k \geq l.
\ee
If 
$k < l$
then the right hand side is zero, according to the convention from the
preceding lemma.
\end{cor}

We still need two results from \cite{Kr6}. The first one is elementary.

\begin{lemma}
Let $q \geq 1$ and
$\rho \in \Omega^{r}_{q}$.
Then in the associated chart
$(V^{r+1},\psi^{r+1})$
the following formula is valid:
\be
(\pi^{r+1,r})^{*}\rho = \sum_{s=0}^{q} {1 \over s!(q-s)!} 
\sum_{|I_{1}|,...,|I_{s}| \leq r} 
B^{I_{1},...,I_{s}}_{\sigma_{1},...,\sigma_{s},i_{s+1},...,i_{q}}
\omega^{\sigma_{1}}_{I_{1}} \wedge \cdots \wedge \omega^{\sigma_{s}}_{I_{s}} 
\wedge dx^{i_{s+1}} \cdots \wedge dx^{i_{q}}
\ee
where the coefficients
$B^{I_{1},...,I_{s}}_{\sigma_{1},...,\sigma_{s},i_{s+1},...,i_{q}}$
are smooth functions on the chart $V^{r}$ and verify symmetry properties of the
type (\ref{symmetry}). Moreover, the form
$p_{k}\rho$
is given by the terms corresponding to $s = k$ in the sum above.
\end{lemma}

Next, we have

\begin{lemma}
Let $q \geq 1$ and
$\rho \in \Omega^{r}_{q}$.
Suppose that in the associated chart
$(V^{r},\psi^{r})$
the form $\rho$ has the generic expression:
\be
\rho = \sum_{s=0}^{q} {1 \over s!(q-s)!} 
\sum_{|I_{1}|,...,|I_{s}| \leq r} 
A^{I_{1},...,I_{s}}_{\sigma_{1},...,\sigma_{s},i_{s+1},...,i_{q}}
dy^{\sigma_{1}}_{I_{1}}\wedge \cdots \wedge dy^{\sigma_{s}}_{I_{s}} \wedge
dx^{i_{s+1}} \cdots \wedge dx^{i_{q}}
\ee
where 
$A^{I_{1},...,I_{s}}_{\sigma_{1},...,\sigma_{s},i_{s+1},...,i_{q}}$
are smooth functions on $V^{r}$ verifying the symmetry property
(\ref{symmetry}). Then on the chart
$(V^{r+1},\psi^{r+1})$
we have
\be
p_{k}\rho = {1 \over k!(q-k)!} \sum_{|I_{1}|,...,|I_{k}| \leq r} 
B^{I_{1},...,I_{k}}_{\sigma_{1},...,\sigma_{k},i_{k+1},...,i_{q}}
\omega^{\sigma_{1}}_{I_{1}} \wedge \cdots \wedge \omega^{\sigma_{k}}_{I_{k}} 
\wedge dx^{i_{k+1}} \cdots \wedge dx^{i_{q}}
\ee
where
\be
B^{I_{1},...,I_{k}}_{\sigma_{1},...,\sigma_{k},i_{k+1},...,i_{q}} =
{\cal S}^{-}_{i_{k+1},...,i_{q}} \sum_{s=k}^{q} {q-k\choose q-s}
\sum_{|I_{k+1}|,...,|I_{s}| \leq r} 
A^{I_{1},...,I_{s}}_{\sigma_{1},...,\sigma_{s},i_{s+1},...,i_{q}}
y^{\sigma_{k+1}}_{I_{k+1}i_{k+1}}\cdots y^{\sigma_{s}}_{I_{s}i_{s}}.
\ee
\label{p-k-ro}
\end{lemma}

{\bf Proof:}
We use the definition (\ref{o'}) to write
\begin{eqnarray}
(\pi^{r+1,r})^{*}\rho = \sum_{s=0}^{q} {1 \over s!(q-s)!} 
\sum_{|I_{1}|,...,|I_{s}| \leq r} 
A^{I_{1},...,I_{s}}_{\sigma_{1},...,\sigma_{s},i_{s+1},...,i_{q}} \nonumber \\
(\omega^{\sigma_{1}}_{I_{1}} + y^{\sigma_{1}}_{I_{1}i_{1}} dx^{i_{1}}) \wedge 
\cdots \wedge 
(\omega^{\sigma_{s}}_{I_{s}} + y^{\sigma_{s}}_{I_{s}i_{s}} dx^{i_{s}})
\wedge dx^{i_{s+1}} \cdots \wedge dx^{i_{q}}
\end{eqnarray}
and now we can apply lemma \ref{decomp} with
$L \mapsto \Omega^{r}, \quad M \mapsto \Omega^{r}_{s}$
and
$\omega \mapsto \Lambda$
where
$$
\Lambda(\omega^{\sigma_{1}}_{I_{1}},...,\omega^{\sigma_{s}}_{I_{s}}) =
\sum_{|I_{1}|,...,|I_{s}| \leq r} 
A^{I_{1},...,I_{s}}_{\sigma_{1},...,\sigma_{s},i_{s+1},...,i_{q}}
\omega^{\sigma_{1}}_{I_{1}} \cdots \wedge \omega^{\sigma_{s}}_{I_{s}}.
$$
Then simple rearrangements leads to the formula from the statement.
$\qed$

Now we can give the structure theorem for strongly contact forms. As in the
preceding section, the proof will be based on Fock space machinery and will
differ from the original proof from \cite{Kr4}.

\begin{thm}
Let 
$n+1 \leq q \leq N$ 
and 
$\rho \in \Omega^{r}_{q}$. 
Let 
$(V,\psi)$
be a chart on $Y$ 
Then $\rho$ is a strongly contact form iff it has the following expression in 
the associated chart
$(V^{r},\psi^{r})$:
\be
\rho  = \sum_{p+s=q-n+1} \sum_{|J_{1}|,...,|J_{p}| \leq r-1}
\sum_{|I_{1}|,...,|I_{s}|=r-1} 
\omega^{\sigma_{1}}_{J_{1}} \cdots \wedge \omega^{\sigma_{p}}_{J_{p}} \wedge
d\omega^{\nu_{1}}_{I_{1}} \cdots \wedge d\omega^{\nu_{s}}_{I_{s}} 
\wedge \Phi^{J_{1},...,J_{p},I_{1},...,I_{s}}_{\sigma_{1},...,\sigma_{p},
\nu_{1},...,\nu_{s}}
\label{str-contact-q}
\ee
where
$\Phi^{J_{1},...,J_{p},I_{1},...,I_{s}}_{\sigma_{1},...,\sigma_{p},
\nu_{1},...,\nu_{s}}
$
are differential forms of degree $n-1-s$ on $V^{r}$. (This imposes 
that the first sum runs in fact only for $s \leq n-1$).
\label{str-contact}
\end{thm}

{\bf Proof:}
If $\rho$ has the expression from the statement, the corollary above gives us
$p_{q-n}\rho = 0$. We prove the converse statement by induction on $q$.

(i) Let $q = n+1$ and 
$\rho \in \Omega^{r}_{n+1}$
such that
$p_{1}\rho = 0$.
We start from the same decomposition of the form $\rho$ as in theorem
\ref{contact} i.e. that given by the formul\ae\/ (\ref{ro})-(\ref{canonical}). 
Using (\ref{ro}) and the corollary above we have from the preceding equation:
\be
\sum_{|J| \leq r-1} \omega^{\sigma}_{J} \wedge h\Phi^{J}_{\sigma} +
p_{1}\rho' = 0.
\label{1}
\ee

But the preceding lemma gives us the following very explicit formula for the
second contribution:
$$
p_{1}\rho' = {1\over (q-1)!} \sum_{|I|=r} B^{I}_{\sigma,i_{2},...,i_{q}}
\omega^{\sigma}_{I} \wedge dx^{i_{2}} \wedge \cdots \wedge dx^{i_{q}}
$$
where
$$
B^{I_{1}}_{\sigma_{1},i_{2},...,i_{q}} = {\cal S}^{-}_{i_{2},...,i_{q}} 
\sum_{s=1}^{q} {n \choose n+1-q} \sum_{|I_{2}|,...,|I_{s}|=r-1} 
A^{I_{1},...,I_{s}}_{\sigma_{1},...,\sigma_{s},i_{s+1},...,i_{q}}
y^{\sigma_{2}}_{I_{2}i_{2}}\cdots y^{\sigma_{s}}_{I_{s}i_{s}}.
$$

Because the 1-form
$\omega^{\sigma}_{I}$
appears only in the second term of (\ref{1}) the two terms must vanish
separately i.e. we have
\be
p_{1}\rho' = 0.
\label{1a}
\ee
and
\be
\sum_{|J| \leq r-1} \omega^{\sigma}_{J} \wedge h\Phi^{J}_{\sigma} = 0.
\label{1b}
\ee

From (\ref{1a}) we get
$$
B^{I}_{\sigma,i_{2},...,i_{q}} = 0
$$
which can be transformed, as in the proof of theorem \ref{contact} into
$$
B_{2}\cdots B_{s} A_{\sigma_{1},...,\sigma_{s}} = 0 \quad (s = 1,...,q).
$$

In fact, because of the symmetry property (\ref{symmetry}) we have:
\be
B_{1}\cdots \hat{B}_{\alpha} \cdots B_{s} A_{\sigma_{1},...,\sigma_{s}} = 0 
\quad (\alpha = 1,...,s; s = 1,...,q). 
\label{p1}
\ee

A consequence of this relation is
$$
B_{1}\cdots B_{s} A_{\sigma_{1},...,\sigma_{s}} = 0 \quad (s = 1,...,q)
$$
which implies (see lemmas \ref{s=q} and \ref{s<q})) that
\be
A_{\sigma_{1},...,\sigma_{q}} = 0 
\ee
and
\be
A_{\sigma_{1},...,\sigma_{s}} = \sum_{\alpha =1}^{s} B_{\alpha}
A^{\alpha}_{\sigma_{1},...,\sigma_{s}} \quad (s = 1,...,q-1)
\ee
for some tensors
$A^{\alpha}_{\sigma_{1},...,\sigma_{s}}$.

If we substitute the expression of 
$A_{\sigma_{1},...,\sigma_{s}}$
above into the initial equation (\ref{p1}) we get immediately
$$
B_{1}\cdots B_{s} A^{\alpha}_{\sigma_{1},...,\sigma_{s}} = 0 
$$
so lemma \ref{s<q} can again be applied to produce the following expression 
\be
A_{\sigma_{1},...,\sigma_{s}} = \sum_{\alpha,\beta =1}^{s} B_{\alpha} B_{\beta}
A^{\alpha\beta}_{\sigma_{1},...,\sigma_{s}}  \quad (s = 1,...,q-1)
\ee
for some tensors
$A^{\alpha\beta}_{\sigma_{1},...,\sigma_{s}} = 0$.
The last expression identically verifies the equation (\ref{p1}) so it is the
general solution of it. As in the end of theorem \ref{contact} it this the 
time to revert to full index notations. Because in the formula above we have
{\it two} $B$-type operators we will obtain that the functions
$A^{I_{1},...,I_{s}}_{\sigma_{1},...,\sigma_{s},i_{s+1},...,i_{q}}$
are sums of terms containing {\it two} delta factors, so in the end {\it two}
factors of the type
$d\omega^{\nu}_{I}$
we show up. Explicitly, $\rho'$ must necessarily have the following structure:
\be
\rho' = \sum_{|I_{1}|=|I_{2}|=r-1} d\omega^{\nu_{1}}_{I_{1}} \wedge
d\omega^{\nu_{2}}_{I_{2}} \wedge \Phi^{I_{1}I_{2}}_{\nu_{1}\nu_{2}}.
\ee

On the other hand, it is easy to see that (\ref{1b}) is equivalent to
$$h\Phi^{J}_{\sigma} = 0$$
and theorem \ref{contact} can be applied. Combining with the formula above, we
obtain the structure formula from the statement (\ref{str-contact-q}) for 
$q = n + 1$. 

(ii) We suppose that if
$p_{q'-n}\rho = 0$
then $\rho'$ has the expression (\ref{str-contact-q}) for 
$q' = n + 1,...,q-1$
and we prove the same statement for $q$. So we have
$$p_{q-n}\rho = 0.$$

Because $\rho$ is a polynomial of degree $q$ (with respect to the wedge product
$\wedge$) in the differentials
$\omega^{\sigma}_{J} \quad |J| \leq r-1, \quad dy^{\sigma}_{I}, \quad |I| = r$
and
$dx^{i}$
one can write it uniquely as follows:
\be
\rho = \sum_{s=0}^{q} {1\over q!} \sum_{|J_{1}|,...,|J_{s}| \leq r-1} 
\omega^{\sigma_{1}}_{J_{1}} \cdots \wedge \omega^{\sigma_{s}}_{J_{s}} \wedge
\Phi^{J_{1},...J_{s}}_{\sigma_{1},...,\sigma_{s}}
\ee
where
$\Phi^{J_{1},...J_{s}}_{\sigma_{1},...,\sigma_{s}}$
are polynomials of degree $q - s$ in the differentials
$dy^{\sigma}_{I}, \quad |I| = r$
and
$dx^{i}$.
Using the corollary above one obtains the following equation:
$$
\sum_{s=0}^{q} {1\over q!} \sum_{|J_{1}|,...,|J_{s}| \leq r-1} 
\omega^{\sigma_{1}}_{J_{1}} \cdots \wedge \omega^{\sigma_{s}}_{J_{s}} \wedge
p_{q-n-s}\Phi^{J_{1},...J_{s}}_{\sigma_{1},...,\sigma_{s}} = 0
$$
which is equivalent to:
\be
p_{q-n-s}\Phi^{J_{1},...J_{s}}_{\sigma_{1},...,\sigma_{s}} = 0 
\quad (s = 0,...,q).
\ee

For $s = 1,...,q$ one can apply the induction hypothesis and obtain that the
forms 
$\Phi^{J_{1},...J_{s}}_{\sigma_{1},...,\sigma_{s}}$
are sums of the type (\ref{str-contact-q}). It remains to analyse the case 
$s = 0$ i.e. the equation
\be
p_{q-n}\Phi = 0
\ee
with $\Phi$ having a structure similar to (\ref{canonical}):
\be
\Phi = \sum_{s=q-n}^{q} {1 \over s!(q-s)!} \sum_{|I_{1}|,...,|I_{s}|=r}
A^{I_{1},...,I_{s}}_{\sigma_{1},...,\sigma_{s},i_{s+1},...,i_{q}}
dy^{\sigma_{1}}_{I_{1}} \wedge \cdots \wedge dy^{\sigma_{s}}_{I_{s}} \wedge
dx^{i_{s+1}} \wedge \cdots \wedge dx^{i_{q}}.
\label{q}
\ee
 
Using the preceding lemma one obtains that
$$
p_{q-n}\Phi =  {1 \over (q-n)!n!} \sum_{|I_{1}|,...,|I_{q-n}|=r} 
B^{I_{1},...,I_{q-n}}_{\sigma_{1},...,\sigma_{q-n},i_{q-n+1},...,i_{q}}
\omega^{\sigma_{1}}_{I_{1}} \wedge \cdots \wedge 
\omega^{\sigma_{q-n}}_{I_{q-n}} 
\wedge dx^{i_{q-n+1}} \cdots \wedge dx^{i_{q}}
$$
where
\be
B^{I_{1},...,I_{q-n}}_{\sigma_{1},...,\sigma_{q-n},i_{q-n+1},...,i_{q}} =
{\cal S}^{-}_{i_{q-n+1},...,i_{q}} \sum_{s=q-n}^{q} {n\choose q-s}
\sum_{|I_{q-n+1}|,...,|I_{s}|=r} 
A^{I_{1},...,I_{s}}_{\sigma_{1},...,\sigma_{s},i_{s+1},...,i_{q}}
y^{\sigma_{q-n+1}}_{I_{q-n+1}i_{q-n+1}}\cdots y^{\sigma_{s}}_{I_{s}i_{s}}.
\ee

The condition on $\Phi$ translates into
$$
B^{I_{1},...,I_{q-n}}_{\sigma_{1},...,\sigma_{q-n},i_{q-n+1},...,i_{q}} = 0
$$
and this can be shown to be equivalent to 
$$
B_{q-n+1}\cdots B_{s} A_{\sigma_{1},...,\sigma_{s}} = 0 \quad (s = q-n,...,q).
$$

In fact, because of the symmetry property (\ref{symmetry}) we have more
generally: 
\be
B_{\alpha_{1}}\cdots B_{\alpha_{s-q+n}} A_{\sigma_{1},...,\sigma_{s}} = 0 
\quad \forall \alpha_{1},...,\alpha_{s-q+n}, (s = q-n,...,q).
\label{pq}
\ee

This relation can be investigated following the ideas from (i) (see rel. 
(\ref{p1})) and the general solution of (\ref{pq}) can be found in the form:
\be
A_{\sigma_{1},...,\sigma_{s}} = \sum_{\alpha_{1},...,\alpha_{q-n+1}=1}^{s}
B_{\alpha_{1}} \cdots B_{\alpha_{q-n+1}} 
A^{\alpha_{1},...,\alpha_{q-n+1}}_{\sigma_{1},...,\sigma_{s}} 
\ee
for some tensors
$A^{\alpha_{1},...,\alpha_{q-n+1}}_{\sigma_{1},...,\sigma_{s}}$.

If we use full index notations, this time $q-n+1$ factors of the type
$d\omega^{\nu}_{I}, \quad |I|=r-1$
will appear in every term of $\Phi$. Collecting all terms we get for
$\rho$ the formula (\ref{str-contact-q}).
$\qed$

\newpage

\section{Euler-Lagrange and Helmholtz-Sonin forms}

An interesting problem in differential geometry is the following one. Suppose
we have a differential form $\rho$ on a given manifold $Y$. What other 
(globally defined) differential forms can be attached to it? This problem can 
be rigorously formulated \cite {Kr8} and the answer is that there is 
essentially only one possibility, namely the exterior differential $d\rho$ 
of $\rho$. In other words the condition of correct behaviour with respect to
all possible charts transformations limits drastically the possible solution to
this kind of problem. But what happens when the manifold $Y$ has a
supplementary structure, say is a fibre bundle? Then, there will be some
restrictions on the charts transformation so other solutions can appear. In
this section we will prove that in this case indeed new possibilities can
appear, as for instance the Euler-Lagrange and Helmholtz-Sonin form. We will
follow essentially \cite{An} making the observation that much of the line of
the argument can be adapted from infinite jet bundle extensions to our case
i.e. finite bundle extensions.

\subsection{Lie-Euler Operators}

The central combinatorial trick used in \cite{An} to prove the existence of the
Euler-Lagrange form is the concept of {\it total differential operator} which,
by definition, is any linear map
$P: {\cal E}(J^{r}Y) \rightarrow \Omega^{s}$
covering the identity map
$id: J^{r}Y \rightarrow J^{r}Y$ with $s \geq r$.
(One considers, of course, ${\cal E}(J^{r}Y)$ and $J^{r}Y$ as fibre bundles 
over
$J^{r}Y$). We will consider in the following that $s$ is sufficiently great; 
in fact one needs that $s > 2r + 2$. Suppose that $\xi$ is an evolution having
the local structure (\ref{ev}) in the chart 
$(V^{r},\psi^{r})$.
Then the image $P(\xi) \in \Omega^{s}$ must have the expression:
\be
P(\xi) = \sum_{|I| \leq r} (d_{I}\xi^{\sigma}) P^{I}_{\sigma} =
\sum_{k=0}^{r} (d_{j_{1}}\cdots d_{j_{k}}\xi^{\sigma}) P^{j_{1},...,j_{k}}_{I},
\label{t-diff-op}
\ee
where $P^{I}_{\sigma}$ are (local) differential forms in the chart
$(V^{s},\psi^{s})$
and, as usual, $d_{j} = d_{j}^{s}$ (see (\ref{formal})).

Then one has the following combinatorial lemma \cite{An}:

\begin{lemma}
In the conditions above, the following formula is true:
\be
P(\xi) = \sum_{|I| \leq r} d_{I} (\xi^{\sigma} Q^{I}_{\sigma})
\label{PQ}
\ee
where 
\be
Q^{I}_{\sigma} \equiv \sum_{|J| \leq r-|I|} (-1)^{|J|} 
{|I|+|J| \choose |I|} d_{J} P^{IJ}_{\sigma}
\label{Q}
\ee
and one assumes that the action of a formal derivative $d_{j}$ on a form is
realized by its action on the function coefficients.
\label{trick}
\end{lemma}

{\bf Proof:}
One starts from the right hand side of (\ref{PQ}) and uses Leibnitz rule:
$$
\sum_{|I| \leq r} d_{I} (\xi^{\sigma} Q^{I}_{\sigma}) =
\sum_{|I| \leq r} \sum_{(J,K)} (d_{J} \xi^{\sigma}) (d_{K} Q^{I}_{\sigma})
$$
where the sum over $(J,K)$ is over all partitions of the set $I$. One can
rearrange this as follows:
\begin{eqnarray}
\sum_{|I| \leq r} d_{I} (\xi^{\sigma} Q^{I}_{\sigma}) =
\sum_{|J|+|K| \leq r} {|J|+|K| \choose |J|} (d_{J} \xi^{\sigma})
(d_{K} Q^{JK}_{\sigma}) =\nonumber \\
= \sum_{|J| \leq r} (d_{J} \xi^{\sigma}) \sum_{|K| \leq r-|J|} 
{|J|+|K| \choose |J|} d_{K} Q^{JK}_{\sigma}.\nonumber
\end{eqnarray}

Now one proves by elementary computations that
\be
P^{I}_{\sigma} = \sum_{|J| \leq r-|I|} {|J|+|K| \choose |J|} 
d_{J} \xi^{\sigma} Q^{IJ}_{\sigma}.
\label{IJ}
\ee
and that finishes the proof.
$\qed$

\begin{rem}
One notices that the relation  (\ref{PQ}) {\bf uniquely} determines the forms
$Q^{I}_{\sigma}$.
\label{unicity}
\end{rem}

We proceed now to formulate the main result of this subsection. The proof is
an easy adaptation to the finite jet bundle extension case of the proof from
\cite{An}. 

\begin{thm}
Let $q \geq n$ and 
$P: {\cal E}(J^{r}Y) \rightarrow \Omega^{s}_{q,X}$
a total differential operator. Let 
$(V,\psi)$
and
$(\bar{V},\bar{\psi})$
two overlapping charts on Y and let us construct on the intersection of the 
corresponding associated charts the forms 
$Q^{I}_{\sigma}$ and $\bar{Q}^{I}_{\sigma}$ according to the preceding lemma.
Then the following relation is true on the intersection
$V^{s} \cap \bar{V}^{s}$:
\be
Q_{\sigma} = (\partial_{\sigma} \bar{y}^{\nu}) \bar{Q}_{\nu}.
\label{Q-Q-bar}
\ee

In particular, there exists a globally defined form, denoted by
$E(P)(\xi)$
such that in the chart
$(V^{r},\psi^{r})$
we have
\be
E(P)(\xi) = Q_{\sigma} \xi^{\sigma}.
\label{E(P)}
\ee
\end{thm}

{\bf Proof:}
The generic expression for $Q^{I}_{\sigma}$ in the chart $V^{s}$ is (see
(\ref{hor2})):
\be
Q^{I}_{\sigma} = \sum_{|J_{1}|,...,|J_{l}| \leq s}   
Q^{I,J_{1},...,J_{l}}_{\sigma,\nu_{1},...,\nu_{l}} dy^{\nu_{1}}_{J_{1}} 
\wedge \cdots \wedge dy^{\nu_{l}}_{J_{l}} \wedge \theta_{0}.
\label{V}
\ee

Here $l = q-n$ and 
$Q^{I,J_{1},...,J_{l}}_{\sigma,\nu_{1},...,\nu_{l}}$
are smooth functions on $V^{s}$ having appropriate antisymmetry properties.

In the other chart $\bar{V}^{s}$ we have a similar expression:

\begin{eqnarray}
\bar{Q}^{I}_{\sigma} = \sum_{|J_{1}|,...,|J_{l}| \leq s}   
\bar{Q}^{I,J_{1},...,J_{l}}_{\sigma,\nu_{1},...,\nu_{l}} 
d\bar{y}^{\nu_{1}}_{J_{1}} \wedge \cdots \wedge d\bar{y}^{\nu_{l}}_{J_{l}} 
\wedge \bar{\theta}_{0} = \nonumber \\
= {\cal J} \sum_{|J_{1}|,...,|J_{l}| \leq s}   
\tilde{Q}^{I,J_{1},...,J_{l}}_{\sigma,\nu_{1},...,\nu_{l}} 
dy^{\nu_{1}}_{J_{1}} \wedge \cdots \wedge dy^{\nu_{l}}_{J_{l}} 
\wedge \theta_{0}
\label{bar-V}
\end{eqnarray}
where ${\cal J}$ is the Jacobian of the chart transformation on $X$ (see 
(\ref{J})).

If we define for any $I$ with $|I| \leq r-1$ 
$$
R^{I}_{\sigma} \equiv i_{d_{j}} Q^{jI}_{\sigma}
$$
then, using the formula above, it is easy to prove that:
$$
Q^{j_{1},...,j_{k}}_{\sigma} = {\cal S}^{+}_{j_{1},...,j_{k}} dx^{j_{1}} \wedge
R^{j_{2},...,j_{k}}_{\sigma} = {1 \over k} \sum_{p=1}^{k} dx^{j_{p}} \wedge
R^{j_{1},...,\hat{j}_{p},...,j_{k}}_{\sigma}.\quad k = 0,...,r.
$$

As a consequence one can rewrite the local formula (\ref{PQ}) as follows:
\be
P(\xi) = \xi^{\sigma} Q_{\sigma} + dx^{i} \wedge d_{i} R(\xi),
\ee
where we have defined
\be
R(\xi) \equiv \sum_{|I| \leq r-1} d_{I} (\xi^{\sigma} R^{I}_{\sigma}).
\ee

So, in the overlap 
$V^{s} \cap \bar{V}^{s}$
we have
$$
\xi^{\sigma} Q_{\sigma} - \bar{\xi}^{\sigma} \bar{Q}_{\sigma} =
d\bar{x}^{i} \wedge \bar{d}_{i} \bar{R}(\bar{\xi}) - 
dx^{i} \wedge d_{i} R(\xi).
$$

Let us remark now that from the definition of the forms
$Q^{j_{1},...,j_{k}}_{\sigma}$
(see (\ref{Q})) it follows that its function coefficients depend only on the
variables 
$(x^{i},y^{\sigma},y^{\sigma}_{j},...,y^{\sigma}_{j_{1},...,j_{2r-k}})$;
as a consequence, the function coefficients of the form 
$R(\xi)$
depend only on the variables
$(x^{i},y^{\sigma},y^{\sigma}_{j},...,y^{\sigma}_{j_{1},...,j_{2r-1}})$.
In this case we can apply formula (\ref{d-bar-d}) to the preceding relation 
and we obtain:
$$
\xi^{\sigma} Q_{\sigma} - \bar{\xi}^{\sigma} \bar{Q}_{\sigma} =
dx^{i} \wedge d_{i} \tilde{R}(\xi)
$$
where
$$
\tilde{R}(\xi) = \bar{R}(\bar{\xi}) - R(\xi) =
\sum_{|J_{1}|,...,|J_{l}| \leq s}   
\tilde{R}^{i,J_{1},...,J_{l}}_{\nu_{1},...,\nu_{l}}(\xi) 
dy^{\nu_{1}}_{J_{1}} \wedge \cdots \wedge dy^{\nu_{l}}_{J_{l}} \wedge 
\theta_{i};
$$
here we have defined
\be
\theta_{i} \equiv (-1)^{i-1} x^{1} \wedge \cdots \wedge dx^{i-1} \wedge
dx^{i+1} \cdots \wedge dx^{n}
\label{theta-i}
\ee
and
$\tilde{R}^{i,J_{1},...,J_{l}}_{\nu_{1},...,\nu_{l}}(\xi)$
are smooth functions on the overlap
$V^{s} \cap \bar{V}^{s}$.

If we also use (\ref{ev}) we obtain 
$$
\xi^{\sigma} \left[ Q^{\emptyset,J_{1},...,J_{l}}_{\sigma,\nu_{1},...,\nu_{l}} 
- {\cal J} (\partial_{\sigma} \bar{y}^{\zeta}) \bar{\xi}^{\sigma} 
\tilde{Q}^{\emptyset,J_{1},...,J_{l}}_{\zeta,\nu_{1},...,\nu_{l}}\right]  =
d_{i} \tilde{R}^{i,J_{1},...,J_{l}}_{\nu_{1},...,\nu_{l}}(\xi).
$$ 

Now one proves that both sides are zero in a standard way: one picks a section
$\gamma$ with support in 
$W \subset \pi(V) \cap \pi(\bar{V})$
such that the closure $\bar{W}$ of $W$ is compact, takes $\xi^{\sigma}$ with
support in the open set 
$U \subset \bar{W}$
and integrates on 
$\bar{W}$
the following relation (which follows from the preceding one):
$$
\xi^{\sigma} \left[ Q^{\emptyset,J_{1},...,J_{l}}_{\sigma,\nu_{1},...,\nu_{l}}
- {\cal J} (\partial_{\sigma} \bar{y}^{\zeta}) \bar{\xi}^{\sigma} 
\tilde{Q}^{\emptyset,J_{1},...,J_{l}}_{\zeta,\nu_{1},...,\nu_{l}}\right] \circ
j^{s}\gamma =
d_{i} \tilde{R}^{i,J_{1},...,J_{l}}_{\nu_{1},...,\nu_{l}}(\xi) \circ 
j^{s}\gamma .
$$

Use of Stokes theorem is made and of the arbitrariness of $\gamma$ and it 
follows that:
$$
Q^{\emptyset,J_{1},...,J_{l}}_{\sigma,\nu_{1},...,\nu_{l}} =
{\cal J} (\partial_{\sigma} \bar{y}^{\zeta}) \bar{\xi}^{\sigma} 
\tilde{Q}^{\emptyset,J_{1},...,J_{l}}_{\zeta,\nu_{1},...,\nu_{l}}
$$

If we introduce this equality in (\ref{V}) and (\ref{bar-V}) we obtain the 
relation (\ref{Q-Q-bar}).
$\qed$

The operator $E(P)$ defined by (\ref{E(P)}) is called the {\it Euler operator}
associated to the total differential operator $P$; it has the local expression:
\be
E(P)(\xi) = \xi^{\sigma} E_{\sigma}(P)
\ee
where
\be
E_{\sigma}(P) = \sum_{|I|=0}^{r} (-1)^{|I|} d_{I} P^{I}_{\sigma}.
\ee

Now one takes 
$\lambda \in \Omega^{r}_{n,X}$
and constructs the total differential operator
$P_{\lambda}$:
\be
P_{\lambda}(\xi) \equiv L_{pr(\xi}) \lambda.
\ee

Suppose that $\lambda$ has the local expression (\ref{L}). 
Then lemma \ref{trick} can be applied and gives the following formula:
\be
P_{\lambda}(\xi) = \sum_{|I|=0}^{r} d_{I} \left( \xi^{\sigma}
E^{I}_{\sigma}(L)\right) \theta_{0}
\ee
where
\be
E^{I}_{\sigma}(L) \equiv \sum_{|J|\leq r-|I|} (-1)^{|J|} {|I|+|J| \choose |I|}
d_{J} \partial^{IJ}_{\sigma} L
\label{Lie-Euler}
\ee
are the so-called {\it Lie-Euler operators}.

In particular, 
\be
Q_{\sigma} = E_{\sigma}(L) \theta_{0}
\label{QL}
\ee
and the Euler operator associated to $P_{\lambda}$ has the following
expression: 
\be
E(P_{\lambda}) = \xi^{\sigma} E_{\sigma}(L) \theta_{0}
\ee
where
\be
E_{\sigma}(L) \equiv \sum_{|J|\leq r} (-1)^{|J|} d_{J} \partial^{J}_{\sigma} L
\label{EL1}
\ee
are the {\it Euler-Lagrange expressions}.

The theorem above leads to
\begin{prop}
If
$\lambda \in \Omega^{r}_{n,X}$
is a Lagrange form, then there exists a globally defined $n + 1$-form, 
denoted by
$E(\lambda)$
such that we have in the chart $V^{s}$:
\be
E(\lambda) = E_{\sigma}(L) \omega^{\sigma} \wedge \theta_{0}.
\label{EL2}
\ee
\end{prop}

{\bf Proof:}
Indeed, one combines the expression (\ref{QL}) with the transformation 
properties (\ref{o-bar-o}) and (\ref{Q-Q-bar}) to obtain 
\be
E_{\sigma}(L) = {\cal J} (\partial_{\sigma}\bar{y}^{\nu}) 
\bar{E}_{\nu}(\bar{L});
\label{E-bar-E}
\ee
the globallity of 
$E(\lambda)$ follows immediately. 
$\qed$

One calls this form the {\it Euler-Lagrange form} associated to $\lambda$
and notes that it is a differential equation (see the beginning of subsection
3.3 more precisely the formula (\ref{T})). In general, a differential equation 
$T \in \Omega^{s}_{n+1,X}$ 
is called {\it (locally) variational} or {\it of Euler-Lagrange type} 
{\it iff} there exists a (local) Lagrange form 
$\lambda \in \Omega^{r}_{n,X} \quad (s \geq 2r)$ 
such that 
\be
T = E(\lambda).
\label{EL-type}
\ee

One notices that in this case the general form of a differential equation
(\ref{eq-diff}) coincides with the well-known form of the Euler-Lagrange
equations. 

\begin{rem}
There are other ways of proving the globallity of the Euler-Lagrange form.
One can use the existence of the so-called Lepagean equivalents \cite{Kr2},
but the combinatorial analysis seems to be more complicated. Also, an argument
based on the connection between the action functional and the Euler-Lagrange
expression is available \cite{Kr4}
\end{rem}


\subsection{Some Properties of the Euler-Lagrange Form}

We collect for further use some properties of the Euler-Lagrange form. We
follow, essentially \cite{An}, with the appropriate modifications.

First, one finds rather easily from (\ref{E-bar-E}) that the Euler-Lagrange
form behaves naturally with respect to bundle morphisms. More precisely, if
$\phi \in Diff(Y)$
is a bundle morphism, then one has:
\be
(j^{s}\phi)^{*} E(\lambda) = E((j^{r}\phi)^{*} \lambda)
\ee
for any Lagrange form $\lambda$. From here, we obtain by differentiation:
\be
L_{j^{s}\xi} E(\lambda) = E(L_{j^{r}\xi} \lambda)
\label{Lie-EL}
\ee
for any projectable vector field $\xi$ on $Y$.

Now we have 
\begin{lemma}
If $f$ is a smooth function on $V^{r}$ then we have in 
$V^{s}, \quad s > 2r + 2$
the following formul\ae:
\be
E^{Ij}_{\sigma}(d_{l}f) = {\cal S}^{+}_{Ij} \delta^{j}_{l} E^{I}_{\sigma}(f),
\quad |I| = 0,...,r
\ee
and
\be
E_{\sigma}(d_{l}f) = 0.
\ee
\end{lemma}

{\bf Proof:}
From lemma \ref{trick} we have in $V^{s}$:
$$
\sum_{|I| \leq r} (d_{I}\xi^{\sigma}) (\partial^{I}_{\sigma}L) =
\sum_{|I| \leq r} d_{I} \left( \xi^{\sigma} E^{I}_{\sigma}(L)\right)
$$
for any smooth function $L$ on $V^{r}$. We make 
$r \rightarrow r + 1$ and $L \rightarrow d_{j}f$
and we have:
$$
\sum_{|I| \leq r+1} (d_{I}\xi^{\sigma}) 
\left(\partial^{I}_{\sigma}d_{j}f\right)
 = \sum_{|I| \leq r+1} d_{I} \left( \xi^{\sigma} E^{I}_{\sigma}(d_{j}f)\right).
$$

The left hand side can be rewritten using the commutation formula 
(\ref{commutator}) and one obtains:
$$
\sum_{k=1}^{r+1} d_{i_{1}}...d_{i_{k}}\left( \xi^{\sigma} 
{\cal S}^{+}_{i_{1},...,i_{k}} \delta^{i_{1}}_{j} 
E^{i_{2},...,i_{k}}_{\sigma}(f) \right)  = 
\sum_{k=0}^{r+1} d_{i_{1}}...d_{i_{k}} \left( \xi^{\sigma} 
E^{i_{1},...,i_{k}}_{\sigma}(d_{j}f)\right).
$$

Using remark \ref{unicity} we obtain the relations from the statement.
$\qed$

\begin{cor}
Let $A^{I}, \quad |I| = l \geq 1$ be some smooth functions on $V^{r}$ and
$f \equiv d_{I} A^{I}$. Then we have on $V^{s}, \quad s > 2(r + l)$:
\be
E^{J}_{\sigma}(f) = 0, \quad |J| \leq |I| -1.
\ee
\label{tot-div}
\end{cor}

\begin{cor}
Let $\xi$ be an evolution on $Y$ and $\lambda \in \Omega^{r}_{n,X}$ a Lagrange
form. Then the following formula is true:
\be
E\left( L_{j^{r}\xi}\lambda\right) = E\left( i_{j^{s}\xi} E(\lambda)\right).
\ee
\end{cor}

{\bf Proof:}
One has by direct computation and use on lemma \ref{trick}:
$$
E\left( L_{j^{r}\xi}\lambda - i_{j^{s}\xi} E(\lambda)\right) =
\sum_{|I|=1}^{r} E_{\nu}\left( d_{I}(\xi^{\sigma} E^{I}_{\sigma}(L))\right)
\omega^{\nu} \wedge \theta_{0};
$$
but the right hand side is zero, according to the preceding corollary.
$\qed$

\begin{cor}
Let $\lambda \in \Omega^{r}_{n,X}$ a Lagrange form. Then the following formula
is true:
\be
L_{j^{s}\xi} E(\lambda) = E\left( i_{j^{s}\xi} E(\lambda)\right).
\label{E-pr-xi}
\ee
\end{cor}

Indeed, one combines the preceding corollary with (\ref{Lie-EL}) and obtains
this formula.
\newpage

\subsection{Helmholtz-Sonin Forms}

In this section, we follow an idea of \cite{AP2} to prove the existence of the
(globally) defined Helmholtz-Sonin form. We have the following central result.

\begin{thm}
Let $T \in \Omega^{s}_{n+1,X}$ be a differential equation with the local form
given by (\ref{T}). We define the following expressions in any chart $V^{t},
\quad t > 2s$:
\be
H^{J}_{\sigma\nu} \equiv \partial^{J}_{\nu} T_{\sigma} - (-1)^{|J|}
E^{J}_{\sigma}(T_{\nu}), \quad |J| \leq s.
\label{HS-local}
\ee

Then there exists a globally defined $(n+2)$-form, denoted by $H(T)$ such that 
in any chart $V^{t}$ we have:
\be
H(T) = \sum_{|J| \leq s} H^{J}_{\sigma\nu} \omega^{\nu}_{J} \wedge
\omega^{\sigma} \wedge \theta_{0}.
\label{HS}
\ee
\end{thm}

{\bf Proof:}
(i) We begin with a construction from \cite{AP2}. Let $\xi$ be an evolution; we
define a (global) $n$-form $H_{\xi}(T)$ according to:
\be
H_{\xi}(T) \equiv L_{j^{s}\xi} T - E\left( i_{j^{s}\xi} T\right)
\label{H-xi}
\ee

Elementary computations and use of corollary \ref{tot-div} leads to the 
following local expression:
\be
H_{\xi}(T) = \sum_{|I| \leq s} (d_{I} \xi^{\sigma}) H^{I}_{\sigma\nu}
\omega^{\sigma} \wedge \theta_{0}.
\ee

(ii) Now one determines the transformation formula at a change of charts for
the expressions 
$d_{I} \xi^{\sigma}$. One considers the evolution 
$\xi$
on 
$Y$
and writes the expression of the vector field
$pr(\xi)$
on the overlap
$V^{t} \cap \bar{V}^{t}$;
the following formula easily emerges:
\be
\bar{d}_{I} \bar{\xi}^{\sigma} = \sum_{|J| \leq |I|} \left(\partial^{J}_{\nu} 
\bar{y}^{\sigma}_{I} \right) (d_{J} \xi^{\nu}), \quad |I| = 0,...,s.
\label{d-xi}
\ee
 
(iii) Using the transformation formula (\ref{d-xi}) one can obtain the
transformation formula for the expressions $H^{I}_{\sigma\nu}$: one has in
the overlap $V^{t} \cap \bar{V}^{t}$:
\be
H^{J}_{\mu\zeta} = {\cal J} \sum_{|I| \geq |J|} 
\left( \partial^{J}_{\mu} \bar{y}^{\nu}_{I} \right) 
\left( \partial_{\zeta} \bar{y}^{\sigma} \right) \bar{H}^{I}_{\sigma\nu}.
\ee

This transformation formula should now be combined with (\ref{o-bar-o}) and one
obtains that $H(T)$ has an invariant meaning; we have
$$
\sum_{|I| \leq s} \bar{H}^{I}_{\sigma\nu} \bar{\omega}^{\nu}_{I} \wedge
\bar{\omega}^{\sigma} \wedge \bar{\theta}_{0} =
\sum_{|J| \leq s} H^{J}_{\mu\zeta} \omega^{\mu}_{J} \wedge \omega^{\zeta} 
\wedge \theta_{0}
$$
on $V^{t} \cap \bar{V}^{t}$ and the proof is finished.
$\qed$
 
$H(T)$ is called the {\it Helmholtz-Sonin form} associated to $T$ and 
$H^{I}_{\sigma\nu}$ 
are the {\it Helmholtz-Sonin expressions} associated to $T$. 

A well-known corollary of the theorem above is:

\begin{cor}
The differential equation $T$ is locally variational {\it iff} $H(T) = 0$.
\end{cor}

{\bf Proof:}
(i) If there exists a Lagrange form $\lambda$ such that locally 
$T = E(\lambda)$, then from (\ref{E-pr-xi}) it follows that
$H_{\xi}(T) = 0$ for any evolution $\xi$. This implies that the Helmholtz-Sonin
expressions associated to $T$ are zero i.e. $H(E(\lambda)) = 0$.

(ii) The converse of this statement is done by some explicit construction.
Suppose that $T$ is such that $H(T) = 0$. Then we define the following
Lagrangian on the chart $V^{s}$:
\be
L = \int_{0}^{1} \quad dt y^{\sigma} T_{\sigma} \circ \chi^{s}_{t}, 
\qquad 
\chi^{s}_{t}(x^{i},y^{\sigma},y^{\sigma}_{j},...,y^{\sigma}_{j_{1},...,j_{s}})
= (x^{i},ty^{\sigma},ty^{\sigma}_{j},...,ty^{\sigma}_{j_{1},...,j_{s}}). 
\label{Tonti}
\ee

Then one proves by direct computation that $E_{\sigma}(L) = {\cal T}_{\sigma}$
so, $T$ is locally variational.
$\qed$

The (local) expression (\ref{Tonti}) is called the {\it Tonti-Vainberg} 
Lagrangian.
\newpage

\section{The Exact Variational Sequence}

This section is divided in two parts. The first one includes in a brief way
the standard proof of the exactness of the variational sequence following the
lines of \cite{Kr4}. The second part is devoted to the characterization of some
elements of the variational sequence by (globally) defined forms. As it is
pointed out in \cite{Kr4}, this can be done using the Euler-Lagrange and the
Helmholtz-Sonin forms defined previously. In this part more details are given
because the proofs from the literature are rather sketchy. In particular, the
proof we offer for the characterization of the $n + 1$ term in the variational
sequence by the Helmholtz-Sonin form is new.

\subsection{The Exactness of the Variational Sequence}

If $\pi: Y \rightarrow X$ is a fibre bundle and $U, V \subset Y$ are charts
such that $U \subset V$, we denote by 
$\i_{U,V}: U^{r} \rightarrow V^{r}$
the canonical inclusion. Then, the collection
$\{\Omega^{r}_{q}(V)\}\quad (q \geq 0), \quad \{\i^{*}_{U,V}\}$
is a presheaf denoted by $\Omega^{r}_{q}$.
Next, one introduces the subspaces 
$\theta^{r}_{q} \in \Omega^{r}_{q(c)}$ 
by:
\be
\theta^{r}_{1} \equiv \Omega^{r}_{1(c)}, \quad
\theta^{r}_{q} \equiv d\Omega^{r}_{q-1(c)} + \Omega^{r}_{q(c)} \quad 
(q= 2,...,N = dim(J^{r}Y));
\label{theta-q}
\ee
this is also a presheaf and one can easily verify that:
\be
\theta^{r}_{q} = \Omega^{r}_{q(c)} \quad (q = 2,...,n), \quad
\theta^{r}_{q} = 0 \quad (q > P ), \quad
d\theta^{r}_{q} \subset \theta^{r}_{q+1}; 
\quad P \equiv m {n+r-1 \choose n} + 2(n-1).
\ee

Next, one introduces the so-called contact homotopy operator. The
construction is the following. Let $U \subset \R^{n}$ (resp. 
$V \subset \R^{m}$) an open set (resp. a ball centred in $0 \in \R^{m}$) and
$W = U \times V$. One considers the operator $\chi_{r}$ as a map
$\chi^{r}_{t}: [0,1] \times J^{r}W \rightarrow J^{r}W$
given by:
\be
\chi_{r}(t,(x^{i},y^{\sigma},y^{\sigma}_{j},...,y^{\sigma}_{j_{1},...,j_{s}})) 
= (x^{i},ty^{\sigma},ty^{\sigma}_{j},...,ty^{\sigma}_{j_{1},...,j_{s}}). 
\label{chi}
\ee
Then for any
$\rho \in \Omega^{r}_{q}(W)$
we have the {\it unique} decomposition
\be
(\chi^{r}_{t})^{*} \rho = dt \wedge \rho_{0}(t) + \rho_{1}(t)
\label{chi-ro}
\ee
where $\rho_{0}(t)$ (resp. $\rho_{1}(t)$) are $q-1$ (resp. $q$) forms which do
not contain the differential $dt$. Then the {\it contact homotopy operator} is
by definition the map 
$A: \Omega^{r}_{q}(V) \rightarrow \Omega^{r}_{q-1}(V)$
given by:
\be
A \rho \equiv \int_{0}^{1} \rho_{0}(t).
\label{A-rho}
\ee

Moreover one has:
\be
\rho_{1}(1) = \rho, \quad \rho_{1}(0) = (\tau_{r})^{*} (\zeta_{0})^{*} \rho
\label{ro(1)}
\ee
where
$\tau_{r}: V^{r} \rightarrow V$ is the canonical projection on the first 
component:
$$
\tau_{r}(x^{i},y^{\sigma},y^{\sigma}_{j},...,y^{\sigma}_{j_{1},...,j_{s}}) 
\equiv (x^{i})
$$
and
$\zeta_{0}: U \rightarrow J^{y}W$
is the zero section given by:
$$
\zeta_{0}(x^{i}) \equiv (x^{i},0,...,0).
$$

Then we have (see \cite{Kr4}):
\begin{lemma}
(i) Let $\rho \in \Omega^{r}W$ be arbitrary. The following formula is true:
\be
\rho = Ad\rho + dA\rho + (\tau_{r})^{*} (\zeta_{0})^{*} \rho
\label{a-ro}
\ee

(ii)
If $\rho$ is contact, then:
\be
\rho = Ad\rho + dA\rho
\label{poincare}
\ee
and
\be
p_{k-1}A\rho = Ap_{k}\rho \quad (k = 1,...,q).
\label{p-A-p}
\ee
\label{ro-A-ro}
\end{lemma}

{\bf Proof:}
(i) The proof of the first formula is standard. First one finds out from 
(\ref{chi-ro}) that:
$$
(\chi_{r})^{*} d\rho = d_{0}\rho_{1} + dt \wedge 
\left( {\partial \rho_{1} \over \partial t} - d_{0}\rho_{0}\right)
$$
where $d_{0}$ is the exterior differentiation with respect to all variables
except $t$. This produces by integration:
$$
A d\rho = \rho_{1}(1) - \rho_{1}(0) - dA\rho
$$
and using (\ref{ro(1)}) we find the formula.

(ii) For the first relation one uses the structure theorem \ref{contact} for 
the contact forms  and sees that last term in (\ref{a-ro}) is zero. For the 
second relation one applies 
$(id \times \pi^{r+1,r})^{*}$
to the relation (\ref{chi-ro}) and after some simple manipulations one gets:
$$
(\pi^{r+1,r})^{*}\rho_{0} = \sum_{k=0}^{q} (p_{k}\rho)_{0} \Rightarrow
(\pi^{r+1,r})^{*} A\rho = \sum_{k=0}^{q} p_{k}A\rho =
\sum_{k=1}^{q} Ap_{k}\rho.
$$

Now use is made of lemma \ref{p-k-ro} and the second formula follows.
$\qed$

The central result of \cite{Kr4} follows. We will insist only on the part of
the proof which can be simplified with the Fock space tricks.

\begin{thm}
We consider the maps 
$d: \theta^{r}_{q} \longrightarrow \theta^{r}_{q+1}$; 
then the long sequence of sheaves
\be
0 \longrightarrow \theta^{r}_{1} \harr{d}{} \theta^{r}_{2} 
\harr{d}{} \cdots \harr{d}{} \theta^{r}_{P} \longrightarrow 0
\ee

is exact.
\end{thm}

{\bf Proof:}

(i) The exactness in $\theta^{r}_{1}$ is elementary \cite{Kr4}. If 
$\beta \in \theta^{r}_{1} = \Omega^{r}_{1(c)}$
then the structure theorem \ref{contact} can be used to write in $V^{r}$:
$$
\beta = \sum_{|J|\leq r-1} \Phi^{J}_{\sigma} \omega^{\sigma}_{J}  
$$
with $\Phi^{J}_{\sigma}$ some smooth functions on $V^{r}$. Then
$$
d\beta = - \sum_{|J| = r-1} \Phi^{J}_{\sigma} dy^{\sigma}_{Ji} \wedge dx^{i} +
\cdots 
$$
where by $\cdots$ we mean terms without the differentials
$dy^{\sigma}_{I} \quad (|I| = r)$;
so, the closedness of $\beta$ gives us 
$\Phi^{J}_{\sigma} = 0 \quad (|J| = r-1)$. This means that, in fact, in the
expression above of $\beta$ the sum finishes at $r - 2$. Continuing by
recurrence we arrive at $\beta = 0$.

(ii) We prove the exactness in $\theta^{r}_{q}, \quad (q = 2,...,q)$. If
$\beta \in \theta^{r}_{q} = \Omega^{r}_{q(c)}$ 
then again we can apply the structure theorem \ref{contact} to write (after 
simple rearrangements): 
$$
\beta = 
\sum_{|J|\leq r-1} \omega^{\sigma}_{J} \wedge \tilde{\Phi}^{J}_{\sigma} +
\sum_{|I|= r-1} d\omega^{\sigma}_{I} \wedge \Psi^{I}_{\sigma} =
\sum_{|J|\leq r-1} \omega^{\sigma}_{J} \wedge \Phi^{J}_{\sigma} +
\sum_{|I|= r-1} d(\omega^{\sigma}_{I} \wedge \Psi^{I}_{\sigma})
$$
so the closedness condition is
$$
\sum_{|J|\leq r-1} d\omega^{\sigma}_{J} \wedge \Phi^{J}_{\sigma} -
\sum_{|J| \leq r-1} \omega^{\sigma}_{I} \wedge d\Phi^{I}_{\sigma} = 0.
$$

Applying $p_{1}$ to this relation we find out with lemma \ref{p-k}:
\be
\sum_{|J|\leq r-1} d\omega^{\sigma}_{J} \wedge h\Phi^{J}_{\sigma} -
\sum_{|J| \leq r-1} \omega^{\sigma}_{I} \wedge hd\Phi^{I}_{\sigma} = 0.
\label{Phi-q}
\ee

One uses now for
$\Phi^{J}_{\sigma}$ 
the standard form
$$
\Phi^{J}_{\sigma} = \chi^{J}_{\sigma} + \zeta^{J}_{\sigma} 
$$
where
$\chi^{J}_{\sigma}$
is generated by 
$\omega^{\sigma}_{K},\quad (|K| \leq r - 1)$ 
and $\zeta^{J}_{\sigma}$
is a polynomial of degree $q-1$ in $dx^{i}$ and 
$dy^{\sigma}_{I}, \quad (|I| = r - 1)$:
$$
\zeta^{I_{1}}_{\nu_{1}} = \sum_{s=1}^{q} {1 \over s!(q-s)!} 
\sum_{|I_{2}|,...,|I_{s}|=r}
A^{I_{1},...,I_{s}}_{\nu_{1},...,\nu_{s},i_{s+1},...,i_{q}}
dy^{\nu_{2}}_{I_{2}} \wedge \cdots dy^{\nu_{s}}_{I_{s}} \wedge
dx^{i_{s+1}} \wedge \cdots \wedge dx^{i_{q}}
$$
where the symmetry properties of the type (\ref{symmetry}) leave out the
indices $I_{1}$ and $\nu_{1}$. Then (\ref{Phi-q}) becomes:
$$
\sum_{|J| = r-1} dy^{\sigma}_{Ji} \wedge dx^{i} \wedge h\zeta^{J}_{\sigma} +
\cdots = 0
$$
where by $\cdots$ we mean an expression which does not contain the
differentials $dy^{\sigma}_{I}, \quad (|I| = r)$.
This relations leads to
$$
{\cal S}^{+}_{Ii} dx^{i} \wedge h\zeta^{I}_{\sigma} = 0;
$$
like in the proof of theorem \ref{contact} this is equivalent to
$$
{\cal S}^{+}_{I_{1}i_{1}} \cdots {\cal S}^{+}_{I_{s}i_{s}} 
{\cal S}^{-}_{i_{1},...,i_{q}} 
\delta^{p_{1}}_{i_{1}} \cdots \delta^{p_{s}}_{i_{s}} 
A^{I_{1},...,I_{s}}_{\sigma_{1},...,\sigma_{s},i_{s+1},...,i_{q}} = 0 \quad
(s = 1,...,q),
$$
or, in tensor notations:
$$
B_{1} \cdots B_{s} A_{\nu_{1},...,\nu_{s}} = 0 \quad (s = 1,...,q).
$$

Using the usual argument, it follows that the forms 
$\zeta^{I}_{\sigma} \quad (|I| = r - 1)$
are generated by the differentials
$d\omega^{\nu}_{K}, \quad (|K| = r - 1)$. It follows that 
$\Phi^{J}_{\sigma}, \quad (|J| = r - 1)$ is are contact forms. By recurrence
it follows that $\Phi^{J}_{\sigma}, \quad (|J| \leq r - 1)$ is are contact
forms. This information must be inserted back in the initial expression of
the form $\beta$; one has that
$$
\beta = \beta_{0} + d\gamma
$$
where $\beta_{0}$ (resp. $\gamma$) are 2-contact (resp. contact) forms. The
closedness condition reduces now to $d\beta_{0} = 0$ and (\ref{poincare}) 
can be applied; one gets
$\beta_{0} = dA\beta_{0}$ so finally it follows that $\beta$ is given by
$\beta = d(A\beta_{0} + \gamma)$. But, using (\ref{p-A-p}) one gets that 
$h(A\beta_{0} + \gamma) = A p_{1} \beta_{0} = 0$ i.e. 
$A\beta_{0} + \gamma \in \theta^{r}_{q-1}$. In other words 
$\beta \in {\rm Im}(d)$.

(iii) The proof of the exactness in $\theta^{r}_{q}\quad (q > n)$ is also
standard \cite{Kr4}. Let
$\beta \in \theta^{r}_{q}$
i.e.
$\beta = \beta_{0} + d\gamma$
where
$\beta_{0} \in \Omega^{r}_{q(c)}$
and
$\gamma \in \Omega^{r}_{q-1(c)}$
such that
$d\beta = 0$. 
Then 
$d\beta_{0} = 0$
and we can apply (\ref{poincare}) to obtain
$\beta_{0} = dA\beta_{0}$.
As a consequence
$\beta = d(A\beta_{0} + \gamma)$.
But formula (\ref{p-A-p}) also implies:
$p_{q-1-n}(A\beta_{0} + \gamma) = Ap_{q-n}\beta_{0} = 0$
so in fact
$A\beta_{0} + \gamma \in \Omega^{r}_{q-1(c)} \subset \theta^{r}_{q-1}$;
this gives 
$\beta \in d\theta^{r}_{q-1} \subset {\rm Im}(d)$.
$\qed$

This theorem has the following consequence. 

\begin{thm}
Let 
$
E_{q}: \Omega^{r}_{q}/\theta^{r}_{q} \rightarrow 
\Omega^{r}_{q+1}/\theta^{r}_{q+1}
$
be given by
\be
E_{q}([\rho]) \equiv [d\rho]
\label{Eq}
\ee
where $[\rho]$ is the class of $\rho$ modulo $\theta^{r}_{q}$.
Then the quotient sequence
\be
0 \longrightarrow \R \longrightarrow \Omega^{r}_{0} \harr{E_{0}}{}
\Omega^{r}_{1}\left/\theta^{r}_{1} \right. \harr{E_{1}}{}\cdots 
\harr{E_{P-1}}{} \Omega^{r}_{P}\left/\theta^{r}_{P}\right. \harr{E_{P}}{}
\Omega^{r}_{P+1} \harr{d}{} \cdots \harr{d}{}
\Omega^{r}_{N} \longrightarrow 0
\ee
is a acyclic resolution of the constant sheaf $\R$. In particular it is exact.
\label{sequence}
\end{thm}

One calls this the {\it variational sequence of order $r$ over $Y$} and 
denotes for simplicity:
${\cal V}^{r}_{q} \equiv \Omega^{r}_{q}/\theta^{r}_{q}$.

Some special classes have distinct names. So, if 
$\lambda \in \Omega^{r}_{n}$ 
then the class
$[\lambda] \in \Omega^{r}_{n+1}/\theta^{r}_{n+1}$
is called the {\it Euler-Lagrange class} of $\lambda$. If
$T \in \Omega^{r}_{n+1}$
then
$[T] \in \Omega^{r}_{n+2}/\theta^{r}_{n+2}$
is called the {\it Helmholtz-Sonin class} of $T$.

Let us note in the end of this subsection that for any
$q= n + 1,...,P, \quad s > r$
there exists a canonical isomorphism
$$
i_{s,r}: \Omega^{r}_{q}\left/ \theta^{r}_{q}\right. \rightarrow 
{\rm Im}\left( \tau^{s}_{q} \circ (\pi^{s,r+1})^{*} \circ p_{q-n}\right), 
$$
where
$
\tau^{s}_{q}: \Omega^{s}_{q} \rightarrow 
\Omega^{s}_{q}\left/ \theta^{s}_{q}\right.
$
is the canonical projection. The explicit expression is
\be
i_{s,r} ([\rho]) = \tau^{s}_{q} \circ (\pi^{s,r+1})^{*} \circ p_{q-n}(\rho)
\ee
and this definition is consistent \cite{Kr4}.
\newpage

\subsection{Characterisation of the Variational Sequence by Forms}

In this subsection we give characterizations of 
${\cal V}^{r}_{q}, \quad q = n, n + 1$ 
using Euler-Lagrange and Helmholtz-Sonin forms. 
First we have 

\begin{thm}
If 
$\lambda \in \Omega^{r}_{n,X}$
is any Lagrange form, then we have for any $s > 2r$
\be
i_{s,r}(E_{n}([\lambda])) = [E(\lambda)]
\ee
where $E(\lambda)$ is the Euler-Lagrange form associated to $\lambda$ (see
(\ref{EL1}) and (\ref{EL2})).
\label{En}
\end{thm}

{\bf Proof:}
We take the proof from \cite{Kr4}. If in the chart
$(V^{r},\psi^{r})$
the expression of $\lambda$ is given by (\ref{L}), then we have on
$(V^{s},\psi^{s})$:
$$
(\pi^{s,r})^{*} d\lambda = \sum_{k=0}^{r} \left(
\partial^{j_{1},...,j_{k}}_{\sigma} L\right) 
\omega_{j_{1},...,j_{k}}^{\sigma} \wedge \theta_{0}.
$$

Now we note that using (\ref{theta-i}) we have
$$
d(\omega_{j_{1},...,j_{k-1}}^{\sigma} \wedge \theta_{j_{k}}) =
- \omega_{j_{1},...,j_{k}}^{\sigma} \wedge \theta_{0} 
$$
for $k \geq 1$
so the preceding relation can be rearranged as follows:
$$
(\pi^{s,r})^{*} d\lambda = 
(\partial_{\sigma} L) \omega^{\sigma} \wedge \theta_{0} +
\sum_{k=1}^{r} \left( d_{j_{k}} \partial^{j_{1},...,j_{k}}_{\sigma} L\right) 
\omega_{j_{1},...,j_{k-1}}^{\sigma} \wedge \theta_{0} + dF_{1} + G_{1}
$$
where $F_{1}$ (resp. $G_{1}$) is a contact (resp. 2-contact form).

One can iterate the procedure and proves by induction on $l$ that:
\begin{eqnarray}
(\pi^{s,r})^{*} d\lambda = 
\left[ \sum_{k=0}^{l} (-1)^{k} \left( d_{j_{1}}\cdots d_{j_{k}} 
\partial^{j_{1},...,j_{k}}_{\sigma} L\right) \right] \omega^{\sigma} \wedge
\theta_{0} + \nonumber \\ + 
(-1)^{l} \sum_{k=l+1}^{r} \left( d_{j_{k-l+1}}\cdots d_{j_{k}} 
\partial^{j_{1},...,j_{k}}_{\sigma} L\right)
\omega^{\sigma}_{j_{1},...,j_{k-l}} \wedge \theta_{0} + dF_{l} + G_{l}
\quad (l = 1,2,...,r) \nonumber
\end{eqnarray}
where $F_{l}$ (resp. $G_{l}$) is a contact (resp. 2-contact form).

In particular if we take $l = r$ we get
$$
(\pi^{s,r})^{*} d\lambda = 
\left[ \sum_{k=0}^{r} (-1)^{k} \left( d_{j_{1}}\cdots d_{j_{k}} 
\partial^{j_{1},...,j_{k}}_{\sigma} L\right) \right] \omega^{\sigma} \wedge
\theta_{0} + dF_{r} + G_{r} =
E_{\sigma}(L) \omega^{\sigma} \wedge \theta_{0} + dF_{r} + G_{r}
$$
i.e.
$$
(\pi^{s,r})^{*} d\lambda =
E_{\sigma}(L) \omega^{\sigma} \wedge \theta_{0} \quad ({\rm mod}\quad
\theta^{r}_{n+1}). 
$$
This is exactly the formula from the statement.
$\qed$

Similarly we have:

\begin{thm}
Let 
$T \in \Omega^{s}_{n+1,X}$ 
be a differential equation. The for any $t >2s$ the following formula is
true: 
\be
i_{t,s} (E_{n+1}([T])) = \left[ {1\over 2} H(T)\right]
\ee
where $H(T)$ is the Helmholtz-Sonin form associated to $T$.
\label{En+1}
\end{thm}

{\bf Proof:}

Suppose that in the chart 
$(V^{s},\psi^{s})$
we have the local expression (\ref{T}) for $T$; then we have on
$(V^{t},\psi^{t})$:
$$
(\pi^{t,s})^{*} dT = \sum_{k=0}^{s} \left( 
\partial^{j_{1},...,j_{k}}_{\nu} T_{\sigma}\right) 
\omega_{j_{1},...,j_{k}}^{\nu} \wedge \omega^{\sigma} \wedge \theta_{0}.
$$

Using the same trick as in the preceding theorem we can rewrite this as
follows: 
\begin{eqnarray}
(\pi^{t,s})^{*} dT =  
\partial_{\nu} T_{\sigma} \omega^{\sigma} \wedge \omega^{\nu} \wedge \theta_{0}
 - \sum_{k=1}^{s}\left( d_{j_{k}} \partial^{j_{1},...,j_{k}}_{\nu} 
T_{\sigma}\right) \omega^{\sigma} \wedge \omega_{j_{1},....,j_{k-1}}^{\nu} 
\wedge \theta_{0} \nonumber \\
- \sum_{k=1}^{s}\left( \partial^{j_{1},...,j_{k}}_{\nu} 
T_{\sigma}\right) \omega^{\sigma}_{j_{k}} \wedge 
\omega_{j_{1},...,j_{k-1}}^{\nu} \wedge \theta_{0} + dF_{1} + G_{1} \nonumber
\end{eqnarray}
where $F_{1}$ (resp. $G_{1}$) is 2-contact (resp. 3-contact).

We want to use the same idea as in the preceding theorem. We have to find the 
proper induction hypothesis. After some thought this proves to be:
\begin{eqnarray}
(\pi^{t,s})^{*} dT =  
\sum_{k=0}^{l} (-1)^{p} \sum_{p=k}^{l} {p \choose k}
\left( d_{j_{k+1}} \cdots d_{j_{p}} \partial^{j_{1},...,j_{p}}_{\nu} 
T_{\sigma}\right) \omega^{\sigma}_{j_{1},...,j_{k}} \wedge \omega^{\nu} 
\wedge \theta_{0}  + \nonumber \\
+ (-1)^{l} \sum_{k=l+1}^{s} \sum_{p=k-l}^{k} {l \choose k-p} 
\left( d_{j_{p+1}} \cdots d_{j_{k}} \partial^{j_{1},...,j_{k}}_{\nu} 
T_{\sigma}\right) \omega^{\sigma}_{j_{k-l+1},...,j_{p}} \wedge 
\omega_{j_{1},...,j_{k-l}}^{\nu} \wedge \theta_{0} + dF_{l} + G_{l} \nonumber
\end{eqnarray}
where $F_{l}$ (resp. $G_{l}$) is 2-contact (resp. 3-contact) and $l = 1,...,s$.

The induction from $l$ to $l + 1$ is accomplished with the same trick only the
computations become much more involved. For $l = s$ the preceding formula 
gives:
$$
(\pi^{t,s})^{*} dT =  
\sum_{k=0}^{s} (-1)^{k} \sum_{p=k}^{l} (-1)^{p-k} {p \choose k}
\left( d_{j_{k+1}} \cdots d_{j_{p}} \partial^{j_{1},...,j_{p}}_{\nu} 
T_{\sigma}\right) \omega^{\sigma}_{j_{1},...,j_{k}} \wedge \omega^{\nu} 
\wedge \theta_{0} + dF_{s} + G_{s}
$$
where $F_{s}$ (resp. $G_{s}$) is 2-contact (resp. 3-contact).
Using the definition of the Lie-Euler expressions (\ref{Lie-Euler}) we can 
write this formula in a more compact way:
$$
(\pi^{t,s})^{*} dT =  
\sum_{|J| \leq s} (-1)^{|J|} E^{J}_{\sigma}(T_{\sigma}) \omega^{\sigma}_{J} 
\wedge \omega^{\nu} \wedge \theta_{0} + dF_{s} + G_{s}.
$$

We add this result to the initial expression for 
$(\pi^{t,s})^{*} dT$
and divide by 2. If we use the definition for the Helmholtz-Sonin forms
(\ref{HS-local}) we obtain:
$$
(\pi^{t,s})^{*} dT = {1\over 2} H(T) + dF_{s} + G_{s} = 
{1\over 2} H(T) \quad ({\rm mod} \quad \theta^{r}_{n+2}). \qquad
\qed
$$
\newpage 

\section{Variationally Trivial Lagrangians}

A {\it variationally trivial Lagrange form of order r} is any
$\lambda \in \Omega^{r}_{n,X}$
such that
$E(\lambda) = 0$.
In other words, the corresponding Euler-Lagrange expressions are identically
zero, or the Euler-Lagrange equations are identities for any section
$\gamma$.
Here we give the general form of such a Lagrange form. We follow,
essentially \cite{Kr5} and \cite{AD} but as before we simplify considerably the
proofs using the techniques developed in Section 3.2.

First, we remind the reader that 
$\rho \in \Omega^{r+1}_{q}$
is called $\pi^{r+1,r}$-{\it projectable iff} there exists 
$\rho' \in \Omega^{r}_{q}$
such that
\be
\rho = (\pi^{r+1,r})^{*} \rho'.
\ee

One can easily see that locally this amounts to the condition that if $\rho$
is written as a polynomial in 
$dx^{i}$ and $dy^{\sigma}_{J} \quad (|J| \leq r+1)$
then, the differentials
$dy^{\sigma}_{J} \quad (|J| = r+1)$
must be absent and moreover, the coefficient functions must not depend on
$y^{\sigma}_{J} \quad (|J| = r + 1)$.

We start with the following result \cite{Kr4}:

\begin{prop}
Let 
$\eta \in \Omega^{r}_{q}\quad (q = 1,...,n-1)$
such that
$hd\eta$
is a $\pi^{r+1,r}$-projectable $(q+1)$-form. Then one can write $\eta$ as 
follows:
\be
\eta = \nu + d\phi + \psi
\ee
where
$\nu \in \Omega^{r}_{q,X}$ 
is a basic $q$-form and 
$\psi \in \Omega^{r}_{q(c)}$
is a contact $q$-form.
\label{eta-nu}
\end{prop}

{\bf Proof:}

Let
$(V,\psi)$ be a chart on $Y$. Then in the associated chart
$(V^{r},\psi^{r})$
one can write $\eta$
in the standard form
\be
\eta = \eta_{0} + \eta_{1}
\label{eta}
\ee
where in $\eta_{0}$ we collect all terms containing at least one factor
$\omega^{\sigma}_{J} \quad (|J| \leq r-1)$
and $\eta_{1}$ is a polynomial only in 
$dx^{i}$ and $dy^{\sigma}_{J} \quad (|J| = r+1)$:
\be
\eta_{1} = \sum_{s=0}^{q}  \sum_{|I_{1}|,...,|I_{s}|=r}
A^{I_{1},...,I_{s}}_{\sigma_{1},...,\sigma_{s},i_{s+1},...,i_{q}}
dy^{\sigma_{1}}_{I_{1}} \wedge \cdots \wedge dy^{\sigma_{s}}_{I_{s}} \wedge
dx^{i_{s+1}} \wedge \cdots \wedge dx^{i_{q}},
\label{eta1}
\ee
where the coefficients
$A^{I_{1},...,I_{s}}_{\sigma_{1},...,\sigma_{s},i_{s+1},...,i_{q}}$
have antisymmetry properties of type (\ref{symmetry}). 

In particular $\eta_{0}$ is a contact form so we have
$hd\eta = hd\eta_{1}$ 
i.e. the form $hd\eta_{1}$ is, by hypothesis, $\pi^{r+1,r}$-projectable.
One first computes in general:
\begin{eqnarray}
d\eta_{1} = \sum_{s=0}^{q} \sum_{|I_{1}|,...,|I_{s}|=r}
\sum_{|J| \leq r-1} \left( \partial^{J}_{\nu} 
A^{I_{1},...,I_{s}}_{\sigma_{1},...,\sigma_{s},i_{s+1},...,i_{q}} \right)
\omega^{\nu}_{J} \wedge dy^{\sigma_{1}}_{I_{1}} \wedge \cdots 
\wedge dy^{\sigma_{s}}_{I_{s}} \wedge dx^{i_{s+1}} \wedge \cdots \wedge 
dx^{i_{q}} \nonumber \\
+ \sum_{s=0}^{q+1} \sum_{|I_{1}|,...,|I_{s}|=r}
\tilde{A}^{I_{1},...,I_{s}}_{\sigma_{1},...,\sigma_{s},i_{s+1},...,i_{q+1}}
dy^{\sigma_{1}}_{I_{1}} \wedge \cdots \wedge dy^{\sigma_{s}}_{I_{s}} \wedge
dx^{i_{s+1}} \wedge \cdots \wedge dx^{i_{q+1}}, \nonumber
\label{d-eta1-q}
\end{eqnarray}
where
\begin{eqnarray}
\tilde{A}^{I_{1},...,I_{s}}_{\sigma_{1},...,\sigma_{s},i_{s+1},...,i_{q+1}} 
\equiv {1\over s} \sum_{k=1}^{s} (-1)^{k-1} \partial^{I_{k}}_{\sigma_{k}} 
A^{I_{1},...,\hat{I_{k}},...,I_{s}}_{\sigma_{1},...,\hat{\sigma_{k}},...,
\sigma_{s},i_{s+1},...,i_{q+1}} +  \nonumber \\ 
{1\over q+1-s} \sum_{k=s+1}^{q+1} (-1)^{k-1} d_{i_{k}} 
A^{I_{1},...,I_{s}}_{\sigma_{1},...,\sigma_{s},
i_{s+1},...,\hat{i_{k}},...,i_{q+1}}. \nonumber
\end{eqnarray}
or, equivalently
\begin{eqnarray}
\tilde{A}^{I_{1},...,I_{s}}_{\sigma_{1},...,\sigma_{s},i_{s+1},...,i_{q+1}} 
\equiv {\cal S}^{-}_{(I_{1},\sigma_{1}),...,(I_{s}\sigma_{s})}
\partial^{I_{1}}_{\sigma_{1}} A^{I_{2},...,I_{s}}_{\sigma_{2},...,\sigma_{s},
i_{s+1},...,i_{q+1}} +  \nonumber \\
(-1)^{s} {\cal S}^{-}_{i_{s+1},...,i_{q+1}} d_{i_{s+1}} 
A^{I_{1},...,I_{s}}_{\sigma_{1},...,\sigma_{s},i_{s+2},...,i_{q+1}}
\label{tilde-A}
\end{eqnarray}
where
${\cal S}^{-}_{(I_{1},\sigma_{1}),...,(I_{s}\sigma_{s})}$
is the antisymmetrization projector in the corresponding couples of indices.

Then one gets
$$
hd\eta = hd\eta_{1} = \sum_{s=0}^{q+1} \sum_{|I_{1}|,...,|I_{s}|=r}
\tilde{A}^{I_{1},...,I_{s}}_{\sigma_{1},...,\sigma_{s},i_{s+1},...,i_{q+1}}
y^{\sigma_{1}}_{I_{1}i_{1}} \cdots y^{\sigma_{s}}_{I_{s}i_{s}} 
dx^{i_{1}} \wedge \cdots \wedge dx^{i_{q+1}}
$$
which is $\pi^{r+1,r}$-projectable {\it iff} the coefficient functions do not
depend on 
$y^{\sigma}_{I}\quad (|I| = r+1)$.

(ii) By repeatedly applying derivative operators, one obtains, as in 
\cite{Gr4}, that the condition of $\pi^{r+1,r}$-projectability is equivalent to
$$
{\cal S}^{-}_{i_{1},...,i_{q+1}} {\cal S}^{+}_{I_{1}p_{1}} \cdots 
{\cal S}^{+}_{I_{s}p_{s}} 
\tilde{A}^{I_{1},...,I_{s}}_{\sigma_{1},...,\sigma_{s},i_{s+1},...,i_{q+1}}
\delta^{p_{1}}_{i_{1}} \cdots \delta^{p_{s}}_{i_{s}} = 0, \quad
(s = 1,...,q+1)
$$
or, using a familiar argument from \cite{Gr4}:
\be
B_{1} \cdots B_{s} \tilde{A}_{\sigma_{1},...,\sigma_{s}} = 0, \quad
(s = 1,...,q+1).
\label{var}
\ee

It is possible to apply lemma \ref{BXB} and one obtains the following generic
expression: 
$$
\tilde{A}_{\sigma_{1},...,\sigma_{s}} = \sum_{\alpha =1}^{s}
B_{\alpha} \tilde{A}_{\sigma_{1},...,\sigma_{s}}^{\alpha}
$$
or, in index notations, this means the function
$\tilde{A}^{I_{1},...,I_{s}}_{\sigma_{1},...,\sigma_{s},i_{s+1},...,i_{q+1}}$
is a sum of terms containing at least a factor of the type
$\delta^{j}_{i}$
where
$j$
belongs to one of the multi-indices
$I_{p}, \quad p=1, \dots, s$
and
$i \in \{i_{s+1},...,i_{q+1}\}.$

If we insert back into the expression of 
$d\eta_{1}$
we get after minor prelucrations that 
\be
d\eta_{1} = \nu_{1} + d\phi + \sum_{|J| \leq r-1} \omega^{\sigma}_{J} \wedge
\Phi^{J}_{\sigma}
\label{deta1}
\ee
with
$\nu_{1} \in \Omega^{r}_{q,X}$
a basic $q$-form,
$\phi \in \Omega^{r}_{q(c)}$
and
$\Phi^{J}_{\sigma}$
a polynomial of degree $q$ in 
$dx^{i}$
and
$dy^{\sigma}_{I} \quad |I| = r.$

This relation implies that
\be
d\nu_{1} + \sum_{|J| \leq r-1} \left( d\omega^{\sigma}_{J} \wedge
\Phi^{J}_{\sigma} - \omega^{\sigma}_{J} \wedge d\Phi^{J}_{\sigma}\right) = 0.
\ee

One applies the operator $p_{1}$ to this equality and uses (\ref{p-wedge}); a 
new relation is obtained, namely:
\be
p_{1} d\nu_{1} + \sum_{|J| \leq r-1} \left( d\omega^{\sigma}_{J} \wedge
h\Phi^{J}_{\sigma} - \omega^{\sigma}_{J} \wedge hd\Phi^{J}_{\sigma}\right) = 0.
\label{dnu1}
\ee

One must consider now the generic forms for
$\nu_{1}$
and
$\Phi^{J}_{\sigma}$
namely: 
\be
\nu_{1} = A_{i_{1},...,i_{q+1}} dx^{i_{1}} \wedge \cdots \wedge dx^{i_{q+1}}
\label{nu1}
\ee
and
\be
\Phi^{J}_{\sigma} = \sum_{s=0}^{q}  \sum_{|I_{1}|,...,|I_{s}|=r}
\Phi^{J,I_{1},...,I_{s}}_{\sigma,\nu_{1},...,\nu_{s},i_{s+1},...,i_{q}}
dy^{\nu_{1}}_{I_{1}} \wedge \cdots \wedge dy^{\nu_{s}}_{I_{s}} \wedge
dx^{i_{s+1}} \wedge \cdots \wedge dx^{i_{q}};
\label{Phi-J-sigma}
\ee
here we can assume a (partial) symmetry property of the type (\ref{symmetry}):
\be
\Phi^{J,I_{P(1)},...,I_{P(s)}}_{\sigma,\nu_{P(1)},...,\nu_{P(s)},
i_{Q(s+1)},...,i_{Q(q)}} = (-1)^{|P|+|Q|}
\Phi^{J,I_{1},...,I_{s}}_{\sigma,\nu_{1},...,\nu_{s},i_{s+1},...,i_{q}}
\label{sym-partial}
\ee
and we make the convention that
\be
\Phi^{J,I_{1},...,I_{s}}_{\sigma,\nu_{1},...,\nu_{s},i_{s+1},...,i_{q}} = 0,
\quad \forall J \quad s.t. \quad |J| \geq r.
\ee

These expression must be plugged into the equation (\ref{dnu1}).
One finds out from this equation the following consequence:
\begin{eqnarray}
{\cal S}^{-}_{i_{1},...,i_{q+1}} {\cal S}^{+}_{I_{1}p_{1}} \cdots 
{\cal S}^{+}_{I_{s}p_{s}} 
\left[ \tilde{\Phi'}^{\{j_{1},...,j_{k}\},I_{1},...,I_{s}}_{\sigma,\nu_{1},...,
\nu_{s},i_{s+1},...,i_{q+1}} + (-1)^{s} {\cal S}_{j_{1},....,j_{k}}
\delta^{j_{1}}_{i_{s+1}} \Phi^{\{j_{2},...,j_{k}\},I_{1},...,I_{s}}_{\sigma,
\nu_{1},...,\nu_{s},i_{s+2},...,i_{q+1}} \right] \times \nonumber \\
\delta^{p_{1}}_{i_{1}} \cdots \delta^{p_{s}}_{i_{s}} = 0, 
\quad (s = 1,...,q+1, k = 0,...,r).
\label{dnu1a}
\end{eqnarray}

Here we have defined, in analogy with (\ref{tilde-A}):
\begin{eqnarray}
\tilde{A'}^{I_{0},...,I_{s}}_{\sigma_{0},...,\sigma_{s},i_{s+1},...,i_{q+1}} 
\equiv {\cal S}^{-}_{(I_{1},\sigma_{1}),...,(I_{s}\sigma_{s})}
\partial^{I_{1}}_{\sigma_{1}} A^{I_{0},I_{2},...,I_{s}}_{\sigma_{0},\sigma_{2},
...,\sigma_{s},i_{s+1},...,i_{q+1}} +  \nonumber \\
(-1)^{s} {\cal S}^{-}_{i_{s+1},...,i_{q+1}} d_{i_{s+1}} 
A^{I_{0},I_{1},...,I_{s}}_{\sigma_{0},...,\sigma_{s},i_{s+2},...,i_{q+1}}. 
\label{tilde-A'}
\end{eqnarray}

To simplify the analysis of the relation (\ref{dnu1a}), one observes that it is 
possible to define a map
$$
\Delta': \oplus_{s=0}^{q} {\cal H}_{q-s,k,\underbrace{r,...,r}_{s-times}}
\rightarrow
\oplus_{s=0}^{q+1} {\cal H}_{q+1-s,k,\underbrace{r,...,r}_{s-times}}
$$
according to
\be
(\Delta' \Phi)^{J,I_{1},...,I_{s}}_{\sigma,\nu_{1},...,\nu_{s},
i_{s+1},...,i_{q+1}} 
= \tilde{\Phi'}^{J,I_{1},...,I_{s}}_{\sigma,\nu_{1},...,\nu_{s},
i_{s+1},...,i_{q+1}}.
\label{delta'}
\ee

Then the equation above takes the form:
\be
B_{1} \cdots B_{s} \left((\Delta \Phi)^{k}_{s} + B_{0} \Phi^{k-1}_{s}\right) 
= 0, \quad (s = 1,...,q+1, k = 0,...,r)
\label{descent}
\ee
where
$\Phi^{k}_{s} \in {\cal H}_{q-s,k,\underbrace{r,...,r}_{s-times}}$
is the tensor of components
$\Phi^{J,I_{1},...,I_{s}}_{\sigma,\nu_{1},...,\nu_{s},i_{s+1},...,i_{q}},
\quad |J| = k.$

We remind the reader that we have made the convention
\be
\Phi^{-1}_{s} = 0, \quad \Phi^{r}_{s} = 0.
\label{conv-psi}
\ee

To solve the preceding equation one first establishes by direct computation 
that
\be
\{\Delta', B_{0}\} = 0.
\label{delta-B}
\ee
and
\be
(\Delta')^{2} = 0.
\label{delta2}
\ee

Now we can solve the system (\ref{descent}) by a procedure similar to the
descent procedure applied in the BRST quantization scheme.

We take in (\ref{descent}) $k = r$ and, taking into account the convention
(\ref{conv-psi}), we obtain:
\be
B_{0} \cdots B_{s} \Phi^{r-1}_{s} = 0, \quad (s = 1,...,q).
\label{r-1}
\ee

Lemma \ref{BXB} can be applied and it follows that we have a generic expression
of the following form:
$$
\Phi^{r-1}_{s} = \sum_{\alpha =0}^{s} B_{\alpha} \Phi^{r-1,\alpha}_{s}, 
\quad (s = 1,...,q).
$$

If we substitute this expression into the initial expression
(\ref{Phi-J-sigma}) of 
$\Phi^{J}_{\sigma}, \quad |J| = r-1$
we find out that the contributions corresponding to
$\alpha = 1,...,s$
are producing contact terms. So, it appears that we can redefine the
expressions 
$\Phi^{J}_{\sigma}$
such that instead of (\ref{deta1}) we have 
\be
d\eta_{1} = \nu_{1} + d\phi + \sum_{|J| \leq r-1} \omega^{\sigma}_{J} \wedge
\Phi^{J}_{\sigma} + \psi
\label{deta1a}
\ee
with $\psi$ some 2-contact form and moreover
\be
\Phi^{r-1}_{s} = B_{0} C^{r-1}_{s}, \quad (s = 1,...,q)
\ee
for some tensors
$C^{r-1}_{s} \in {\cal H}_{q-s-1,r-2,\underbrace{r,...,r}_{s-times}}.$

The procedure can now be iterated taking into (\ref{descent}) $k = r-1$, etc.
The result of this descent procedure is the following: one can redefine the
tensors 
$\Phi^{J}_{\sigma}$
such that we have (\ref{deta1a}) and
\be
\Phi^{k}_{s} = B_{0} C^{k}_{s} + (\Delta C^{k+1})_{s}, \quad (s = 1,...,q; 
\quad k = 0,...,r-1)
\label{psi-C}
\ee
for some tensors
$C^{k}_{s} \in {\cal H}_{q-s-1,k-1,\underbrace{r,...,r}_{s-times}};$
by convention
$$
C^{0}_{s} = 0, \quad C^{r}_{s} = 0, \quad (s = 1,...,q).
$$

In full index notations this means that we have
\begin{eqnarray}
\Phi^{\{j_{1},...,j_{k}\},I_{1},...,I_{s}}_{\sigma,\nu_{1},...,\nu_{s},
i_{s+1},...,i_{q}} = 
(-1)^{s} {\cal S}^{+}_{j_{1},...,j_{k}} {\cal S}^{-}_{i_{s+1},...,i_{q}}
\delta^{j_{1}}_{i_{s+1}} 
C^{\{j_{2},...,j_{k}\},I_{1},...,I_{s}}_{\sigma,\nu_{1},...,\nu_{s},
i_{s+2},...,i_{q}} + \nonumber \\
\tilde{C}'^{\{j_{1},...,j_{k}\},I_{1},...,I_{s}}_{\sigma,\nu_{1},...,\nu_{s},
i_{s+1},...,i_{q}} \quad (s = 1,...,q; k = 0,...,r).
\label{phi-C-full}
\end{eqnarray}

Now one can substitute this expression for the tensors
$\Phi^{...}_{...}$
into the original expression for the forms (\ref{Phi-J-sigma}). After some
algebra one obtains that
\be
\sum_{|J| \leq r-1} \omega^{\sigma}_{J} \wedge \Phi^{J}_{\sigma} =
\sum_{|J| \leq r-1} \Phi^{J}_{\sigma,i_{1},...,i_{q}} \omega^{\sigma}_{J} 
\wedge dx^{i_{1}} \wedge \cdots \wedge dx^{i_{q}} + d\phi_{1}
\ee
where
$\phi_{1} \in \Omega^{r}_{q(c)}$
is a contact form.

It emerges that the relation (\ref{deta1a}) becomes:
\be
d\eta_{1} = \nu_{1} + d\phi_{1} + 
\sum_{|J| \leq r-1} \Phi^{J}_{\sigma,i_{1},...,i_{q}} \omega^{\sigma}_{J} 
\wedge dx^{i_{1}} \wedge \cdots \wedge dx^{i_{q}} + \psi
\label{deta1b}
\ee
where $\phi_{1}$ (resp. $\psi$) is some contact (resp. 2-contact) form
of degree $q$ (resp. $q+1$.)

(iii) The last step of our proof consists in using the contact homotopy
operator $A$ (see the definition contained in the relations (\ref{chi}) -
(\ref{A-rho})). We apply the relation (\ref{a-ro}) to the form
$\rho = \eta_{1} - \phi_{1}$:
\be
\eta_{1} - \phi_{1} = A \left( \nu_{1} + \sum_{|J| \leq r-1} 
\Phi^{J}_{\sigma,i_{1},...,i_{q}} \omega^{\sigma}_{J} 
\wedge dx^{i_{1}} \wedge \cdots \wedge dx^{i_{q}} + \psi \right) +
dA(\eta_{1} - \phi_{1}) + \eta_{1}(0) - \phi_{1}(0).
\label{eta-A}
\ee

One proves by direct computation the following facts:

- $A \nu_{1} = 0$;

- if $\psi$ is 2-contact form, then $A\psi$ is a contact form;

- $A(\sum \cdots )$ and $\eta_{1}(0)$ are basic forms;

- if $\phi_{1}$ is a contact form, then $\phi_{1}(0)$ is also a contact form.

Inserting this information in the preceding relation one obtains that the 
result from the statement of the proposition is true for the form
$\eta_{1}$. Taking into account (\ref{eta}) we obtain the same result for the 
form $\eta$.
$\qed$

Let us note that the converse of this statement is not true. In fact the 
condition of $\pi^{r+1,r}$-projectability imposes additional constraints on the
basic form $\nu$ which will be analysed in the next lemma. A complete proof of
the following result appears in \cite{AD}; here we offer an alternative proof
which seems to be much more simpler.

\begin{prop}
Let 
$\nu \in \Omega^{r}_{q,X} \quad (q = 1,...,n-1)$
Then the basic form $hd\nu$ is $\pi^{r+1,r}$-projectable {\it iff} there exists
for any chart 
$(V,\psi)$
a form 
$\tilde{\nu}_{V} \in \Omega^{r-1}_{q}$
such that we have in the associated chart
$(V^{r},\psi^{r})$ the equality 
$\nu = h\tilde{\nu}_{V}$.
\label{nu-nu-tilde}
\end{prop}

{\bf Proof:}
According to (\ref{hor1}) we have in the associated chart 
$(V^{r},\psi^{r})$:
$$
\nu = A_{i_{1},...,i_{q}} dx^{i_{1}} \wedge \cdots \wedge dx^{i_{q}}.
$$

Then one obtains by direct computation that:
$$
hd\nu = \left[ \left( d_{i_{0}} A_{i_{1},...,i_{q}}\right) + \sum_{|I|=r} 
\left( \partial^{I}_{\sigma} A_{i_{1},...,i_{q}}\right) y^{\sigma}_{Ii_{0}}
\right] dx^{i_{0}} \wedge \cdots \wedge dx^{i_{q}}.
$$
(here $d_{j}$ is the formal derivative on $V^{r}$).

This form is $\pi^{r+1,r}$-projectable {\it iff} the square bracket does not
depend on 
$y^{\sigma}_{I} \quad (|I| = r+1)$, i.e.
\be
{\cal S}^{-}_{i_{0},...,i_{q}} {\cal S}^{+}_{Ip} \delta^{p}_{i_{0}}
\left( \partial^{I}_{\sigma} A_{i_{1},...,i_{q}}\right) = 0.
\label{projectability}
\ee

It is clear that this relation is very similar to those already analysed. To be
able to apply the central result contained in lemma \ref{BXB} one must do a
little trick. From the previous relation we obtain by derivation:
\be
{\cal S}^{-}_{i_{0},...,i_{q}} {\cal S}^{+}_{I_{0}p_{0}} \delta^{p_{0}}_{i_{0}}
\left( \partial^{I_{0}}_{\sigma_{0}} \cdots \partial^{I_{q}}_{\sigma_{q}} 
A_{i_{1},...,i_{q}}\right) = 0.
\label{nu2}
\ee

If we define the tensor 
$A_{\sigma_{0},...,\sigma_{q}} \in 
{\cal H}_{q,\underbrace{r,...,r}_{(q+1)-times}}$
by
$$
A^{I_{0},...,I_{q}}_{\sigma_{0},...,\sigma_{q},i_{1},...,i_{q}} \equiv
\partial^{I_{0}}_{\sigma_{0}} \cdots \partial^{I_{q}}_{\sigma_{q}} 
A_{i_{1},...,i_{q}} 
$$
then the preceding equation (\ref{nu2}) can be compactly written as follows:
$$
B_{0} A_{\sigma_{0},...,\sigma_{q}} = 0.
$$

In fact, due to symmetry properties of the type (\ref{symmetry}), namely:
\be
A^{I_{P(1)},...,I_{P(s)}}_{\sigma_{P(1)},...,\sigma_{P(s)},
i_{Q(1)},...,i_{Q(q)}} = (-1)^{|Q|}
A^{I_{1},...,I_{s}}_{\sigma_{1},...,\sigma_{s},i_{1},...,i_{q}}
\label{new-sym}
\ee
we have
\be
B_{\alpha} A_{\sigma_{0},...,\sigma_{q}} = 0 \quad (\alpha = 0,...,q).
\label{nu3}
\ee

As a consequence we obtain:
$$
B_{0} \cdots B_{q} A_{\sigma_{0},...,\sigma_{q}} = 0
$$
and lemma \ref{BXB} can be applied; it follows that 
$A_{...}$
has the generic form
$$
A_{...} = \sum_{\alpha =0}^{q} B_{\alpha} A^{\alpha}_{...}
$$
for some tensors
$A^{\alpha}_{...} \in {\cal H}_{q,r,...,r,r-1,r,...,r}$
where the entry $r - 1$ is on the position $\alpha$.

One can plug this relation into the initial relation (\ref{nu3}) and a similar
relation is obtained for the tensors 
$A^{\alpha}_{...}$.
By recurrence, one gets:
$$
A_{...} = \sum_{\alpha_{0},...,\alpha_{s}=0}^{q} B_{\alpha_{0}} \cdots
B_{\alpha_{s}} A^{\alpha_{0},...,\alpha_{s}}_{...} \quad (s = 0,...,q)
$$
with
$A^{\alpha_{0},...,\alpha_{s}}_{...}$
some tensors in
${\cal H}_{q-s-1,r_{0},...,r_{q}}$.
In particular, for $s = q$ we obtain that in fact:
$A_{...} = 0$
i.e. 
\be
\partial^{I_{0}}_{\sigma_{0}} \cdots \partial^{I_{q}}_{\sigma_{q}} 
A_{i_{1},...,i_{q}} = 0.
\ee

In other words, the functions
$A_{i_{1},...,i_{q}} $
are polynomials in 
$y^{\sigma}_{I} \quad (|I| = r)$
of maximal degree $q$. So, the generic form of these functions is:
\be
A_{i_{1},...,i_{q}}  = \sum_{s=0}^{q} {1 \over s!(q-s)!} 
\sum_{|I_{1}|,...,|I_{s}|=r}
C^{I_{1},...,I_{s}}_{\sigma_{1},...,\sigma_{s},i_{1},...,i_{q}}
y^{\sigma_{1}}_{I_{1}} \cdots y^{\sigma_{s}}_{I_{s}}
\label{AA}
\ee
with
$C^{I_{1},...,I_{s}}_{\sigma_{1},...,\sigma_{s},i_{1},...,i_{q}}$
some smooth functions on $V^{r}$ having symmetry properties of the type
(\ref{new-sym}). Remark the fact that these functions do not depend on the
variables 
$y^{\sigma}_{I} \quad (|I| = r)$
and this justifies the fact that they live on the chart $V^{r-1}$.

Now we insert this generic expression into the projectability condition 
(\ref{projectability}) and we get, in the same way as before 
\be
{\cal S}^{-}_{i_{0},...,i_{q}} {\cal S}^{+}_{I_{1}p} \delta^{p}_{i_{0}}
C^{I_{1},...,I_{s}}_{\sigma_{1},...,\sigma_{s},i_{1},...,i_{q}} = 0 \quad
(s = 1,...,q-1)
\label{nu4}
\ee
or, in tensor notations:
$$
B_{1} C_{\sigma_{1},...,\sigma_{q}} = 0 \quad (s = 1,...,q-1).
$$

In fact, due to the symmetry properties one has from here:
\be
B_{\alpha} C_{\sigma_{1},...,\sigma_{q}} = 0 \quad (\alpha = 1,...,s;\quad 
s = 1,...,q-1). 
\ee

Using an argument familiar by now we get from here:
$$
C_{\sigma_{1},...,\sigma_{q}} = B_{1} \cdots B_{s} 
\tilde{C}_{\sigma_{1},...,\sigma_{q}}  \quad (s = 1,...,q-1).
$$

In full index notations this means that we have the following generic
expression:
\be
C^{J_{1}p_{1},...,J_{s}p_{s}}_{\sigma_{1},...,\sigma_{s},i_{1},...,i_{q}} =
{\cal S}^{+}_{J_{1}p_{1}} \cdots {\cal S}^{+}_{J_{s}p_{s}} 
{\cal S}^{-}_{i_{1},...,i_{q}} 
\delta^{p_{1}}_{i_{1}} \cdots \delta^{p_{s}}_{i_{s}} 
A^{J_{1},...,J_{s}}_{\sigma_{1},...,\sigma_{s},i_{1},...,i_{q}} 
\quad (s = 1,...,q)
\ee
where
$|J_{1}| = \cdots = |J_{s}| = r-1$
and
$A^{J_{1},...,J_{s}}_{\sigma_{1},...,\sigma_{s},i_{1},...,i_{q}}$
are some smooth function on $V^{r-1}$ completely antisymmetric in the indices
$i_{1},...,i_{q}.$

Inserting this expression in (\ref{AA}) one immediately obtains
\be
A_{i_{1},...,i_{q}}  = \sum_{s=0}^{q} {1 \over s!(q-s)!} 
\sum_{|J_{1}|,...,|J_{s}|=r-1} {\cal S}^{-}_{i_{1},...,i_{q}}
A^{J_{1},...,J_{s}}_{\sigma_{1},...,\sigma_{s},i_{1},...,i_{q}}
y^{\sigma_{1}}_{J_{1}i_{1}} \cdots y^{\sigma_{s}}_{J_{s}i_{s}}.
\label{AAA}
\ee

Let us remark that the expression
$
{\cal S}^{-}_{i_{1},...,i_{q}}
y^{\sigma_{1}}_{J_{1}i_{1}} \cdots y^{\sigma_{s}}_{J_{s}i_{s}}
$
is completely antisymmetric in the couples
$(I_{1},\sigma_{1}),...,(I_{s},\sigma_{s})$
so one can consider that the tensors
$A^{J_{1},...,J_{s}}_{\sigma_{1},...,\sigma_{s},i_{1},...,i_{q}}$
have the same property. In the end, it follows that they have the symmetry
property (\ref{symmetry}). 

The expression for the form $\nu$ becomes
\be
\nu = \sum_{s=0}^{q} {1 \over s!(q-s)!} 
\sum_{|J_{1}|,...,|J_{s}|=r-1} 
A^{J_{1},...,J_{s}}_{\sigma_{1},...,\sigma_{s},i_{s+1},...,i_{q}}
y^{\sigma_{1}}_{J_{1}i_{1}} \cdots y^{\sigma_{s}}_{J_{s}i_{s}}
dx^{i_{1}} \wedge \cdots \wedge dx^{i_{q}} = h\tilde{\nu}_{V}
\ee
where
\be
\tilde{\nu}_{V} \equiv \sum_{s=0}^{q} {1 \over s!(q-s)!} 
\sum_{|J_{1}|,...,|J_{s}|=r-1} 
A^{J_{1},...,J_{s}}_{\sigma_{1},...,\sigma_{s},i_{s+1},...,i_{q}}
dy^{\sigma_{1}}_{J_{1}}\wedge \cdots \wedge dy^{\sigma_{s}}_{J_{s}} \wedge
dx^{i_{s+1}} \wedge \cdots \wedge dx^{i_{q}}.
\label{nu-tilde}
\ee
This finishes the proof if we take into account that the relation 
(\ref{projectability}) is equivalent to the projectability condition.
$\qed$

Now we are ready to obtain the most general form of a variationally trivial
Lagrangian. First we note that we have

\begin{prop}
The Lagrange form
$\lambda \in \Omega^{r}_{n,X}$
is variationally trivial {\it iff} on has:
\be
E_{n}([\lambda]) = 0.
\ee
\end{prop}

{\bf Proof:}
It is an immediate consequence of the definition of a variationally trivial
Lagrange form and of theorem \ref{En}.
\label{e-n}
$\qed$

The central result now follows:
\begin{thm}
Let
$\lambda \in \Omega^{r}_{n,X}$
be a variationally trivial Lagrange form. Then for every point 
$j^{r}_{x}\gamma \in J^{r}Y$
there exists a neighbourhood $V$ of $\gamma(x) \in Y$
and a $n$-form $\rho_{V}$ defined in the chart
$(V^{r-1},\psi^{r-1})$
which is basic and we also have: (1) 
$\lambda = h\rho_{V}$ 
in the chart
$(V^{r},\psi^{r})$;
(2) $d\rho_{V} = 0$.
Conversely, if such a local form $\rho_{V}$ exists, then the form
$\lambda$ is variationally trivial.
\label{trivial}
\end{thm}

{\bf Proof:}
(i) According to the previous proposition and applying the exactness of the 
variational sequence (theorem \ref{sequence}), there exists 
$\eta \in \Omega^{r}_{n-1}$
such that
$$
[\lambda] = [d\eta],
$$
or, equivalently
\be
\lambda -d\eta \in \theta^{r}_{n} = \Omega^{r}_{n(c)}.
\label{lambda}
\ee

As a consequence we have 
$$
(\pi^{r+1,r})^{*} \lambda = h\lambda = h d\eta.
$$

In particular, it follows that the $n$-form $h d\eta$ must be
$\pi^{r+1,r}$-projectable. We can apply proposition \ref{eta-nu} and 
(\ref{lambda}) rewrites as follows:
\be
\lambda -d\nu \in \Omega^{r}_{n(c)}.
\ee
for some basic form $\nu$. From here we get 
\be
h\lambda = hd\nu
\label{lambda-nu}
\ee
or using (\ref{pi-ro-pi})
\be
(\pi^{r+1,r})^{*} \lambda = h d\nu.
\label{lambda-h-nu}
\ee

Let us note that that this relation is completely equivalent to the initial
condition of variationally triviality.

(ii) Next, one sees that from (\ref{lambda-nu}) it follows in particular that
the $n$-form $d\nu$ is $\pi^{r+1.r}$-projectable. We can apply proposition
\ref{nu-nu-tilde} and obtain that $\nu = h\tilde{\nu}_{V}$ for some basic form 
in the chart
$(V^{r-1},\psi^{r-1})$.
If we define
\be
\rho_{V} \equiv d\tilde{\nu}_{V}
\label{ro-V}
\ee
then we obtain after some computation that 
$$
(\pi^{r+1,r})^{*}(\lambda - h\rho_{V}) = 0
$$
and (1) from the statement follows. The definition (\ref{ro-V}) guarantees that
we also have (2).
$\qed$

\begin{rem}
A statement of the type appearing in this theorem is, in fact, valid for every
$\lambda \in \Omega^{r}_{q,X} \quad (q \leq n)$
such that
$E_{q}([\lambda]) = 0$.
\end{rem}

The theorem we just have proved has the following consequence (see also 
\cite{AD} thm. 4.3)

\begin{cor}
Any variationally trivial Lagrange form
$\lambda \in \Omega^{r}_{n,X}$
can be locally written as a total exterior derivative of a local form
$\omega_{V} \in {\cal J}^{r}_{n-1}$:
\be
\lambda = D\omega_{V}.
\ee
\end{cor}

{\bf Proof:}
In the proof of the preceding theorem we restrict, eventually, the chart
$V^{r-1}$
and we have
$\rho_{V} = d\eta_{V}$
for some $(n-1)$-form on
$V^{r-1}$.
Then, according to the definition (\ref{D}) of the total exterior derivative we
have the formula from the statement with 
$\omega_{V} \equiv h\eta_{V} \in {\cal J}^{r}_{n-1}$.
$\qed$

One can now obtain the most general form of a variationally trivial local 
Lagrangian. 

\begin{thm}
Any variationally trivial local Lagrangian of order $r$ has the following form
in the chart
$(V^{r},\psi^{r})$:
\be
L = \sum_{s=0}^{n} {1\over s!(n-s)!} \sum_{|I_{1}|,...,|I_{s}|=r-1} 
{\cal L}^{I_{1},...,I_{s}}_{\sigma_{1},...,\sigma_{s},i_{s+1},...,i_{n}}
{\cal J}_{I_{1},...,I_{s}}^{\sigma_{1},...,\sigma_{s},i_{s+1},...,i_{n}}.
\label{L-hiper}
\ee

Here we have defined
\be
{\cal J}_{I_{1},...,I_{s}}^{\sigma_{1},...,\sigma_{s},i_{s+1},...,i_{n}} 
\equiv \varepsilon^{i_{1},...,i_{n}} 
\prod_{l=1}^{s} y^{\sigma_{l}}_{I_{l}i_{l}} \quad (s = 0,...,n)
\label{hyperJacobians}
\ee

and the function
${\cal L}^{...}_{...}$
are given by
\be
{\cal L}^{I_{1},...,I_{s}}_{\sigma_{1},...,\sigma_{s}i_{s+1},...,i_{n}} \equiv
\sum_{k=1}^{s} (-1)^{k-1} \partial^{I_{k}}_{\sigma_{k}} 
A^{I_{1},...,\hat{I_{k}},...,I_{s}}_{\sigma_{1},...,\hat{\sigma_{k}},...,
\sigma_{s},i_{s+1},...,i_{n}} + 
\sum_{k=s+1}^{q}  (-1)^{k-1} d_{i_{k}} A^{I_{1},...,I_{s}}_{\sigma_{1},...,
\sigma_{s},i_{s+1},...,\hat{i_{k}},...,i_{n}};
\label{AL}
\ee
the expressions
$A^{I_{1},...,I_{s}}_{\sigma_{1},...,\sigma_{s}i_{s+1},...,i_{n-1}}$
are smooth functions on 
$V^{r-1}$
verifying the symmetry property  (\ref{symmetry}) and
$d_{j} = d_{j}^{r-1}$
is the corresponding formal derivative on
$V^{r-1}$.
\end{thm}

{\bf Proof:}
It is convenient to introduce the following forms:
$$
\chi \equiv \sum_{s=0}^{n-1} {1 \over s!(n-1-s)!} 
\sum_{|I_{1}|,...,|I_{s}|=r-1}
A^{I_{1},...,I_{s}}_{\sigma_{1},...,\sigma_{s},i_{s+1},...,i_{n-1}}
dy^{\sigma_{1}}_{I_{1}} \wedge \cdots \wedge dy^{\sigma_{s}}_{I_{s}} \wedge
dx^{i_{s+1}} \wedge \cdots \wedge dx^{i_{n-1}}
$$
and
$$
\theta \equiv \sum_{s=0}^{n} {1 \over s!(n-s)!} 
\sum_{|I_{1}|,...,|I_{s}|=r-1}
{\cal L}^{I_{1},...,I_{s}}_{\sigma_{1},...,\sigma_{s},i_{s+1},...,i_{n}}
dy^{\sigma_{1}}_{I_{1}} \wedge \cdots \wedge dy^{\sigma_{s}}_{I_{s}} \wedge
dx^{i_{s+1}} \wedge \cdots \wedge dx^{i_{n}}.
$$

Then one finds out by direct computation that
$$\theta = d\chi + {\rm contact~ terms}.$$

Next, one takes in theorem \ref{trivial} 
$\tilde{\nu}_{V} = \chi$ and it follows that 
$\rho_{V} = \theta + {\rm contact~ terms}$. 
Finally, by direct computation one discovers that 
$\lambda = h\rho_{V}$
has the expression (\ref{L-hiper}).
$\qed$

\begin{rem}
Let us note that the expressions (\ref{AL}) are of the same type as those 
given by (\ref{tilde-A}). 
\end{rem}

The expressions (\ref{hyperJacobians}) defined above are called {\it
hyper-Jacobians} \cite{BCO} \cite{Ol} (see these references for similar 
results). It is immediate that they have antisymmetry properties of the type 
(\ref{symmetry}).   

Now we give another argument for the converse statement from the preceding 
theorem is true. First we have: 

\begin{cor}
The local expression of a variationally trivial Lagrangian (\ref{L-hiper})
can be rewritten as follows:
\be
L = d_{j} V^{j}
\ee
where 
$V^{j}$
are some smooth functions on 
$V^{r}$.
\end{cor}

{\bf Proof:}
Using the notations introduced in the proof above let us define the local 
expressions on
$V^{r}$:
$$
V^{j} \equiv \varepsilon^{j,i_{1},...,i_{n-1}} \sum_{s=0}^{n} 
{1 \over s!(n-1-s)!} \sum_{|I_{1}|,...,|I_{s}|=r-1}
A^{I_{1},...,I_{s}}_{\sigma_{1},...,\sigma_{s},i_{s+1},...,i_{n-1}}
y^{\sigma_{1}}_{I_{1}i_{1}} \cdots y^{\sigma_{s}}_{I_{s}i_{s}}. 
$$
One checks now that the formula from the statement is true.
$\qed$

Now we indeed have:

\begin{thm}
The expression (\ref{L-hiper}) is  variationally trivial.
\end{thm}

{\bf Proof:}
We have according to the preceding corollary
$E_{\sigma}(L) = E_{\sigma}(d_{j} V^{j}) = 0$
because we can apply  (\ref{tot-div}).
$\qed$

\begin{rem}
Some globalisation of the results above can be found in \cite{Kr5} (see
theorem 5 and corollary 1 from this reference). In particular, for
$r = 1$ one obtains known results \cite{Ru}, \cite{Ed}, \cite{Kr9}, \cite{GP}.
\end{rem}

We close this Section with the analysis of the following problem. We have 
discovered that a variationally trivial Lagrangian depends on the highest order
derivatives through some very particular polynomial expressions. The problem is
to obtain a system of partial differential equations which is compatible only
with this structure. More precisely, we have the following result.

\begin{thm}

Let us suppose that the local Lagrangian 
$L$ on $V^{r}$
verifies the system of partial differential equations:
\be
{\cal S}^{+}_{p_{1},...,p_{r},j_{r}} \partial^{p_{1},...,p_{r}}_{\rho}
\partial^{j_{1},...,j_{r}}_{\sigma} L  = 0.
\label{system}
\ee

Then $L$ is a polynomial in
$y^{\sigma}_{I}, \quad |I| = r$
of the following form:
\be
L = \sum_{s=0}^{n} {1\over s!(n-s)!} \sum_{|I_{1}|,...,|I_{s}|=r-1} 
{\cal L}^{I_{1},...,I_{s}}_{\sigma_{1},...,\sigma_{s},i_{s+1},...,i_{n}}
{\cal J}_{I_{1},...,I_{s}}^{\sigma_{1},...,\sigma_{s},i_{s+1},...,i_{n}}
\label{LJ}
\ee
where 
${\cal L}^{I_{1},...,I_{s}}_{\sigma_{1},...,\sigma_{s},i_{s+1},...,i_{n}}$
are some smooth functions on 
$V^{r-1}$
verifying symmetry properties of the type (\ref{symmetry}). Conversely, if
$L$ is of the form above, then the system (\ref{system}) is identically
fulfilled. 
\end{thm}

{\bf Proof:}

(i) We will first prove that if the local Lagrangian $L$ verifies the
identities: 
\be
{\cal S}^{+}_{p_{1},...,p_{r-l+1},j_{k-l+1},...,j_{r}} 
\partial^{p_{1},...,p_{r-l+1}}_{\rho}
\partial^{j_{1},...,j_{k}}_{\sigma} L  = 0 \quad ( 1\leq l \leq k \leq r)
\label{syst1}
\ee
and
\be
\sum_{k=0}^{r} (-1)^{k} d_{j_{1}}^{r} \cdots d_{j_{k}}^{r} 
\partial^{j_{1},...,j_{k}}_{\sigma} L = 0
\label{syst2}
\ee
then we have 
$E_{\sigma}(L) = 0$
in $V^{s} \quad (s > 2r).$

Indeed one starts directly form the definition
\be
E_{\sigma}(L) = \sum_{k=0}^{r} (-1)^{k} d_{j_{k}}^{s} \cdots d_{j_{1}}^{s} 
\partial^{j_{1},...,j_{k}}_{\sigma} L = 
\sum_{k=0}^{r} (-1)^{k} (d_{j_{k}}^{r} + \sum_{l=r}^{2r-1} 
y^{\nu}_{p_{1},...,p_{l},k_{k}} \partial^{p_{1},...,p_{l}}_{\nu})
d_{j_{k-1}}^{s} \cdots  d_{j_{1}}^{s} \partial^{j_{1},...,j_{k}}_{\sigma} L 
\label{ind}
\ee
and commutes
$\partial^{p_{1},...,p_{l}}_{\nu}$
over
$d_{j_{k-1}}^{s} \cdots  d_{j_{1}}^{s}$.

We must use the hypothesis (\ref{syst1}) and a generalization of the formula 
(\ref{commutator}), namely 
\be
\left[ \partial^{p_{1},...,p_{k}}_{\sigma}, d_{j_{1}} \cdots d_{j_{l}}\right] =
{\cal S}^{+}_{p_{1},...,p_{k}} {\cal S}^{+}_{j_{1},...,j_{l}} 
\sum_{t=1}^{l} c_{k,l,t} \delta^{p_{1}}_{j_{1}} \cdots \delta^{p_{t}}_{j_{t}} 
d_{j_{t+1}} \cdots d_{j_{l}}  \partial^{p_{t+1},...,p_{k}}_{\sigma}
\quad (k \geq l)
\ee
where 
$c_{k,l,t} \in \R_{+}$;
this formula can be proved by induction on $l$.

Now we easily obtain that in fact the terms corresponding to the sum over $l$ 
in (\ref{ind}) give a null contribution, so we are left with:
$$
E_{\sigma}(L) = \sum_{k=0}^{r} (-1)^{k} d_{j_{k}}^{r} d_{j_{k-1}}^{s} \cdots  
d_{j_{1}}^{s} \partial^{j_{1},...,j_{k}}_{\sigma} L.
$$

We continue in the same way by recurrence and finally get
$$
E_{\sigma}(L) = \sum_{k=0}^{r} (-1)^{k} d_{j_{1}}^{r} \cdots d_{j_{k}}^{r} 
\partial^{j_{1},...,j_{k}}_{\sigma} L 
$$
which is zero, according to (\ref{syst2}).

It appears that (\ref{syst1}) and (\ref{syst2}) imply that the Lagrangian $L$
is variationally trivial. According to theorem \ref{trivial}, $L$ has the form
(\ref{LJ}) from the statement, but with some restrictions on the functions
${\cal L}^{...}_{...}$. On the other hand, it is clear that the dependence on
the highest order derivatives must follow only from the equations containing
only the partial derivatives of order $r$ i.e. (\ref{syst1}) for 
$k = r$
and
$l = 1$
i.e. the system (\ref{system}) from the statement.

(ii) The converse statement follows rather easily. One shows by direct
computations that if $L$ is given by the formula (\ref{LJ}), then:
$$
\partial^{I_{1}}_{\sigma_{1}} L = 
\sum_{s=1}^{n} (-1)^{s-1} {1\over (s-1)!(n-s)!} \sum_{|I_{2}|,...,|I_{s}|=r-1} 
(B_{1}{\cal L})^{I_{1},...,I_{s}}_{\sigma_{1},...,\sigma_{s},i_{s+1},...,i_{n}}
{\cal J}_{I_{1},...,I_{s}}^{\sigma_{1},...,\sigma_{s},i_{s+1},...,i_{n}}.
$$

Iterating the derivation procedure we have
$$
\partial^{I_{1}}_{\sigma_{1}} \partial^{I_{2}}_{\sigma_{2}}L = -
\sum_{s=2}^{n} {1\over (s-2)!(n-s)!} \sum_{|I_{3}|,...,|I_{s}|=r-1} 
(B_{1}B_{2}{\cal L})^{I_{1},...,I_{s}}_{\sigma_{1},...,\sigma_{s},
i_{s+1},...,i_{n}}
{\cal J}_{I_{1},...,I_{s}}^{\sigma_{1},...,\sigma_{s}i_{s+1},...,i_{n}}
$$
and one proves that (\ref{system}) is fulfilled by direct computation.
$\qed$
\newpage

\section{Locally Variational Differential Equations}

The definitions for a general differential equation and for a locally
variational differential equation have been given previously (see the formul\ae
\quad (\ref{T}) and resp. (\ref{EL-type}).) We want to analyse the general 
structure of a locally variational differential equation along the same lines 
of argument as in the previous Section. We will not be able to obtain the most 
general expression for such an object (as we have been able to obtain in the 
case of variationally trivial Lagrangians) but we will produce a generic 
expression of the same type as (\ref{LJ}). We mention again that this result 
has been already obtained in \cite{An} with a completely different method.

Our starting point is a analogue of the proposition \ref{eta-nu} for the case
$q = n$; it is natural to also make the replacement
$h \rightarrow p_{1}.$

\begin{prop}
Let 
$\eta \in \Omega^{r}_{n}$
such that
$p_{1}d\eta$
is a $\pi^{r+1,r}$-projectable $(n+1)$-form. Then one can write $\eta$ as 
follows:
\be
\eta = \nu + \sum_{|J| \leq r-1} \Phi^{J}_{\sigma,i_{1},...,i_{n-1}} 
\omega^{\sigma}_{J} \wedge dx^{i_{1}} \wedge \cdots \wedge dx^{i_{n-1}} + 
d\phi + \psi
\ee
where
$\nu \in \Omega^{r}_{n,X}$ 
is a basic $n$-form,
$\psi \in \Omega^{r}_{n+1}$
is a 2-contact form and
$
\Phi^{J}_{\sigma,i_{1},...,i_{n-1}} 
$
are smooth functions on $V^{r}$ completely antisymmetric in the indices
$i_{1},...,i_{n-1}.$
\label{eta-nu-n}
\end{prop}

{\bf Proof:}

We proceed in analogy with proposition \ref{eta-nu}. Let
$(V,\psi)$ be a chart on $Y$. Then in the associated chart
$(V^{r},\psi^{r})$
one can write $\eta$
in the standard form
\be
\eta = \eta_{0} + \eta_{1}
\label{eta-n}
\ee
where in $\eta_{0}$ we collect all terms containing at least one factor
$\omega^{\sigma}_{J} \quad (|J| \leq r-1)$
and $\eta_{1}$ is a polynomial only in 
$dx^{i}$ and $dy^{\sigma}_{J} \quad (|J| = r+1)$:
\be
\eta_{1} = \sum_{s=0}^{n}  \sum_{|I_{1}|,...,|I_{s}|=r}
A^{I_{1},...,I_{s}}_{\sigma_{1},...,\sigma_{s},i_{s+1},...,i_{n}}
dy^{\sigma_{1}}_{I_{1}} \wedge \cdots \wedge dy^{\sigma_{s}}_{I_{s}} \wedge
dx^{i_{s+1}} \wedge \cdots \wedge dx^{i_{n}},
\label{eta1-n}
\ee
where the coefficients
$A^{I_{1},...,I_{s}}_{\sigma_{1},...,\sigma_{s},i_{s+1},...,i_{n}}$
have antisymmetry properties of type (\ref{symmetry}). 

The generic form of $\eta_{0}$ is
\be
\eta_{0} = \sum_{|J| \leq r-1} \omega^{\sigma}_{J} \wedge \Phi^{J}_{\sigma}
\label{eta0}
\ee
with
$\Phi^{J}_{\sigma}$
some $(n-1)$-forms. It follows then that
\be
d\eta = d\eta_{1} + \sum_{|J| \leq r-1} \left[ (d\omega^{\sigma}_{J}) \wedge
\Phi^{J}_{\sigma}- \omega^{\sigma}_{J} \wedge (d\Phi^{J}_{\sigma}) \right].
\label{d-eta-n}
\ee

One applies the operator $p_{1}$ to this equality and uses (\ref{p-wedge}); a 
new relation is obtained, namely:
\be
p_{1} d\eta = p_{1} d\eta_{1} + \sum_{|J| \leq r-1} 
\left[ (d\omega^{\sigma}_{J}) \wedge h\Phi^{J}_{\sigma} - 
\omega^{\sigma}_{J} \wedge (hd\Phi^{J}_{\sigma})\right].
\label{p1-d-eta-n}
\ee

The expression of 
$d\eta_{1}$
has been computed quite generally before and is given by (\ref{d-eta1-q}). 
Applying the operator $p_{1}$ and using (\ref{p-k-ro}) one obtains:
\begin{eqnarray}
p_{1} d\eta_{1} = \sum_{s=0}^{n} \sum_{|I_{1}|,...,|I_{s}|=r}
\sum_{|J| \leq r-1} \left( \partial^{J}_{\nu} 
A^{I_{1},...,I_{s}}_{\sigma_{1},...,\sigma_{s},i_{s+1},...,i_{n}} \right)
y^{\sigma_{1}}_{I_{1}} \cdots y^{\sigma_{s}}_{I_{s}} 
\omega^{\nu}_{J} \wedge dx^{i_{s+1}} \wedge \cdots \wedge dx^{i_{n}} 
\nonumber \\
+ \sum_{|I_{1}|=r} \tilde{B}^{I_{1}}_{\sigma_{1},i_{2},...,i_{n+1}}
\omega^{\sigma_{1}}_{I_{1}} 
\wedge dx^{i_{2}} \wedge \cdots \wedge dx^{i_{n+1}}, \nonumber
\label{p1-eta1-n}
\end{eqnarray}
where 
\be
\tilde{B}^{I_{1}}_{\sigma_{1},i_{2},...,i_{n+1}} =
{\cal A}_{i_{2},...,i_{n+1}} \sum_{s=1}^{n+1} s
\sum_{|I_{1}|,...,|I_{s}|=r} 
\tilde{A}^{I_{1},...,I_{s}}_{\sigma_{1},...,\sigma_{s},i_{s+1},...,i_{n+1}}
y^{\sigma_{2}}_{I_{2}i_{2}}\cdots y^{\sigma_{s}}_{I_{s}i_{s}}.
\label{p1-n}
\ee

The expression (\ref{p1-d-eta-n}) becomes
\be
p_{1} d\eta = - \sum_{|J|=r-1} \omega^{\sigma}_{Ji} \wedge dx^{i} \wedge 
h\Phi^{J}_{\sigma} + \sum_{|I_{1}|=r} 
\tilde{B}^{I_{1}}_{\sigma_{1},i_{2},...,i_{n+1}} \omega^{\sigma_{1}}_{I_{1}} 
\wedge dx^{i_{2}} \wedge \cdots \wedge dx^{i_{n+1}} + \cdots , \nonumber
\label{p1-r}
\ee
where by $\cdots$ we mean contributions which do not contain the differentials
$\omega^{\sigma}_{I}, \quad |I| = r.$

The first term has the generic form 
$$
\sum_{|I_{1}|=r} C^{I_{1}}_{\sigma_{1},i_{2},...,i_{n+1}} 
\omega^{\sigma_{1}}_{I_{1}} \wedge dx^{i_{2}} \wedge \cdots \wedge dx^{i_{n+1}}
$$
where 
$
C^{I_{1}}_{\sigma_{1},i_{2},...,i_{n+1}}
$
are some polynomials in
$y^{\sigma}_{I}, \quad |I| = r+1$
of maximal degree $(n-1)$;
because of the presence of the combination
$
\omega^{\sigma}_{Ji} \wedge dx^{i}
$
these polynomials are of delta-type i.e. they are obtained by applying 
$B_{1}$ on some other tensors; it follows that we have in compact tensor
notations: 
\be
B_{1}C_{\sigma} = 0.
\label{B1-C}
\ee

The expression (\ref{p1-r}) does not depend on 
$y^{\sigma}_{I}, \quad |I| = r+1$
by hypothesis. This is equivalent with the independence of 
$\tilde{B}_{\sigma} - C_{\sigma}$
on
$y^{\sigma}_{I}, \quad |I| = r+1.$
In particular, the same thing must be true for 
$
B_{1} (\tilde{B}_{\sigma} - C_{\sigma}) = B_{1} \tilde{B}_{\sigma}
$
where use of (\ref{B1-C}) has been made. 

If the expression (\ref{p1-n}) is used one obtains 
$$
{\cal S}^{-}_{i_{2},...,i_{n+1}} 
(B_{1}\tilde{A})^{I_{1},...,I_{s}}_{\sigma_{1},...,\sigma_{s},i_{s+1},...,
i_{n+1}}
y^{\sigma_{2}}_{I_{2}i_{2}}\cdots y^{\sigma_{s}}_{I_{s}i_{s}} = 0, \quad
(s = 2,...,n+1)
$$
or, in compact tensor notations:
$$
B_{2} \cdots B_{s} B_{1} \tilde{A}_{\sigma_{1},...,\sigma_{s}} = 0, \quad
(s = 2,...,n+1).
$$

One can apply lemma \ref{BXB} and obtains that
$\tilde{A}_{\sigma_{1},...,\sigma_{s}}$
has the following form
$$
\tilde{A}_{\sigma_{1},...,\sigma_{s}} = \sum_{\alpha =1}^{s} B_{\alpha}
A^{\alpha}_{\sigma_{1},...,\sigma_{s}}, \quad (s = 2,...,n+1).
$$

This expression must be substituted into the formula for
$d\eta_{1}$
and, after some elementary prelucrations, the following generic form is
produced: 
\be
d\eta_{1} = \nu_{1} + d\phi + \sum_{|J| \leq r-1} \omega^{\sigma}_{J} \wedge 
\Phi^{J}_{\sigma}.
\label{d-eta1-n}
\ee
Here
$\phi \in \Omega^{r}_{n(c)}$
is a contact form,
$\nu_{1}$
has the form
\be
\nu_{1} = \sum_{|I_{1}|=r} A^{I_{1}}_{\sigma_{1},i_{2},...,i_{n+1}}
dy^{\sigma_{1}}_{I_{1}} \wedge dx^{i_{2}} \wedge \cdots \wedge dx^{i_{n+1}}
\label{nu1-n}
\ee
and
$\Phi^{J}_{\sigma}$
is a polynomial of degree $n$ in
$dx^{i}$ and $dy^{\sigma}_{I}, \quad |I| = r.$
Explicitly:
\be
\Phi^{J}_{\sigma} = \sum_{s=0}^{n}  \sum_{|I_{1}|,...,|I_{s}|=r}
\Phi^{J,I_{1},...,I_{s}}_{\sigma,\nu_{1},...,\nu_{s},i_{s+1},...,i_{n}}
dy^{\nu_{1}}_{I_{1}} \wedge \cdots \wedge dy^{\nu_{s}}_{I_{s}} \wedge
dx^{i_{s+1}} \wedge \cdots \wedge dx^{i_{n}};
\label{Phi-J-sigma-n}
\ee
here we can assume a (partial) symmetry property of the type 
(\ref{sym-partial}) and we make the convention 
\be
\Phi^{J,I_{1},...,I_{s}}_{\sigma,\nu_{1},...,\nu_{s},i_{s+1},...,i_{n}} = 0,
\quad \forall J \quad s.t. \quad |J| \geq r.
\ee

(ii) From the expression (\ref{d-eta1-n}) one obtains after exterior
differentiation and application of the operator $p_{2}$ the following
condition: 
\be
p_{2}d\nu_{1} - \sum_{|J| \leq r-1} (\omega^{\sigma}_{Ji} \wedge dx^{i} \wedge
p_{1}\Phi^{J}_{\sigma} - \omega^{\sigma}_{J} \wedge p_{1}d\Phi^{J}_{\sigma})
= 0.
\label{d-eta1-n-a}
\ee
 
As in the proof of the proposition \ref{eta-nu}, one finds out from this 
equation the following consequences:
\begin{eqnarray}
{\cal S}^{-}_{i_{2},...,i_{n+1}} {\cal S}^{+}_{I_{1}p_{1}} \cdots 
{\cal S}^{+}_{I_{s}p_{s}} 
\left[ \tilde{\Phi'}^{\{j_{1},...,j_{k}\},I_{1},...,I_{s}}_{\sigma,\nu_{1},...,
\nu_{s},i_{s+1},...,i_{n+1}} + (-1)^{s} {\cal S}_{j_{1},...,j_{k}}
\delta^{j_{1}}_{i_{s+1}} \Phi^{\{j_{2},...,j_{k}\},I_{1},...,I_{s}}_{\sigma,
\nu_{1},...,\nu_{s},i_{s+2},...,i_{n+1}} \right] \times \nonumber \\
\delta^{p_{2}}_{i_{2}} \cdots \delta^{p_{s}}_{i_{s}} = 0, 
\quad (s = 2,...,n+1, k = 0,...,r-1)
\label{dnu1-n-a}
\end{eqnarray}
and
\begin{eqnarray}
{\cal S}^{-}_{i_{2},...,i_{n+1}} {\cal S}^{+}_{I_{1}p_{1}} \cdots 
{\cal S}^{+}_{I_{s}p_{s}} {\cal S}^{+}_{j_{1},...,j_{k}} 
\delta^{j_{1}}_{i_{s+1}} \left( 
\Phi^{\{j_{2},...,j_{r}\},I_{1},...,I_{s}}_{\sigma,\nu_{1},...,\nu_{s},
i_{s+2},...,i_{n+1}} +
\Phi^{I_{1},\{j_{2},...,j_{k}\},I_{2},...,I_{s}}_{\nu_{1},\sigma,\nu_{2},...,
\nu_{s},i_{s+2},...,i_{n+1}} \right) \times
\nonumber \\
\delta^{p_{2}}_{i_{2}} \cdots \delta^{p_{s}}_{i_{s}} = 0,
\quad (s = 2,...,n+1).
\label{dnu1-n-b}
\end{eqnarray}

(Here use have been made of the definition (\ref{tilde-A'}).)

In compact notations, these two relations can be written in a similar way and
with similar notations (using, in particular, the convenient operator 
(\ref{delta'})):
\be
B_{2} \cdots B_{s} \left((\Delta' \Phi)^{k}_{s} + B_{0} \Phi^{k-1}_{s}\right) 
= 0, \quad (s = 2,...,n+1, k = 0,...,r-1)
\label{descent-n-a}
\ee
and
\be
B_{2} \cdots B_{s} \left(B_{0} \Phi^{r-1}_{s} + B_{1} \Psi^{r-1}_{s}\right) 
= 0, \quad (s = 2,...,n+1);
\label{descent-n-b}
\ee
here we have defined:
\be
\Psi^{I_{0},I_{1},...,I_{s}}_{\nu_{0},\nu_{1},...,\nu_{s},i_{s+1},...,i_{n+1}}
\equiv
\Phi^{I_{1},I_{0},...,I_{s}}_{\nu_{1},\nu_{0},...,\nu_{s},i_{s+1},...,i_{n+1}}.
\ee

It is clear that we must apply again the descent procedure from the proof of
the proposition \ref{eta-nu} in a case which is a little more complicated. We
first deal with the equation (\ref{descent-n-b}) by applying the operator
$B_{1}$; then one gets:
$$
B_{0} \cdots B_{s} \Phi^{r-1}_{s} = 0, \quad (s = 2,...,n+1) 
$$
from where one obtains, with lemma \ref{BXB}
$$
\Phi^{r-1}_{s} = \sum_{\alpha =0}^{s} B_{\alpha} \Phi^{r-1,\alpha}_{s}. 
\quad (s = 2,...,n+1).
$$

As before, one can redefine 
$\Phi^{r-1}_{s}$
such that, instead of the previous formula we have a generic expression of the
type: 
\be
\Phi^{r-1}_{s} = B_{0} C^{r-1}_{s}, \quad (s = 2,...,n+1)
\label{n-(r-1)}
\ee
and, instead of (\ref{d-eta1-n}) we have
\be
d\eta_{1} = \nu_{1} + d\phi + \sum_{|J| \leq r-1} \omega^{\sigma}_{J} \wedge 
\Phi^{J}_{\sigma} + \psi
\label{d-eta1-n-b}
\ee
where $\psi$ is a 2-contact form.

One can take from the beginning the previous argument, based on acting on this
relation with the operator $p_{2}d$, and obtains this time only the relation 
(\ref{descent-n-a}) together with (\ref{n-(r-1)}); the tensors 
$C^{r-1}_{s}$ 
stay arbitrary. Now we are back in the same case we have previously analysed in
proposition \ref{eta-nu}. The descent procedure can be applied and instead of
(\ref{d-eta1-n-b}) one gets the following expression:
\begin{eqnarray}
d\eta_{1} = \nu_{1} + d\phi + \psi +
\sum_{|J| \leq r-1} \Phi^{J}_{\sigma,i_{1},...,i_{n}} 
\omega^{\sigma}_{J} \wedge dx^{i_{1}} \wedge \cdots \wedge dx^{i_{n}} +
\nonumber \\
\sum_{|J| \leq r-1} \sum_{|I_{1}|=r}
\Phi^{JI_{1}}_{\sigma,\nu_{1},i_{2},...,i_{n}} \omega^{\sigma}_{J} \wedge 
dy^{\nu_{1}}_{I_{1}} \wedge dx^{i_{2}} \wedge \cdots \wedge dx^{i_{n}} 
\label{d-eta1-n-c}
\end{eqnarray}
where the notations have the same meaning as before. 

Now we apply for the third time the argument based on acting on this relation 
with the operator $p_{2}d$, and obtain instead of (\ref{descent-n-a}) and
(\ref{descent-n-b}) only the equation
$$
B_{2}(\Delta'\Phi^{k})_{2} = 0, \quad (k = 0,...,r-1).
$$

From symmetry considerations we also have the same relation with
$B_{2} \rightarrow B_{1}$ 
so, we can obtain as in the proof of proposition \ref{nu-nu-tilde}, that the 
following equality stays true:
\be
\Delta'\Phi^{k}_{1} = B_{1} B_{2} C^{k}, \quad (k = 0,...,,r-1)
\label{d-phi}
\ee
where
$C^{k} \in {\cal H}_{n-1,k-2,r-1,r-1}.$

(iii) We concentrate on the analysis of the previous equation (\ref{d-phi}).
The key observation is that there exists a homotopy operator for $\Delta'$. 
Indeed, we have 
\begin{lemma}
Suppose that the tensor 
$A = \oplus_{s=1}^{q} A_{s}$ 
verifies the equation:
\be
\Delta' A = 0.
\ee

Then there exists a tensor
$C = \oplus_{s=0}^{q-1} C_{s}$
such that
\be
A = \Delta' C.
\ee
\end{lemma}

{\bf Proof:}
 
An explicit formula for the tensor $C$ (suggested by \cite{Kr5}) is the 
following one. One defines the map
$\chi: \R \times V^{r} \rightarrow V^{r}$
by
$$
\chi(u,(x^{i},y^{\sigma},y^{\sigma}_{j},...,y^{\sigma}_{j_{1},...,j_{r}})) =
(x^{i},y^{\sigma},y^{\sigma}_{j},...,uy^{\sigma}_{j_{1},...,j_{r}})
$$
and afterwards:
$$
C^{J,I_{1},...,I_{s}}_{\nu,\sigma_{1},...,\sigma_{s},i_{s+1},...,i_{q-1}} 
\equiv 
(s+1) \sum_{|I_{0}|=r} y^{\sigma_{0}}_{I_{0}} \int_{0}^{1} u^{s}
A^{J,I_{0},...,I_{s}}_{\nu,\sigma_{0},...,\sigma_{s},i_{s+1},...,i_{q-1}} 
\circ \chi \quad du, \quad (s=1,...,q-1).
$$

By elementary computations one finds out that:
$$
(\Delta' C)^{J,I_{1},...,I_{s}}_{\nu,\sigma_{1},...,\sigma_{s},
i_{s+1},...,i_{q}} =
A^{J,I_{1},...,I_{s}}_{\nu,\sigma_{1},...,\sigma_{s},i_{s+1},...,i_{q}},
\quad (s = 1,...,q)
$$
and this gives us the desired homotopy formula.
$\nabla$

The previous lemma has a generalization. 
\begin{lemma}
Let us suppose that we have the equation: 
\be
\Delta' A = D.
\ee

Then tensor $D$ verifies the consistency relation: 
\be
\Delta' D = 0
\ee
and a particular solution of the equation above is of the form:
\be
A^{J,I_{1},...,I_{s}}_{\nu,\sigma_{1},...,\sigma_{s}i_{s+1},...,i_{q-1}} =
(s+1) \sum_{|J|=r} y^{\sigma_{0}}_{I_{0}} \int_{0}^{1} u^{s}
D^{J,I_{0},...,I_{s}}_{\nu,\sigma_{0},...,\sigma_{s},i_{s+1},...,i_{q-1}} 
\circ \chi \quad du, \quad (s=1,...,q-1).
\ee
\end{lemma}

{\bf Proof:} The consistency condition follows from (\ref{delta2}) and the last
equality by direct computations.
$\nabla$

(iv) We need the preceding two lemmas only in the case $q = n$. One finds out
that the generic solution of the equation (\ref{d-phi}) is
\be
\Phi^{k}_{1} = \Delta'\Psi^{k} + B_{1} A^{k}, \quad (k = 0,...,r-1)
\ee
where the first term is an arbitrary solution of the homogeneous equation
(obtained with the first lemma)  and the second term is a particular solution
of the non-homogeneous equation (obtained with the second lemma for
$D = B_{1} B_{2} C^{k}$.) 

Now we must substitute this result in the last sum from the formula 
(\ref{d-eta1-n-c}); it is not very hard to regroup the result with the first
sum such that one obtains a more simple expression:
\be
d\eta_{1} = \nu_{1} + d\phi + \psi +
\sum_{|J| \leq r} A^{J}_{\sigma,i_{1},...,i_{n}} 
dy^{\sigma}_{J} \wedge dx^{i_{1}} \wedge \cdots \wedge dx^{i_{n}}.
\label{d-eta1-n-d}
\ee

(v) The formula above can be prelucrated, as in proposition \ref{eta-nu}, using
the contact homotopy operator $A$. Indeed, if we apply to the form
$\eta_{1} - \phi_{1}$
the formula (\ref{a-ro}) we obtain:
\be
\eta_{1} - \phi_{1} = A \left( \nu_{1} + \sum_{|J| \leq r} 
A^{J}_{\sigma,i_{1},...,i_{n}} dy^{\sigma}_{J} 
\wedge dx^{i_{1}} \wedge \cdots \wedge dx^{i_{n}} + \psi \right) +
dA(\eta_{1} - \phi_{1}) + \eta_{1}(0) - \phi_{1}(0).
\label{eta-A-n}
\ee

As before we easily establish the following facts:

- $A\nu_{1} = 0$;

- $A( \sum \cdots )$ and $\eta_{1}(0)$ are basic $n$-forms;

- if $\phi$ is contact, then $\phi_{1}(0)$ is contact;

-  if $\psi$ is a 2-contact form, then $A\psi$ is a contact form.

These results, gives us the following generic expression for $\eta_{1}$:
\be
\eta_{1} = \nu + \eta_{0} + d\phi
\ee
where $\nu$ is a basic $n$-form, $\phi$ is a contact $(n-1)$-form and 
$\eta_{0}$ is a $n$-contact form of the type (\ref{eta0}).
It is obvious that the same type of expression stays true for $\eta$ also.

(vi) It is useful to notice that the forms 
$\Phi^{J}_{\sigma}$
from (\ref{eta0}) can also be decomposed into a contribution having at least a
factor 
$\omega^{\sigma}_{J}$
and a polynomial of degree $(n-1)$ in the differentials
$dx^{i}$ and $dy^{\sigma}_{I}, \quad |I| = r.$
In this way, the formula for $\eta$ takes a more convenient form, namely:
\be
\eta = \nu + \eta_{0} + d\phi + \psi
\label{vi}
\ee
where $\nu$ and $\phi$ have the same properties as above, $\psi$ is a 2-contact
form and $\eta_{0}$ has the expression of the type (\ref{eta0}) but the forms 
$\Phi^{J}_{\sigma}$
are polynomials (of degree $n-1$) only in the differentials
$dx^{i}$ and $dy^{\sigma}_{I}, \quad |I| = r.$
Explicitly:
\be
\Phi^{J}_{\sigma} = \sum_{s=0}^{n-1}  \sum_{|I_{1}|,...,|I_{s}|=r}
\Phi^{J,I_{1},...,I_{s}}_{\sigma,\nu_{1},...,\nu_{s},i_{s+1},...,i_{n-1}}
dy^{\nu_{1}}_{I_{1}} \wedge \cdots \wedge dy^{\nu_{s}}_{I_{s}} \wedge
dx^{i_{s+1}} \wedge \cdots \wedge dx^{i_{n-1}}
\label{Phi-J-sigma-n-1}
\ee
where symmetry properties of the type (\ref{sym-partial}) can be imposed and 
we admit that: 
\be
\Phi^{J,I_{1},...,I_{s}}_{\sigma,\nu_{1},...,\nu_{s},i_{s+1},...,i_{n}} = 0,
\quad \forall J \quad s.t. \quad |J| \geq r.
\ee

We must impose on this expression for $\eta$ the condition from the statement 
of the proposition. The local expression for $\nu$ is given by (\ref{L}).
We obtain instead of (\ref{p1-d-eta-n}):
\be
p_{1} d\eta = \sum_{|J| \leq r} (\partial^{J}_{\sigma} L) \omega^{\sigma}_{J} 
\wedge \theta_{0} + \sum_{|J| \leq r-1} 
\left[ (d\omega^{\sigma}_{J}) \wedge h\Phi^{J}_{\sigma} - 
\omega^{\sigma}_{J} \wedge (hd\Phi^{J}_{\sigma})\right].
\label{p1-d-eta-n-1}
\ee

The condition of $\pi^{r+1,r}$-projectability for this expression amounts to a
set of relation similar to (\ref{dnu1-n-a}), namely:
\begin{eqnarray}
{\cal S}^{-}_{i_{1},...,i_{n}} {\cal S}^{+}_{I_{1}p_{1}} \cdots 
{\cal S}^{+}_{I_{s}p_{s}} 
\left[ \tilde{\Phi'}^{\{j_{1},...,j_{k}\},I_{1},...,I_{s}}_{\sigma,\nu_{1},...,
\nu_{s},i_{s+1},...,i_{n}} + (-1)^{s} {\cal S}_{j_{1},...,j_{k}}
\delta^{j_{1}}_{i_{s+1}} \Phi^{\{j_{2},...,j_{k}\},I_{1},...,I_{s}}_{\sigma,
\nu_{1},...,\nu_{s},i_{s+2},...,i_{n}} \right] \times \nonumber \\
\delta^{p_{1}}_{i_{1}} \cdots \delta^{p_{s}}_{i_{s}} = 0, 
\quad (s = 1,...,n, k = 0,...,r)
\label{dnu1-n-1-a}
\end{eqnarray}

In compact notations this amounts to
\be
B_{1} \cdots B_{s} \left((\Delta' \Phi)^{k}_{s} + B_{0} \Phi^{k-1}_{s}\right) 
= 0, \quad (s = 1,...,n, k = 0,...,r)
\label{descent-n-1-a}
\ee
which can be analysed with the, by now familiar, descent technique. Some terms
can be grouped into the exterior differential of a contact form and in this way
the formula from the statement emerges.
$\qed$

The converse of this proposition is not true. We proceed nevertheless to the 
study of locally variational differential equations by proving an analogue of 
proposition \ref{e-n}, namely:

\begin{prop}
Let $T$ be a differential equation. Then $T$ is locally variational {\it iff} 
\be
E_{n+1}([T]) = 0.
\ee
\end{prop}

{\bf Proof:}
It is an immediate consequence of the definition of a variationally trivial
Lagrange form and of theorem \ref{En+1}.
$\qed$

Then we have the central result of this section:

\begin{thm}
Let
$T$
be a locally variational differential equation. Then $T$ has the following
local form in $V^{r}$ given by (\ref{T}) where:
\be
T_{\sigma} = \varepsilon^{i_{1},...,i_{n}} \sum_{s=0}^{n} 
\sum_{|I_{1}|,...,|I_{s}|=r-1} 
{\cal T}^{I_{1},...,I_{s}}_{\sigma,\nu_{1},...,\nu_{s},i_{s+1},...,i_{n}}
{\cal J}_{I_{1},...,I_{s}}^{\sigma_{1},...,\sigma_{s},i_{s+1},...,i_{n}}
\label{TJ}
\ee
where 
${\cal T}^{I_{1},...,I_{s}}_{\sigma,\nu_{1},...,\nu_{s},i_{s+1},...,i_{n}}$
are some smooth functions on 
$V^{r-1}$
verifying symmetry properties of the type (\ref{symmetry}); this symmetry
property leaves aside the index $\sigma$. 
\end{thm}

{\bf Proof:}
(i) According to the previous proposition and applying the exactness of the 
variational sequence (theorem \ref{sequence}), there exists 
$\lambda \in \Omega^{r}_{n}$
such that
$$
[T] = [d\lambda],
$$
or, equivalently
$$
T -d\lambda \in \theta^{r}_{n+1}.
$$

If we use (\ref{theta-q}) we obtain that
$$
T -d\lambda = \alpha + d\beta
$$
where
$\alpha \in \Omega^{r}_{n+1(c)}$
and
$\beta \in \Omega^{r}_{n(c)}.$
We can redefine
$\lambda \rightarrow \lambda - \beta$
and we obtain
\be
T -d\lambda = \alpha \in \Omega^{r}_{n+1(c)}
\ee
or, equivalently:
\be
p_{1} T = p_{1} d\lambda. 
\label{p1-T}
\ee

(ii) In particular, it follows that the $(n+1)$-form $p_{1}d\lambda$ must be
$\pi^{r+1,r}$-projectable. We can apply proposition \ref{eta-nu-n} and 
obtain that $\lambda$ has the following expression:
\be
\lambda = \nu + \sum_{|J| \leq r-1} \Phi^{J}_{\sigma,i_{1},...,i_{n-1}} 
\omega^{\sigma}_{J} \wedge dx^{i_{1}} \wedge \cdots \wedge dx^{i_{n-1}} + 
d\phi + \psi
\label{lambda-T}
\ee
where
$\nu \in \Omega^{r}_{n,X}$ 
is a basic $n$-form,
$\psi \in \Omega^{r}_{n+1}$
is a 2-contact form and
$
\Phi^{J}_{\sigma,i_{1},...,i_{n-1}} \quad |J| = 1,....,r-1
$
are smooth functions on $V^{r}$ completely antisymmetric in the indices
$i_{1},...,i_{n-1};$
for uniformity of notations we will extend the sum to $|J| = r$ with the 
convention:
\be
\Phi^{J}_{\sigma,i_{1},...,i_{n-1}} = 0, \quad |J| = r.
\ee

One computes the exterior differential of this form and imposes the condition
of projectability. As in the proposition \ref{nu-nu-tilde} one finds out that
$p_{1}d\lambda$ is $\pi^{r+1,r}$-projectable {\it iff} the following equations 
are fulfilled:
\be
{\cal S}^{-}_{i_{1},...,i_{n}} \delta^{p}_{i_{1}}
\tilde{\Phi}^{J,I_{1}}_{\sigma,\nu_{1},i_{2},...,i_{n}} = 0
\label{proj-n}
\ee
or :
\be
{\cal S}^{-}_{i_{1},...,i_{n}} {\cal S}_{Ip} \delta^{p}_{i_{1}}
\left( \partial^{I}_{\sigma} \Phi^{J}_{\sigma,i_{2},...,i_{n}}\right) = 0.
\label{proj-n-a}
\ee

This equation is of the same type as (\ref{nu2}) and the analysis performed 
there can be applied (one notices that the index $q$ can take the value 
$n-1$). As a result, 
$\Phi^{J}_{\sigma}$
has a polynomial of the type (\ref{AAA})
\be
\Phi^{J}_{\sigma,i_{1},...,i_{n-1}}  = \sum_{s=0}^{n-1} 
\sum_{|I_{1}|,...,|I_{s}|=r-1} {\cal S}^{-}_{i_{1},...,i_{n-1}}
C^{J,I_{1},...,I_{s}}_{\sigma,\nu_{1},...,\nu_{s},i_{1},...,i_{n-1}}
y^{\nu_{1}}_{I_{1}i_{1}} \cdots y^{\nu_{s}}_{I_{s}i_{s}}
\label{phi-C}
\ee
where one can suppose that the tensors
$C^{J,I_{1},...,I_{s}}_{\sigma,\nu_{1},...,\nu_{s},i_{1},...,i_{n-1}}$
have (partial) symmetry properties of the type (\ref{sym-partial}). 

(iii) Now one inserts the expressions (\ref{T}) and (\ref{lambda-T}) for $T$ 
and respectively for $\lambda$ into the projectability condition (\ref{p1-T}) 
and obtains after some algebra the following two relations:
\be
T_{\sigma} = \partial_{\sigma} L + \varepsilon^{i_{1},...,i_{n}} 
\tilde{\Phi}^{\emptyset}_{\sigma,i_{1},...,i_{n}}
\label{T-sigma}
\ee
and
\be
\partial^{J}_{\sigma} L = \varepsilon^{i_{1},...,i_{n}}  \left[
\sum_{s=0}^{n-1} \sum_{|I_{1}|,...,|I_{s}|=r-1} 
(B_{0}C)^{J,I_{1},...,I_{s}}_{\sigma,\nu_{1},...,\nu_{s},i_{s+1},...,i_{n}}
y^{\nu_{1}}_{I_{1}i_{1}} \cdots y^{\nu_{s}}_{I_{s}i_{s}} +
\tilde{\Phi}^{J}_{\sigma,i_{1},...,i_{n}} \right], \quad (|J| = 1,...,r).
\label{derivatives-L}
\ee

One must substitute into the last equation the expressions (\ref{phi-C}) for 
$
\tilde{\Phi}^{J}_{\sigma,i_{1},...,i_{n}}
$ and after some computations one arrives at 
\begin{eqnarray}
\partial^{J}_{\sigma} L = F^{J}_{\sigma} \equiv 
\varepsilon^{i_{1},...,i_{n}} 
\sum_{s=0}^{n-1} \sum_{|I_{1}|,...,|I_{s}|=r-1} \left[
(B_{0}C)^{J,I_{1},...,I_{s}}_{\sigma,\nu_{1},...,\nu_{s},i_{s+1},...,i_{n}} +
(\Delta'C)^{J,I_{1},...,I_{s}}_{\sigma,\nu_{1},...,\nu_{s},i_{s+1},...,i_{n}}
\right] \times \nonumber \\
y^{\nu_{1}}_{I_{1}i_{1}} \cdots y^{\nu_{s}}_{I_{s}i_{s}} , 
\quad (|J| = 1,...,r).
\label{der-L-F}
\end{eqnarray}

It is clear that Frobenius conditions of integrability should be fulfilled:
\be
\partial^{K}_{\zeta} F^{J}_{\sigma} = \partial^{J}_{\sigma} F^{K}_{\zeta},
\quad |J|, |K| = 1,...,r.
\label{Frobenius}
\ee

Instead of trying to solve directly these equations we proceed as follows. We
admit that the integrability conditions (\ref{Frobenius}) are fulfilled for
$|J| = |K| = r.$ Then one can obtain from (\ref{der-L-F}) with a well-known 
homotopy formula:
\be
L = L_{0} + \int_{0}^{1} \sum_{|J|=r} y^{\sigma}_{J} 
F^{J}_{\sigma}(x^{i},y^{\sigma},y^{\sigma}_{j},...,
y^{\sigma}_{j_{1},...,j_{r-1}},...,u y^{\sigma}_{j_{1},...,j_{r}}) du
\label{LL}
\ee
where $L_{0}$ does not depend on the highest order derivatives. Using the
explicit expression for $F^{J}_{\sigma}$ from (\ref{der-L-F}) one arrives at:
\be
L = L_{0} + \varepsilon^{i_{1},...,i_{n}} \sum_{s=1}^{n} {1 \over s}
\sum_{|I_{1}|,...,|I_{s}|=r-1} 
A^{I_{1},...,I_{s}}_{\nu_{1},...,\nu_{s},i_{s+1},...,i_{n}}
y^{\nu_{1}}_{I_{1}i_{1}} \cdots y^{\nu_{s}}_{I_{s}i_{s}}
\label{L-T}
\ee
where
\be
A^{I_{0},...,I_{s}}_{\nu_{0},...,\nu_{s},i_{s+1},...,i_{n-1}} \equiv
{\cal S}^{-}_{(I_{0},\nu_{0}),...,(I_{s},\nu_{s})} 
C^{I_{0},...,I_{s}}_{\nu_{0},...,\nu_{s},i_{s+1},...,i_{n-1}}, \quad 
(s = 0,...,n-1)
\ee
is constructed from $C^{...}_{...}$ such that it verifies the complete symmetry
property (\ref{symmetry}). Now one computes from (\ref{LL}) the partial
derivatives of $L$ and, considering the highest order ones, obtains by 
comparison with the relation (\ref{der-L-F}) 
$$
B_{0} \cdots B_{s} (C_{s} - A_{s}) = 0, \quad s = 0,...,n-1).
$$

If we apply lemma \ref{BXB} one obtains, as usual, that
$$
C^{I_{1},...,I_{s}}_{\nu_{1},...,\nu_{s},i_{s+1},...,i_{n}} =
A^{I_{1},...,I_{s}}_{\nu_{1},...,\nu_{s},i_{s+1},...,i_{n}} +
\delta - {\rm terms}.
$$

If we substitute this expression for $C^{...}_{...}$ into the equation 
(\ref{phi-C}) for $|J| = r-1$ we obtain that the delta terms give a null 
contribution; as a consequence, it follows that the expressions
$
C^{I_{1},...,I_{s}}_{\nu_{1},...,\nu_{s},i_{s+1},...,i_{n}}, \quad
|I_{1}|,...,|I_{s}| = r-1
$
can be considered to have the symmetry property (\ref{symmetry}) without
loosing the generality. Then the integrability condition (\ref{Frobenius}) is
fulfilled for $|J| = r.$ 

(iv) Finally, we develop the expression (\ref{T-sigma}) for $T_{\sigma}$ using
the expression (\ref{phi-C}) for 
$
\tilde{\Phi}^{\emptyset}_{\sigma,i_{1},...,i_{n}}
$
and the expression (\ref{LL}) for $L$. It is elementary to prove that both
terms in (\ref{T-sigma}) are polynomials in the hyper-Jacobians of the type
appearing in the statement of the theorem. 
$\qed$

\begin{rem}
The functions 
$
{\cal T}^{I_{1},...,I_{s}}_{\sigma,\nu_{1},...,\nu_{s},i_{s+1},...,i_{n}}
$
appearing in the statement cannot be completely arbitrary because we did not
use all the integrability conditions (\ref{Frobenius}). 
\end{rem}

\begin{rem}
For the case $r = 2$ analysed in detail in \cite{Gr3} it is possible to use
completely the Helmholtz-Sonin equations involving the highest order
derivatives. One obtains that $T_{\sigma}$ has a expression as in the statement
of the theorem, but the coefficients
$
{\cal T}^{l_{1},...,l_{s}}_{\sigma,\nu_{1},...,\nu_{s},i_{s+1},...,i_{n}}
$
can be chosen to be completely symmetric in the indices
$
\sigma,\nu_{1},...,\nu_{s}
$
and completely antisymmetric in the indices
$
l_{1},....,l_{s}.$
\end{rem}

\section{Conclusions}

We have succeeded to give a fairly complete presentation of the most important
results connected with the existence of the variational sequences. Many of the
proofs have been essentially simplified using Fock space techniques. We also
have been able to give very explicit expression for the most general form of a
variationally trivial Lagrangian and the generic expression of a locally
variational system of partial differential equations, completing in this way
known results from the literature.

Some criticism of the approach to the Lagrangian formalism accepted in this
paper is necessary. We think that the weak point of this approach, from the
physical point of view, is in fact, the definition of a differential equation. 
The definition of such an object as a special type of differential form leads
at the following problem. Let
$T_{\sigma}, \quad T_{\sigma}'$
be the components of two differential equations in a given chart. It is 
possible (and examples can be provided) that one can arrange such that:
(1) the hyper-surface 
$T_{\sigma} = 0 \quad \sigma = 1,...,m$
coincides with the hyper-surface
$T_{\sigma}'=  0 \quad \sigma = 1,...,m$
(in this way the two sets of functions are describing in fact the
same set of physical solutions of the equation of motion) and (2) 
$T_{\sigma}$
are locally variational and
$T_{\sigma}'$
are not locally variational. So, in a certain sense, the property of being
locally variational is not intrinsically defined. (See on this point also 
\cite{AT}). In this sense, the most reasonable definition (from the physical
point of view) of a differential equation would be certain hyper-surfaces in
the jet bundle extension with some regularity properties (guaranteeing a well
posed Cauchy problem). The problem would be to attach in an intrinsic way the
property of being locally variational to such a hyper-surface. 

Another interesting and open problem is to find, if possible, a physical
meaning for the elements of the variational sequence of index $q \geq n+2$
and eventually some representatives by forms.

One would also be interested to see to what extent the results of this 
approach to the Lagrangian formalism can be extended to the case when $Y$ is
not a fibre bundle over some ``space-time" manifold $X$ (the typical case being
a relativistic system with $Y$ the Minkowski space). Although it is clear that
the line of argument from this paper depends essentially on the existence of
the fibre bundle structure, some steps in this direction exists in the
literature \cite{Gr2}, \cite{GK}; the proper substitute for jet extensions of
a fibre bundle are the higher order Grassmann bundles.

Finally, there exists some physical interest to extend this formalism to the
situation when anticommuting variables are present. This case appears when one
is studying, for instance, BRST-type symmetries.
\newpage

\section{Appendix}

In this appendix we give the basic definitions of the Fock space concepts we
have used in this paper and provide a fairly simple proof of the so-called 
{\it trace decomposition theorem} \cite{Kr7}. One can simplify somehow all the
proofs in this paper if one uses this more refined decomposition of tensors;
however this simplification is rather modest.

\subsection{Fock Space Notions}

To avoid unnecessary complications we consider
${\cal H}$
to be a finite dimensional real Hilbert space. Then the {\it associated Fock
space} is, by definition: 
$$ 
{\cal F}({\cal H}) \equiv \oplus_{n=0}^{\infty} {\cal H}_{n} 
$$ 
where 
$$ 
{\cal H}_{n} \equiv {\cal H}^{\otimes n}, \quad n > 0, \quad
{\cal H}_{0} \equiv \R.
$$ 

The Hilbert space 
${\cal H} \sim {\cal H}_{1}$
is called the {\it one-particle space} and the element
$(1,0,\dots) \in {\cal F}({\cal H})$
is called the {\it vacuum}.

We can introduce the symmetrization and the antisymmetrization operators in
${\cal H}_{n}$
according to
\be
{\cal S}^{\pm}_{n} \phi_{1} \otimes \dots \otimes \phi_{n} \equiv
{1 \over n!} \sum_{P \in {\cal P}_{n}} \epsilon_{\pm}(P) 
\phi_{P(1)} \otimes \dots \otimes \phi_{P(n)}
\label{sym-antisym}
\ee 
where
\be
\epsilon_{+}(P) \equiv 1, \quad \epsilon_{-}(P) \equiv (-1)^{|P|}, \quad
\forall P \in {\cal P}_{n}.
\ee
Here 
${\cal P}_{n}$
is the permutation group of the numbers
$1,2,\dots,n$
and 
$|P|$
is the signature of the permutation $P$. One can prove easily that these
operators are in fact orthogonal projectors. We also define the following 
projector operators acting in the Fock space:
\be
{\cal S}^{\pm} \equiv \oplus_{n=0}^{\infty} {\cal S}^{\pm}_{n}.
\ee

Next, one defines the following subspaces of
${\cal H}_{n}$:
\be
{\cal H}_{n}^{\pm} \equiv {\cal S}^{\pm}_{n} {\cal H}_{n}
\ee
and of 
${\cal F}({\cal H})$:
\be
{\cal F}^{\pm}({\cal H}) \equiv {\cal S}^{\pm} {\cal F}({\cal H}) =
\oplus_{n=0}^{\infty} {\cal S}^{\pm}_{n} {\cal H}_{n}.
\ee
The subspace
${\cal F}^{\pm}({\cal H})$
is called the {\it bosonic (or symmetric) Fock space} for the sign $+$ and the 
{\it fermionic (or antisymmetric) Fock space} for the sign $-$.

One can define in
${\cal F}({\cal H})$
the so-called {\it particle number operators} according to
\be
N \phi_{1} \otimes \dots \otimes \phi_{n} \equiv 
n \phi_{1} \otimes \dots \otimes \phi_{n}.
\label{number}
\ee

It is clear that these operators can be restricted to the bosonic and to the
fermionic Fock spaces.

To simplify the presentation according to our specific needs, we consider an
orthonormal basis in the 
${\cal H}: e_{1},\dots,e_{k}, \quad k \equiv {\rm dim}({\cal H})$.
Then every element of the Fock space
${\cal F}({\cal H})$
can be represented as a collection
$$
(f,f^{i},\dots,f^{i_{1},\dots,i_{n}},\dots)
$$
where 
$f^{i_{1},\dots,i_{n}}$
are the elements of a (real) tensor of degree $n$.

The operators 
${\cal S}^{\pm}$
are represented by the following formul\ae:
\be
\left( {\cal S}^{\pm} f\right)^{i_{1},\dots,i_{n}} \equiv
{1 \over n!} \sum_{P \in {\cal P}_{n}} \epsilon_{\pm}(P)
f^{i_{P(1)},\dots,i_{P(n)}}, \quad \forall n > 0.
\label{sa}
\ee

Sometimes it is convenient to indicate explicitly the indices affected by the
operation of symmetrization, or antisymmetrization, by writing the preceding
formula as follows:
\be
{\cal S}^{\pm}_{i_{1},\dots,i_{n}} f^{i_{1},\dots,i_{n}} \equiv
{1 \over n!} \sum_{P \in {\cal P}_{n}} \epsilon_{\pm}(P)
f^{i_{P(1)},\dots,i_{P(n)}}, \quad \forall n > 0.
\label{sa-indices}
\ee

This notation is important because we have the following formul\ae~ which are
used many times in the proofs:
\be
{\cal S}^{\pm}_{I} {\cal S}^{\pm}_{J} = {\cal S}^{\pm}_{I}, \quad \forall 
J \subset I.
\ee

The bosonic (resp. fermionic) Fock space is formed by symmetric (resp. 
antisymmetric) tensors.

We are ready to introduce now the so-called {\it creation} and {\it
annihilation} operators. We will not use the usual convention form quantum
mechanics textbooks, but another one more suitable for our purposes.
We also use consistently Bourbaki conventions:
$$
\sum_{\emptyset} \equiv 0, \quad \prod_{\emptyset} \equiv 1.
$$

Consider first the bosonic case. We define in the bosonic Fock space the
{\it creation} and {\it annihilation} operators
$b^{*}_{l}, \quad b^{l}, \quad l = 1,\dots,n$
by: 
\be
\left( b^{*}_{l} f\right)^{i_{1},\dots,i_{n}} \equiv
{1 \over n} \sum_{p=1}^{n} \delta_{l}^{i_{p}} 
f^{i_{1},\dots,\hat{i}_{p},\dots,i_{n}} = 
{\cal S}^{+}_{i_{1},\dots,i_{n}} \delta_{l}^{i_{1}} f^{i_{2},\dots,i_{n}} , 
\quad l = 1,\dots,k, \quad \forall n \geq 0
\label{creation-b}
\ee
(where we prefer to specify explicitly the indices on which the operation of
symmetrization is performed; evidently, for 
$n = 0$
the right hand side must be considered $0$) and respectively by:
\be
\left( b^{l} f\right)^{i_{1},\dots,i_{n}} \equiv (n+1) f^{l,i_{1},\dots,i_{n}},
\quad \forall n \geq 0.
\label{annihilation-b}
\ee

It is easy to see that, in fact,
$b^{*}_{l}$
is the adjoint of
$b^{l}$.
The essential property of these operators is contained in the so-called
{\it canonical commutation relations (CCR)}:
\be
[b^{l}, b^{m}] = 0, \quad [b^{*}_{l}, b^{*}_{m}] = 0, \quad 
[b^{l}, b^{*}_{m}] = \delta^{l}_{m}, \quad \forall l,m = 1,\dots,k
\label{CCR}
\ee
where by
$[\cdot,\cdot ]$
we are designating, obviously, the commutator of the two entries.

We now consider the fermionic case. The definitions are in complete analogy
with the previous ones. We define in the fermionic Fock space the
{\it creation} and {\it annihilation} operators
$a^{*}_{l}, \quad a^{l}, \quad l = 1,\dots,n$
by: 
\be
\left( a^{*}_{l} f\right)^{i_{1},\dots,i_{n}} \equiv
{1 \over n} \sum_{p=1}^{n} (-1)^{p-1} \delta_{l}^{i_{p}}
f^{i_{1},\dots,\hat{i}_{p},\dots,i_{n}} = 
{\cal S}^{-}_{i_{1},\dots,i_{n}} \delta_{l}^{i_{1}} f^{i_{2},\dots,i_{n}} , 
\quad l = 1,\dots,k, \quad \forall n \geq 0
\label{creation-f}
\ee
(where again we specify explicitly the indices on which the operation of
antisymmetrization is performed; for 
$n = 0$
the right hand side must be considered $0$) and respectively by:
\be
\left( a^{l} f\right)^{i_{1},\dots,i_{n}} \equiv (n+1) f^{i_{1},\dots,i_{n},l},
\quad \forall n \geq 0.
\label{annihilation-f}
\ee

Again one sees that
$a^{*}_{l}$
is the adjoint of
$a^{l}$.
We have analogously to (\ref{CCR}) the {\it canonical anticommutation 
relations (CAR)}:
\be
\{a^{l}, a^{m}\} = 0, \quad \{a^{*}_{l}, a^{*}_{m}\} = 0, \quad 
\{a^{l}, a^{*}_{m}\} = \delta^{l}_{m}, \quad \forall l,m = 1,\dots,k
\label{CAR}
\ee
where by
$\{\cdot,\cdot \}$
we are designating the anticommutator of the two entries.

\begin{rem}
It is not very complicated to find the abstract definitions for the creation
and the annihilation operators, i.e. basis independent definitions, but we will
not need them here. It is also noteworthy that the whole formalism works for
any Hilbert space, even infinite dimensional, defined over an arbitrary
commutative division field.
\end{rem}

It is obvious how to extend the definitions of the particle number operator and
of the creation and annihilation operators to the more general case of the
Hilbert space 
${\cal H}_{s}$
appearing in the proof of theorem \ref{contact}. In particular, we have more
particle number operators:
$N_{f}$
corresponding to the fermionic degrees of freedom and
$N_{\alpha}, \quad \alpha = 1,\dots,s$
corresponding to the bosonic degrees of freedom and and giving $k$
(respectively 
$r_{\alpha}$)
when applied 
on
${\cal H}_{k,r_{1},...,r_{s}}$.

In this case, if we define the operators
$B_{\alpha}, \quad \alpha = 1,\dots,s$
according to (\ref{B-alpha}) then it is easy to establish the following 
anticommutation relations:
\be
\{B_{\alpha}, B_{\alpha}^{*}\} = N_{f} - N_{\alpha} + n{\bf 1}, \quad \alpha =
1,\dots,s
\label{BB-star}
\ee
and
\be
\{ B_{\alpha},B_{\beta}\} = 0, \quad \alpha, \beta = 1,\dots,s;
\label{BB}
\ee
in particular we have:
\be
B_{\alpha}^{2} = 0, \quad \alpha = 1,\dots,s.
\label{B2}
\ee

We also have the commutation relations:
\be
[B_{\alpha}, N_{f} - N_{\alpha} + n{\bf 1}] = 0, \quad \alpha = 1,\dots,s.
\label{BN}
\ee

\newpage

\subsection{The Trace Decomposition Formula}

A tensor of the form
$A^{I_{1},...,I_{s}}_{\sigma_{1},...,\sigma_{s},i_{s+1},...,i_{q}}$
is called {\it traceless iff} 
it verifies:
\be
A^{\{j,p_{1},...,p_{l}\},I_{2},...,I_{s}}_{\sigma_{1},...,\sigma_{s},j,
i_{s+1},...,i_{q}} = 0
\label{traceless}
\ee
for all indices left free. If we use compact tensor notations, this can be
written as follows:
\be
B^{*}_{1} A_{\sigma_{1},...,\sigma_{s}} = 0.
\ee

Because of the symmetry properties (\ref{symmetry}) we have in fact:
\be
B^{*}_{\alpha} A_{\sigma_{1},...,\sigma_{s}} = 0 \quad (\alpha = 1,...,s).
\ee

Then we have the following result

\begin{lemma}
Let
$X \in {\cal H}_{p,r_{1},...,r_{s}}\quad (r_{1},...,r_{s} \in \N^{*}, 
\quad s \leq n-p)$.
Then the following decomposition formula is valid:
\be
X = X_{0} + \sum_{\alpha = 1}^{k} B_{\alpha} X_{\alpha}
\label{dec}
\ee
where the tensor 
$X_{0}$
is traceless and it is {\bf uniquely} determined by $X$.
\end{lemma}

{\bf Proof:}

It follows elementary from the relations at the end of the preceding 
subsection that the operators
\be
P_{\alpha} \equiv (N_{f} - N_{\alpha} + n{\bf 1})^{-1} 
B_{\alpha} B_{\alpha}^{*}, \quad
Q_{\alpha} \equiv (N_{f} - N_{\alpha} + n{\bf 1})^{-1} 
B_{\alpha}^{*} B_{\alpha}, \quad \alpha = 1,\dots,s
\ee
are orthogonal projectors and we have
\be
P_{\alpha} + Q_{\alpha} = 1, \quad \alpha = 1,\dots,s.
\ee

We consider now the corresponding subspaces
\be
V_{\alpha} \equiv {\rm Span}(P_{\alpha}{\cal H}_{s}), \quad \alpha = 1,\dots,s.
\ee

It is well known that the set of all linear subspaces in a finite dimensional
Hilbert space is a orthomodular lattice if the order relation 
$<$ 
is given by the inclusion of subspaces 
$\subset$
and the orthogonalization operator is constructed using the scalar product 
\cite{Va}. It follows that one has the following formul\ae~ for the 
{\it infimum} and {\it supremum} operations:
\be
\wedge_{i \in I} V_{i} \equiv \cap_{i \in I} V_{i}, \quad
\vee_{i \in I} V_{i} \equiv \sum_{i \in I} V_{i}.
\ee

The property above can be transported for the set of orthogonal projectors
acting in this Hilbert space. In this new representation, the orthogonalization
operation looks very simple:
\be
P^{\perp} \equiv {\bf 1} - P.
\ee

Now we decompose the generic element $X$ as follows:
\be
X = X_{0} + X_{1}
\label{trace-dec}
\ee
where
\be
X_{0} \equiv \left( \wedge Q_{\alpha} \right) X,\quad X_{1} \equiv X - X_{0}.
\ee

One immediately establishes that 
$X_{0}$
is traceless and that
$$
X_{1} = \left( \wedge Q_{\alpha} \right)^{\perp} X = \vee Q_{\alpha}^{\perp} X 
= \vee P_{\alpha} X
$$
so
$X_{1}$
has the form:
$$
X_{1} = \sum_{\alpha} B_{\alpha} X_{\alpha}.
$$

This finishes the proof.
$\qed$

\begin{rem}
One can prove that the lattice property is valid in a more general case of
linear spaces of Hilbertian type \cite{Va}).
\end{rem}

\begin{rem}
A more general trace decomposition formula appears in \cite{Kr7} in the sense
that one can give up the hypothesis that the Hilbert space is of Fock type.
\end{rem}

As a consequence we have

\begin{thm}
Let
$\rho \in \Omega^{r}_{q}$
and let
$(V,\psi)$
be a chart on $Y$. Then the form $\rho$ admits in the chart
$(V^{r},\psi^{r})$
the following {\bf unique} decomposition:
\be
\rho = \rho_{0} + \rho'
\ee
where 
$\rho_{0} \in \Omega^{r}_{q(c)}$
is a contact form, and $\rho'$ has the expression
\be
\rho' = \sum_{s=0}^{q} {1 \over s!(q-s)!} \sum_{|I_{1}|,...,|I_{s}|=r}
A^{I_{1},...,I_{s}}_{\sigma_{1},...,\sigma_{s},i_{s+1},...,i_{q}}
dy^{\sigma_{1}}_{I_{1}} \wedge \cdots \wedge dy^{\sigma_{s}}_{I_{s}} \wedge
dx^{i_{s+1}} \wedge \cdots \wedge dx^{i_{q}},
\ee
where the tensors
$A^{I_{1},...,I_{s}}_{\sigma_{1},...,\sigma_{s},i_{s+1},...,i_{q}}$
verify the antisymmetry property (\ref{symmetry}) are traceless. 
\end{thm}

{\bf Proof:}

We use the standard decomposition given by (\ref{ro}) - (\ref{canonical}) and
we refine it. 

We apply the preceding lemma to the tensors
$A^{I_{1},...,I_{s}}_{\sigma_{1},...,\sigma_{s},i_{s+1},...,i_{q}}$
appearing in the expression of $\rho'$ (see (\ref{canonical})).
So, this tensor splits into a $\delta$-contribution and a traceless
contribution. 
The first contribution leads to terms containing at least a factor
$d\omega^{\sigma}_{J} \quad |J| = r-1)$
so it is a contact form and can be combined with
$\rho_{0}$;
the second contribution is traceless. This means that we have the desired 
decomposition. The uniqueness assertion is easy to prove.
$\qed$

\begin{rem}
It immediately follows that the equivalence classes
$
\Omega^{r}_{q}(V) \left/ \right. \Omega^{r}_{q(c)}(V) \quad (q = 1,...,n)
$
(for any chart $V$ on $Y$) are indexed by sets of traceless tensors
$$
A^{I_{1},...,I_{s}}_{\sigma_{1},...,\sigma_{s},i_{s+1},...,i_{q}} 
\quad (|I_{1}| = \cdots = |I_{s}| = r, \quad s = 0,...,q).
$$

For the proof of the isomorphism
$
\Omega^{r}_{q}(V) \left/ \right. \Omega^{r}_{q(c)}(V) \cong 
{\cal J}^{r+1}_{q}(V) \quad (q = 1,...,n)
$
see the end of \cite{Kr6}.
\end{rem}
\newpage


\end{document}